%
\newif\ifloadreferences\loadreferencestrue
%
%
%
%
%
\let\myfrac=\frac%
\input eplain %
\let\frac=\myfrac%
\let\myfootnote=\footnote%
\input amstex \input epsf %
\let\footnote=\myfootnote%
%
%
\loadeufm\loadmsam\loadmsbm\message{symbol names}\UseAMSsymbols\message{,}%
%
\font\myfontdefault=cmr10%
\newif\ifmakebiblio%
\newif\ifinappendices%
\newif\ifundefinedreferences%
\newif\ifchangedreferences%
\makebibliofalse%
\undefinedreferencesfalse%
\changedreferencesfalse%
%
%
%
%
%
\def\setcatcodes{\catcode`\!=0 \catcode`\\=11}%
{\global\let\noe=\noexpand%
\catcode`\@=11 \catcode`\_=11 \setcatcodes%
!global!def!_@@internal@@makeref#1{%
!global!expandafter!def!csname #1ref!endcsname##1{%
!csname _@#1@##1!endcsname%
!expandafter!ifx!csname _@#1@##1!endcsname!relax%
    !write16{#1 ##1 not defined - run saving references}%
    !undefinedreferencestrue%
!fi}}%
!global!def!_@@internal@@makelabel#1{%
!global!expandafter!def!csname #1label!endcsname##1{%
!edef!temptoken{!csname #1info!endcsname}%
!ifloadreferences%
!expandafter!ifx!csname _@#1@##1!endcsname!relax%
!write16{#1 ##1 not hitherto defined - rerun saving references}%
!changedreferencestrue%
!else%
!expandafter!ifx!csname _@#1@##1!endcsname!temptoken%
!else%
!write16{#1 ##1 reference has changed - rerun saving references}%
!changedreferencestrue%
!fi%
!fi%
!else%
!expandafter!edef!csname _@#1@##1!endcsname{!temptoken}%
!edef!textoutput{!write!references{\global\def\_@#1@##1{!temptoken}}}%
!textoutput%
!fi}}%
!global!def!makecounter#1{!_@@internal@@makelabel{#1}!_@@internal@@makeref{#1}}%
!unsetcatcodes%
}
%
%
%
%
%
\def\turnintolatin#1{\ifcase #1 _\or i\or ii\or iii\or iv\or v\or vi\or vii\or viii\or ix\or x\or xi\or xii\or xiii\or xiv\or xv\or xvi\or xvii\or xviii\or xix\or xx\or xxi\or xxii\or xxiii\or xxiv\or xxv\or xxvi\fi}%
\def\alphanum#1{\ifcase #1 _\or A\or B\or C\or D\or E\or F\or G\or H\or I\or J\or K\or L\or M\or N\or O\or P\or Q\or R\or S\or T\or U\or V\or W\or X\or Y\or Z\fi}%
\newwrite\references%
\ifloadreferences{\catcode`\@=11 \catcode`\_=11 \global\def\_@citation@AltschulerWu{1}
\global\def\_@citation@ClutterbuckSchnurerSchulze{2}
\global\def\_@citation@ChengZhou{3}
\global\def\_@citation@ChengMejiaZhouI{4}
\global\def\_@citation@ChengMejiaZhouII{5}
\global\def\_@citation@ChengMejiaZhouIII{6}
\global\def\_@citation@DavilaDelPinoNguyen{7}
\global\def\_@citation@GilbTrud{8}
\global\def\_@citation@Halldorsson{9}
\global\def\_@citation@HauswirthPacard{10}
\global\def\_@citation@HoffmanMeeks{11}
\global\def\_@citation@KapouleasI{12}
\global\def\_@citation@KapouleasII{13}
\global\def\_@citation@KapouleasIII{14}
\global\def\_@citation@MartinSavasSmoczyk{15}
\global\def\_@citation@MazzeoPacard{16}
\global\def\_@citation@MazzeoPollack{17}
\global\def\_@citation@Morabito{18}
\global\def\_@citation@Nayatani{19}
\global\def\_@citation@NguyenI{20}
\global\def\_@citation@NguyenII{21}
\global\def\_@citation@Pacard{22}
\global\def\_@citation@Schwarz{23}
\global\def\_@citation@Wang{24}
\global\def\_@citation@Weber{25}
\global\def\_@head@Introduction{1}
\global\def\_@subhead@BackgroundAndMainResult{1.1}
\global\def\_@eqn@IntroMCFSEquation{\relax \unhbox \voidb@x \hbox {{\relax \tenrm (1.1)}}}
\global\def\_@subhead@Techniques{1.2}
\global\def\_@subhead@Notation{1.3}
\global\def\_@eqn@AppendixControl{\relax \unhbox \voidb@x \hbox {{\relax \tenrm (1.2)}}}
\global\def\_@eqn@AppendixAsymptoticControl{\relax \unhbox \voidb@x \hbox {{\relax \tenrm (1.3)}}}
\global\def\_@subhead@Acknowledgements{1.4}
\global\def\_@head@GrimSurfaces{2}
\global\def\_@subhead@TheLargeScale{2.1}
\global\def\_@eqn@PrincipleCurvaturesOfSurfaceOfRevolution{\relax \unhbox \voidb@x \hbox {{\relax \tenrm (2.1)}}}
\global\def\_@eqn@NormalComponentOfSurfaceOfRevolution{\relax \unhbox \voidb@x \hbox {{\relax \tenrm (2.2)}}}
\global\def\_@eqn@IntegratedFormulaForGrimSurface{\relax \unhbox \voidb@x \hbox {{\relax \tenrm (2.3)}}}
\global\def\_@eqn@AsymptoticSolutionAtInfinity{\relax \unhbox \voidb@x \hbox {{\relax \tenrm (2.4)}}}
\global\def\_@proc@AsymptoticSolutionAtInfinity{2.1.1}
\global\def\_@eqn@GrimSurfaceEquation{\relax \unhbox \voidb@x \hbox {{\relax \tenrm (2.5)}}}
\global\def\_@eqn@AnsatzLaurentPolynomial{\relax \unhbox \voidb@x \hbox {{\relax \tenrm (2.6)}}}
\global\def\_@proc@LemmaFormalSolution{2.1.2}
\global\def\_@eqn@RecurrenceRelationAtInfinity{\relax \unhbox \voidb@x \hbox {{\relax \tenrm (2.7)}}}
\global\def\_@eqn@PartialSumAtInfinity{\relax \unhbox \voidb@x \hbox {{\relax \tenrm (2.8)}}}
\global\def\_@eqn@FirstOrderEstimateLargeScale{\relax \unhbox \voidb@x \hbox {{\relax \tenrm (2.9)}}}
\global\def\_@proc@LemmaFirstOrderEstimate{2.1.3}
\global\def\_@eqn@ZeroethOrderAsymptoticRelationLargeScale{\relax \unhbox \voidb@x \hbox {{\relax \tenrm (2.10)}}}
\global\def\_@proc@ZeroethOrderAsymptoticRelationLargeScale{2.1.4}
\global\def\_@eqn@BoundOnDerivativeOfW{\relax \unhbox \voidb@x \hbox {{\relax \tenrm (2.11)}}}
\global\def\_@eqn@AsymptoticRelationLargeScale{\relax \unhbox \voidb@x \hbox {{\relax \tenrm (2.12)}}}
\global\def\_@eqn@SpecialCaseARLargeScale{\relax \unhbox \voidb@x \hbox {{\relax \tenrm (2.13)}}}
\global\def\_@proc@AsymptoticRelationLargeScale{2.1.5}
\global\def\_@subhead@TheSmallScaleFormalSolutions{2.2}
\global\def\_@eqn@ControlOfParameters{\relax \unhbox \voidb@x \hbox {{\relax \tenrm (2.14)}}}
\global\def\_@eqn@DependenceOfICOnLogCoeff{\relax \unhbox \voidb@x \hbox {{\relax \tenrm (2.15)}}}
\global\def\_@eqn@ErrorForExactSolutionWithGivenLogarithmicParameter{\relax \unhbox \voidb@x \hbox {{\relax \tenrm (2.16)}}}
\global\def\_@proc@ErrorForExactSolutionWithGivenLogarithmicParameter{2.2.1}
\global\def\_@eqn@ReparametrisedGrimSurfaceEquation{\relax \unhbox \voidb@x \hbox {{\relax \tenrm (2.17)}}}
\global\def\_@eqn@DefinitionOfD{\relax \unhbox \voidb@x \hbox {{\relax \tenrm (2.18)}}}
\global\def\_@eqn@FormalSum{\relax \unhbox \voidb@x \hbox {{\relax \tenrm (2.19)}}}
\global\def\_@eqn@FormalToFunctionCorrespondence{\relax \unhbox \voidb@x \hbox {{\relax \tenrm (2.20)}}}
\global\def\_@eqn@FormalOperators{\relax \unhbox \voidb@x \hbox {{\relax \tenrm (2.21)}}}
\global\def\_@proc@LemmaFormalSolutionOverSmallScale{2.2.2}
\global\def\_@eqn@EqnRecurrenceForSmallScaleFormalSolutions{\relax \unhbox \voidb@x \hbox {{\relax \tenrm (2.22)}}}
\global\def\_@subhead@TheSmallScaleExactSolutionsI{2.3}
\global\def\_@eqn@FunctionOfPartialSum{\relax \unhbox \voidb@x \hbox {{\relax \tenrm (2.23)}}}
\global\def\_@eqn@DefinitionOfSigmaK{\relax \unhbox \voidb@x \hbox {{\relax \tenrm (2.24)}}}
\global\def\_@eqn@DefinitionOfSigma{\relax \unhbox \voidb@x \hbox {{\relax \tenrm (2.25)}}}
\global\def\_@eqn@FirstOrderNorm{\relax \unhbox \voidb@x \hbox {{\relax \tenrm (2.26)}}}
\global\def\_@eqn@ComparisonOfBothNorms{\relax \unhbox \voidb@x \hbox {{\relax \tenrm (2.27)}}}
\global\def\_@eqn@NormOfDAndItsInverse{\relax \unhbox \voidb@x \hbox {{\relax \tenrm (2.28)}}}
\global\def\_@eqn@DefinitionOfErrorOperator{\relax \unhbox \voidb@x \hbox {{\relax \tenrm (2.29)}}}
\global\def\_@eqn@DerivativeOfErrorOperator{\relax \unhbox \voidb@x \hbox {{\relax \tenrm (2.30)}}}
\global\def\_@eqn@FormulaForErrorOperator{\relax \unhbox \voidb@x \hbox {{\relax \tenrm (2.31)}}}
\global\def\_@eqn@SizeOfErrorOperator{\relax \unhbox \voidb@x \hbox {{\relax \tenrm (2.32)}}}
\global\def\_@eqn@DefinitionOfContractionMapping{\relax \unhbox \voidb@x \hbox {{\relax \tenrm (2.33)}}}
\global\def\_@eqn@ContractionMappingInequality{\relax \unhbox \voidb@x \hbox {{\relax \tenrm (2.34)}}}
\global\def\_@eqn@InitialConditions{\relax \unhbox \voidb@x \hbox {{\relax \tenrm (2.35)}}}
\global\def\_@eqn@GlobalEstimateOfError{\relax \unhbox \voidb@x \hbox {{\relax \tenrm (2.36)}}}
\global\def\_@proc@LemmaGlobalEstimateOfError{2.3.4}
\global\def\_@subhead@TheSmallScaleExactSolutionsII{2.4}
\global\def\_@eqn@DifferenceBetweenFormalAndExactToAllOrders{\relax \unhbox \voidb@x \hbox {{\relax \tenrm (2.37)}}}
\global\def\_@eqn@HigherDerivativesOfDifference{\relax \unhbox \voidb@x \hbox {{\relax \tenrm (2.38)}}}
\global\def\_@eqn@BetterInitialConditions{\relax \unhbox \voidb@x \hbox {{\relax \tenrm (2.39)}}}
\global\def\_@eqn@EstimateOfErrorInLogCoordinates{\relax \unhbox \voidb@x \hbox {{\relax \tenrm (2.40)}}}
\global\def\_@proc@EstimateOfErrorInLogCoordinates{2.4.2}
\global\def\_@subhead@TheSmallScaleJacobiFields{2.5}
\global\def\_@eqn@ErrorForExactJacobiFieldWithGivenLogarithmicParameter{\relax \unhbox \voidb@x \hbox {{\relax \tenrm (2.41)}}}
\global\def\_@proc@ErrorForExactJacobiFieldWithGivenLogarithmicParameter{2.5.1}
\global\def\_@eqn@FormalNDerivative{\relax \unhbox \voidb@x \hbox {{\relax \tenrm (2.42)}}}
\global\def\_@eqn@FormalJacobiField{\relax \unhbox \voidb@x \hbox {{\relax \tenrm (2.43)}}}
\global\def\_@eqn@FormalJacobiEquation{\relax \unhbox \voidb@x \hbox {{\relax \tenrm (2.44)}}}
\global\def\_@eqn@PartialSumOfFormalJacobiField{\relax \unhbox \voidb@x \hbox {{\relax \tenrm (2.45)}}}
\global\def\_@eqn@FormalJacobiFieldSolvesFormalJacobiEquation{\relax \unhbox \voidb@x \hbox {{\relax \tenrm (2.46)}}}
\global\def\_@eqn@FunctionIsApproximateJacobiField{\relax \unhbox \voidb@x \hbox {{\relax \tenrm (2.47)}}}
\global\def\_@eqn@ExactJacobiField{\relax \unhbox \voidb@x \hbox {{\relax \tenrm (2.48)}}}
\global\def\_@eqn@SizeOfErrorForJacobiFields{\relax \unhbox \voidb@x \hbox {{\relax \tenrm (2.49)}}}
\global\def\_@eqn@ImageOfDifference{\relax \unhbox \voidb@x \hbox {{\relax \tenrm (2.50)}}}
\global\def\_@eqn@ErrorForSecondDerivativeOfGrimProfile{\relax \unhbox \voidb@x \hbox {{\relax \tenrm (2.51)}}}
\global\def\_@head@TheGrimParaboloid{3}
\global\def\_@subhead@TheMCFSJacobiOperator{3.1}
\global\def\_@eqn@MetricOverPlane{\relax \unhbox \voidb@x \hbox {{\relax \tenrm (3.1)}}}
\global\def\_@eqn@ApproximateMetricOverPlane{\relax \unhbox \voidb@x \hbox {{\relax \tenrm (3.2)}}}
\global\def\_@eqn@LargeScaleWeight{\relax \unhbox \voidb@x \hbox {{\relax \tenrm (3.3)}}}
\global\def\_@eqn@DefinitionOfWeightedNormsOverG{\relax \unhbox \voidb@x \hbox {{\relax \tenrm (3.4)}}}
\global\def\_@eqn@ApproximateLargeScaleWeight{\relax \unhbox \voidb@x \hbox {{\relax \tenrm (3.5)}}}
\global\def\_@eqn@DefinitionOfWeightedSpacesOverG{\relax \unhbox \voidb@x \hbox {{\relax \tenrm (3.6)}}}
\global\def\_@proc@JIsLinearIsomorphismOverG{3.1.1}
\global\def\_@eqn@PhiWeightedJacobiOperator{\relax \unhbox \voidb@x \hbox {{\relax \tenrm (3.7)}}}
\global\def\_@eqn@FormulaForConjugateOperator{\relax \unhbox \voidb@x \hbox {{\relax \tenrm (3.8)}}}
\global\def\_@eqn@CoefficientsOfConjugateOperator{\relax \unhbox \voidb@x \hbox {{\relax \tenrm (3.9)}}}
\global\def\_@eqn@IntrinsicCurvatureOfCircles{\relax \unhbox \voidb@x \hbox {{\relax \tenrm (3.10)}}}
\global\def\_@eqn@EqnAsymptoticBoundsOnLittleKappa{\relax \unhbox \voidb@x \hbox {{\relax \tenrm (3.11)}}}
\global\def\_@eqn@AsymptoticFormulaForRho{\relax \unhbox \voidb@x \hbox {{\relax \tenrm (3.12)}}}
\global\def\_@eqn@ConjugateJacobiOperatorInPolarCoordinates{\relax \unhbox \voidb@x \hbox {{\relax \tenrm (3.13)}}}
\global\def\_@eqn@ZeroethOrderTermInConjugateJacobiOperatorInPolarCoordinates{\relax \unhbox \voidb@x \hbox {{\relax \tenrm (3.14)}}}
\global\def\_@subhead@InvertibilityOverSobolevSpaces{3.2}
\global\def\_@eqn@EqnSobolevMassConcentration{\relax \unhbox \voidb@x \hbox {{\relax \tenrm (3.15)}}}
\global\def\_@proc@SobolevControlAtInfinity{3.2.1}
\global\def\_@eqn@EqnSobolevEllipticEstimate{\relax \unhbox \voidb@x \hbox {{\relax \tenrm (3.16)}}}
\global\def\_@proc@LemmaSobolevFredholm{3.2.3}
\global\def\_@proc@LemmaNoBoundedJacobiFields{3.2.4}
\global\def\_@proc@LemmaSobolevKernelIsEmpty{3.2.5}
\global\def\_@proc@LemmaInvertibilityOverSobolev{3.2.6}
\global\def\_@subhead@InvertibilityOverHoelderSpaces{3.3}
\global\def\_@eqn@ViscositySolutionsToODEInequality{\relax \unhbox \voidb@x \hbox {{\relax \tenrm (3.17)}}}
\global\def\_@proc@LemmaViscositySolutionsToODEInequality{3.3.1}
\global\def\_@eqn@EqnHolderMassConcentration{\relax \unhbox \voidb@x \hbox {{\relax \tenrm (3.18)}}}
\global\def\_@eqn@EqnHolderEllipticEstimate{\relax \unhbox \voidb@x \hbox {{\relax \tenrm (3.19)}}}
\global\def\_@proc@ThmInvertibilityOverHolder{3.3.4}
\global\def\_@head@GrimEnds{4}
\global\def\_@subhead@TheModifiedMCFSJacobiOperator{4.1}
\global\def\_@eqn@AsymptoticFormulaForGrimEndSmallScale{\relax \unhbox \voidb@x \hbox {{\relax \tenrm (4.1)}}}
\global\def\_@eqn@DefinitionHybridNorm{\relax \unhbox \voidb@x \hbox {{\relax \tenrm (4.2)}}}
\global\def\_@eqn@DefinitionOfWeightSmallScale{\relax \unhbox \voidb@x \hbox {{\relax \tenrm (4.3)}}}
\global\def\_@eqn@SmallScaleModifiedJacobiOperator{\relax \unhbox \voidb@x \hbox {{\relax \tenrm (4.4)}}}
\global\def\_@proc@JIsLinearIsomorphismOverEnds{4.1.1}
\global\def\_@eqn@ExplicitFormulaForModifiedJacobiOperator{\relax \unhbox \voidb@x \hbox {{\relax \tenrm (4.5)}}}
\global\def\_@eqn@ModifiedJacobiOperatorOfGrimEnd{\relax \unhbox \voidb@x \hbox {{\relax \tenrm (4.6)}}}
\global\def\_@eqn@ModifiedJacobiOperatorOfGrimEndErrorTerms{\relax \unhbox \voidb@x \hbox {{\relax \tenrm (4.7)}}}
\global\def\_@proc@ModifiedJacobiOperatorOfGrimEnd{4.1.3}
\global\def\_@subhead@TheRegularComponent{4.2}
\global\def\_@eqn@AsymptoticFormulaForGrimParaboloidSmallScale{\relax \unhbox \voidb@x \hbox {{\relax \tenrm (4.8)}}}
\global\def\_@eqn@ModifiedJacobiOperatorOfGrimParaboloid{\relax \unhbox \voidb@x \hbox {{\relax \tenrm (4.9)}}}
\global\def\_@eqn@ModifiedJacobiOperatorOfGrimParaboloidErrorTerms{\relax \unhbox \voidb@x \hbox {{\relax \tenrm (4.10)}}}
\global\def\_@eqn@PhiWeightedModifiedJacoiOperator{\relax \unhbox \voidb@x \hbox {{\relax \tenrm (4.11)}}}
\global\def\_@eqn@PhiWeightedModifiedJacobiOperatorOfParaboloid{\relax \unhbox \voidb@x \hbox {{\relax \tenrm (4.12)}}}
\global\def\_@eqn@DefinitionOfDAndE{\relax \unhbox \voidb@x \hbox {{\relax \tenrm (4.13)}}}
\global\def\_@eqn@ExtensionOfModifiedJacobiOperator{\relax \unhbox \voidb@x \hbox {{\relax \tenrm (4.14)}}}
\global\def\_@eqn@DifferenceOperatorErrorTerms{\relax \unhbox \voidb@x \hbox {{\relax \tenrm (4.15)}}}
\global\def\_@eqn@CoefficientsOfDConvergeToZeroOnSmallestScale{\relax \unhbox \voidb@x \hbox {{\relax \tenrm (4.16)}}}
\global\def\_@proc@CoefficientsOfDConvergeToZeroOnSmallestScale{4.2.1}
\global\def\_@eqn@CoefficientsOfDConvergeToZeroOnSmallScale{\relax \unhbox \voidb@x \hbox {{\relax \tenrm (4.17)}}}
\global\def\_@proc@CoefficientsOfDConvergeToZeroOnSmallScale{4.2.2}
\global\def\_@eqn@InequalityPreservedOverShortInterval{\relax \unhbox \voidb@x \hbox {{\relax \tenrm (4.18)}}}
\global\def\_@proc@InequalityPreservedOverShortInterval{4.2.3}
\global\def\_@eqn@CoefficientsOfDConvergeToZeroOnLowerCentralScale{\relax \unhbox \voidb@x \hbox {{\relax \tenrm (4.19)}}}
\global\def\_@eqn@CoefficientsOfDConvergeToZeroOnLargeScale{\relax \unhbox \voidb@x \hbox {{\relax \tenrm (4.20)}}}
\global\def\_@eqn@CoefficientsOfDConvergeToZeroOnUpperCentralScale{\relax \unhbox \voidb@x \hbox {{\relax \tenrm (4.21)}}}
\global\def\_@proc@TheOperatorDConvergesToZero{4.2.7}
\global\def\_@subhead@TheSingularPart{4.3}
\global\def\_@eqn@DefinitionOfAForOperatorE{\relax \unhbox \voidb@x \hbox {{\relax \tenrm (4.22)}}}
\global\def\_@eqn@FormulaForAOverCentralBall{\relax \unhbox \voidb@x \hbox {{\relax \tenrm (4.23)}}}
\global\def\_@eqn@HybridProperty{\relax \unhbox \voidb@x \hbox {{\relax \tenrm (4.24)}}}
\global\def\_@proc@HybridProperty{4.3.1}
\global\def\_@proc@TheOperatorEConvergesToZero{4.3.4}
\global\def\_@head@SurgeryAndThePerturbationFamily{5}
\global\def\_@subhead@TheBasicSurgeryOperation{5.1}
\global\def\_@eqn@ProfileOfCatenoidalEnd{\relax \unhbox \voidb@x \hbox {{\relax \tenrm (5.1)}}}
\global\def\_@eqn@ProfileOfGrimEnd{\relax \unhbox \voidb@x \hbox {{\relax \tenrm (5.2)}}}
\global\def\_@eqn@DefinitionOfJoinedSurface{\relax \unhbox \voidb@x \hbox {{\relax \tenrm (5.3)}}}
\global\def\_@eqn@FormulaForJoinedSurface{\relax \unhbox \voidb@x \hbox {{\relax \tenrm (5.4)}}}
\global\def\_@subhead@TheDeformationFamily{5.2}
\global\def\_@eqn@ModifiedNormalVectorField{\relax \unhbox \voidb@x \hbox {{\relax \tenrm (5.5)}}}
\global\def\_@eqn@PerturbationFamily{\relax \unhbox \voidb@x \hbox {{\relax \tenrm (5.6)}}}
\global\def\_@subhead@MicroscopicAndMacroscopicPerturbations{5.3}
\global\def\_@eqn@MCFSFunctionalAgain{\relax \unhbox \voidb@x \hbox {{\relax \tenrm (5.7)}}}
\global\def\_@eqn@PerturbationOperators{\relax \unhbox \voidb@x \hbox {{\relax \tenrm (5.8)}}}
\global\def\_@subhead@ModifiedJacobiOperators{5.4}
\global\def\_@eqn@GeneralModifiedJacobiOperator{\relax \unhbox \voidb@x \hbox {{\relax \tenrm (5.9)}}}
\global\def\_@eqn@ExplicitFormulaOfJacobiOperatorForGraphs{\relax \unhbox \voidb@x \hbox {{\relax \tenrm (5.10)}}}
\global\def\_@eqn@ValuesOfYVAndZW{\relax \unhbox \voidb@x \hbox {{\relax \tenrm (5.11)}}}
\global\def\_@eqn@JacobiOperatorOfCHM{\relax \unhbox \voidb@x \hbox {{\relax \tenrm (5.12)}}}
\global\def\_@eqn@ErrorTermsInJacobiOperatorOfCHM{\relax \unhbox \voidb@x \hbox {{\relax \tenrm (5.13)}}}
\global\def\_@eqn@DefinitionOfWeightPsi{\relax \unhbox \voidb@x \hbox {{\relax \tenrm (5.14)}}}
\global\def\_@eqn@ModifiedMCFSJacobiOperatorOfG{\relax \unhbox \voidb@x \hbox {{\relax \tenrm (5.15)}}}
\global\def\_@proc@DefnsOfModifiedMCFSJacobiOpMatch{5.4.2}
\global\def\_@eqn@JacobiOperatorOfRescaledGrimEnd{\relax \unhbox \voidb@x \hbox {{\relax \tenrm (5.16)}}}
\global\def\_@eqn@ErrorTermsInJacobiOperatorOfRescaledGrimEnd{\relax \unhbox \voidb@x \hbox {{\relax \tenrm (5.17)}}}
\global\def\_@eqn@ErrorFromJacobiOperatorOfJoinedSurface{\relax \unhbox \voidb@x \hbox {{\relax \tenrm (5.18)}}}
\global\def\_@eqn@ErrorTermsInCommutators{\relax \unhbox \voidb@x \hbox {{\relax \tenrm (5.19)}}}
\global\def\_@proc@ErrorTermsInCommutators{5.4.5}
\global\def\_@subhead@ControllingMacroscopicPerturbations{5.5}
\global\def\_@eqn@NormOfDeficiencyTermInLowerRegion{\relax \unhbox \voidb@x \hbox {{\relax \tenrm (5.20)}}}
\global\def\_@eqn@NormOfDeficiencyTermOfJoinedSurface{\relax \unhbox \voidb@x \hbox {{\relax \tenrm (5.21)}}}
\global\def\_@eqn@DifferenceBetweenZAndW{\relax \unhbox \voidb@x \hbox {{\relax \tenrm (5.22)}}}
\global\def\_@head@ConstructingTheGreensOperator{6}
\global\def\_@subhead@TheCylindricalAndGrimNorms{6.1}
\global\def\_@eqn@DefinitionOfCylindricalDerivative{\relax \unhbox \voidb@x \hbox {{\relax \tenrm (6.1)}}}
\global\def\_@eqn@DefinitionOfFractionalCylindricalDerivative{\relax \unhbox \voidb@x \hbox {{\relax \tenrm (6.2)}}}
\global\def\_@eqn@DefinitionOfCylindricalNorms{\relax \unhbox \voidb@x \hbox {{\relax \tenrm (6.3)}}}
\global\def\_@eqn@ScaleFreeNormOverCosta{\relax \unhbox \voidb@x \hbox {{\relax \tenrm (6.4)}}}
\global\def\_@proc@RightInverseOverCosta{6.1.1}
\global\def\_@eqn@DefinitionOfGrimNorms{\relax \unhbox \voidb@x \hbox {{\relax \tenrm (6.5)}}}
\global\def\_@eqn@HybridNormII{\relax \unhbox \voidb@x \hbox {{\relax \tenrm (6.6)}}}
\global\def\_@proc@RightInverseOverGrim{6.1.2}
\global\def\_@eqn@DefinitionOfGrimDerivative{\relax \unhbox \voidb@x \hbox {{\relax \tenrm (6.7)}}}
\global\def\_@eqn@DefinitionOfFractionalGrimDerivative{\relax \unhbox \voidb@x \hbox {{\relax \tenrm (6.8)}}}
\global\def\_@eqn@OtherFormulaForGrimNorm{\relax \unhbox \voidb@x \hbox {{\relax \tenrm (6.9)}}}
\global\def\_@eqn@DefinitionOfCylindricalAndGrimVolumeForms{\relax \unhbox \voidb@x \hbox {{\relax \tenrm (6.10)}}}
\global\def\_@subhead@ThePingPongArgumentPartI{6.2}
\global\def\_@eqn@SeminormsOfJoinedSurface{\relax \unhbox \voidb@x \hbox {{\relax \tenrm (6.11)}}}
\global\def\_@eqn@DefinitionOfUpwardOperator{\relax \unhbox \voidb@x \hbox {{\relax \tenrm (6.12)}}}
\global\def\_@eqn@UsefulFormalaForUpwardOperator{\relax \unhbox \voidb@x \hbox {{\relax \tenrm (6.13)}}}
\global\def\_@eqn@NormOfA{\relax \unhbox \voidb@x \hbox {{\relax \tenrm (6.14)}}}
\global\def\_@proc@NormOfA{6.2.1}
\global\def\_@eqn@HolderNormOfJPhiInUpperRegion{\relax \unhbox \voidb@x \hbox {{\relax \tenrm (6.15)}}}
\global\def\_@eqn@FormulaUsedAlsoInSobolevCase{\relax \unhbox \voidb@x \hbox {{\relax \tenrm (6.16)}}}
\global\def\_@eqn@SobolevNormOfJPhiInUpperRegion{\relax \unhbox \voidb@x \hbox {{\relax \tenrm (6.17)}}}
\global\def\_@eqn@HolderNormOfCommutatorInUpperRegion{\relax \unhbox \voidb@x \hbox {{\relax \tenrm (6.18)}}}
\global\def\_@eqn@SecondFormulaAlsoUsedInSobolevCase{\relax \unhbox \voidb@x \hbox {{\relax \tenrm (6.19)}}}
\global\def\_@eqn@SobolevNormOfCommutatorInUpperRegion{\relax \unhbox \voidb@x \hbox {{\relax \tenrm (6.20)}}}
\global\def\_@eqn@NormOfPhiInUpperRegion{\relax \unhbox \voidb@x \hbox {{\relax \tenrm (6.21)}}}
\global\def\_@subhead@ThePingPongArgumentPartII{6.3}
\global\def\_@eqn@NormOfPsi{\relax \unhbox \voidb@x \hbox {{\relax \tenrm (6.22)}}}
\global\def\_@eqn@DefinitionOfDownwardOperator{\relax \unhbox \voidb@x \hbox {{\relax \tenrm (6.23)}}}
\global\def\_@eqn@UsefulFormulaForDownwardOperator{\relax \unhbox \voidb@x \hbox {{\relax \tenrm (6.24)}}}
\global\def\_@eqn@NormOfB{\relax \unhbox \voidb@x \hbox {{\relax \tenrm (6.25)}}}
\global\def\_@proc@NormOfB{6.3.1}
\global\def\_@eqn@NormOfV{\relax \unhbox \voidb@x \hbox {{\relax \tenrm (6.26)}}}
\global\def\_@proc@NormOfV{6.3.2}
\global\def\_@eqn@NormOfJPsiInLowerRegion{\relax \unhbox \voidb@x \hbox {{\relax \tenrm (6.27)}}}
\global\def\_@eqn@NormOfPsiInLowerRegion{\relax \unhbox \voidb@x \hbox {{\relax \tenrm (6.28)}}}
\global\def\_@eqn@NormOfCommutatorInLowerRegion{\relax \unhbox \voidb@x \hbox {{\relax \tenrm (6.29)}}}
\global\def\_@subhead@ThePingPongArgumentPartIII{6.4}
\global\def\_@eqn@DefinitionOfAdjustmentOperators{\relax \unhbox \voidb@x \hbox {{\relax \tenrm (6.30)}}}
\global\def\_@eqn@PreliminariesDefinitionOfRightInverse{\relax \unhbox \voidb@x \hbox {{\relax \tenrm (6.31)}}}
\global\def\_@eqn@PreliminaryOperatorsAreRightInverses{\relax \unhbox \voidb@x \hbox {{\relax \tenrm (6.32)}}}
\global\def\_@eqn@NormOfMultiplicationByChi{\relax \unhbox \voidb@x \hbox {{\relax \tenrm (6.33)}}}
\global\def\_@eqn@DefinitionOfRightInverse{\relax \unhbox \voidb@x \hbox {{\relax \tenrm (6.34)}}}
\global\def\_@eqn@ConstructedOperatorIsRightInverse{\relax \unhbox \voidb@x \hbox {{\relax \tenrm (6.35)}}}
\global\def\_@eqn@NormOfJoinedLAndV{\relax \unhbox \voidb@x \hbox {{\relax \tenrm (6.36)}}}
\global\def\_@eqn@NormOfW{\relax \unhbox \voidb@x \hbox {{\relax \tenrm (6.37)}}}
\global\def\_@proc@NormOfW{6.4.2}
\global\def\_@eqn@NormOfPInLowerRegion{\relax \unhbox \voidb@x \hbox {{\relax \tenrm (6.38)}}}
\global\def\_@proc@NormOfPInLowerRegion{6.4.3}
\global\def\_@eqn@NormOfPInUpperRegion{\relax \unhbox \voidb@x \hbox {{\relax \tenrm (6.39)}}}
\global\def\_@proc@NormOfPInUpperRegion{6.4.4}
\global\def\_@head@ExistenceAndEmbeddedness{7}
\global\def\_@subhead@TheSchauderFixedPointTheorem{7.1}
\global\def\_@eqn@ModifiedGrimNorms{\relax \unhbox \voidb@x \hbox {{\relax \tenrm (7.1)}}}
\global\def\_@eqn@EstimateOfMCFSOfJoinedSurface{\relax \unhbox \voidb@x \hbox {{\relax \tenrm (7.2)}}}
\global\def\_@eqn@SecondOrderPerturbationOverCHM{\relax \unhbox \voidb@x \hbox {{\relax \tenrm (7.3)}}}
\global\def\_@eqn@SecondOrderPerturbationOverGrim{\relax \unhbox \voidb@x \hbox {{\relax \tenrm (7.4)}}}
\global\def\_@proc@QuadraticError{7.1.2}
\global\def\_@eqn@EstimateOfError{\relax \unhbox \voidb@x \hbox {{\relax \tenrm (7.5)}}}
\global\def\_@proc@ExistenceTheorem{7.1.3}
\global\def\_@proc@EmbeddednessTheorem{7.1.4}
\global\def\_@head@TerminologyConventionsEtc{A}
\global\def\_@subhead@GeneralDefinitions{A.1}
\global\def\_@subhead@SurfaceGeometry{A.2}
\global\def\_@eqn@AppMCFSOperator{\relax \unhbox \voidb@x \hbox {{\relax \tenrm (A.1)}}}
\global\def\_@eqn@AppMCFSJacobiOperator{\relax \unhbox \voidb@x \hbox {{\relax \tenrm (A.2)}}}
\global\def\_@eqn@AppGradientAndHessian{\relax \unhbox \voidb@x \hbox {{\relax \tenrm (A.3)}}}
\global\def\_@eqn@AppCommutator{\relax \unhbox \voidb@x \hbox {{\relax \tenrm (A.4)}}}
\global\def\_@subhead@SurfaceGeometryOfGraphs{A.3}
\global\def\_@eqn@ThirdComponentOfNormalVector{\relax \unhbox \voidb@x \hbox {{\relax \tenrm (A.5)}}}
\global\def\_@eqn@AppBasicFormulaInCoordinateCharts{\relax \unhbox \voidb@x \hbox {{\relax \tenrm (A.6)}}}
\global\def\_@eqn@AppMCFSFunctionalForGraphs{\relax \unhbox \voidb@x \hbox {{\relax \tenrm (A.7)}}}
\global\def\_@subhead@FunctionSpaces{A.4}
\global\def\_@eqn@AppHolderSeminorm{\relax \unhbox \voidb@x \hbox {{\relax \tenrm (A.8)}}}
\global\def\_@eqn@AppTotalVariationBoundedByTwiceSupremum{\relax \unhbox \voidb@x \hbox {{\relax \tenrm (A.9)}}}
\global\def\_@eqn@AppFirstInterpolationInequality{\relax \unhbox \voidb@x \hbox {{\relax \tenrm (A.10)}}}
\global\def\_@eqn@AppSecondInterpolationInequality{\relax \unhbox \voidb@x \hbox {{\relax \tenrm (A.11)}}}
\global\def\_@eqn@AppProductRule{\relax \unhbox \voidb@x \hbox {{\relax \tenrm (A.12)}}}
\global\def\_@eqn@AppHolderNormOverUnionI{\relax \unhbox \voidb@x \hbox {{\relax \tenrm (A.13)}}}
\global\def\_@eqn@AppHolderNormOverUnionII{\relax \unhbox \voidb@x \hbox {{\relax \tenrm (A.14)}}}
\global\def\_@eqn@AppDefinitionOfFractionalDifferentialOperator{\relax \unhbox \voidb@x \hbox {{\relax \tenrm (A.15)}}}
\global\def\_@eqn@AppDefinitionOfHolderNorm{\relax \unhbox \voidb@x \hbox {{\relax \tenrm (A.16)}}}
\global\def\_@eqn@AppDefinitionOfLPNorm{\relax \unhbox \voidb@x \hbox {{\relax \tenrm (A.17)}}}
\global\def\_@eqn@AppDefinitionOfSobolevNorm{\relax \unhbox \voidb@x \hbox {{\relax \tenrm (A.18)}}}
\global\def\_@eqn@AppSobolevEmbeddingTheorem{\relax \unhbox \voidb@x \hbox {{\relax \tenrm (A.19)}}}
\global\def\_@eqn@MaximumOfExponential{\relax \unhbox \voidb@x \hbox {{\relax \tenrm (A.20)}}}
\global\def\_@eqn@IntegralOfExponential{\relax \unhbox \voidb@x \hbox {{\relax \tenrm (A.21)}}}
\global\def\_@subhead@EllipticEstimates{A.5}
\global\def\_@eqn@AppEllipticEstimate{\relax \unhbox \voidb@x \hbox {{\relax \tenrm (A.22)}}}
\global\def\_@proc@ThmEllipticToFredholm{A.5.1}
\global\def\_@head@Bibliography{B}
 }%
\else{\openout\references=references.tex }%
\fi%
%
%
\newcount\headno%
\global\headno=0%
\def\headinfo{\ifinappendices\alphanum\headno\else\the\headno\fi}%
\def\nextheadno{\global\advance\headno by 1 \global\subheadno=0 \global\procno=0 \global\eqnno=0 \headinfo}%
\makecounter{head}%
%
%
\newcount\subheadno%
\global\subheadno=0%
\def\subheadinfo{\headinfo.\the\subheadno}%
\def\nextsubheadno{\global\advance\subheadno by 1 \global\procno=0 \subheadinfo}%
\makecounter{subhead}%
%
%
\newcount\procno%
\global\procno=0%
\def\procinfo{\subheadinfo.\the\procno}%
\def\nextprocno{\global\advance\procno by 1 \procinfo}%
\makecounter{proc}%
%
%
\newcount\figno%
\global\figno=0%
\def\figinfo{\subheadinfo.\the\figno}%
\def\nextfigno{\global\advance\figno by 1 \figinfo}%
\makecounter{fig}%
%
%
\newcount\eqnno%
\global\eqnno=0%
\def\eqninfo{\text{{\rm (\headinfo.\the\eqnno)}}}%
\def\nexteqnno[#1]{\global\advance\eqnno by 1 \eqninfo\hbox{\eqnlabel{#1}}}%
\makecounter{eqn}%
%
%
%
%
%
\def\gobbleeight#1#2#3#4#5#6#7#8{}%
\newcount\citationno%
\global\citationno=0%
\def\citationinfo{\the\citationno}%
\makecounter{citation}%
\newwrite\biblio%
\def\newref#1#2{%
\def\temptext{#2}%
\edef\bibliotextoutput{\expandafter\gobbleeight\meaning\temptext}%
\global\advance\citationno by 1\citationlabel{#1}%
\ifmakebiblio%
    \edef\fileoutput{\write\biblio{\noindent\hbox to 0pt{\hss$[\the\citationno]$}\hskip 0.2em\bibliotextoutput\medskip}}%
    \fileoutput%
\fi}%
\def\cite#1{%
$[\citationref{#1}]$%
\ifmakebiblio%
    \edef\fileoutput{\write\biblio{#1}}%
    \fileoutput%
\fi%
}%
%
%
%
%
\let\mypar=\par%
\edef\Pagetitle={Blank}\headline={\hfil\Pagetitle\hfil}%
\edef\Pagefooter={Blank}\footline={\hfil\Pagefooter\hfil}%
%
%
\newcount\showpagenumflag%
\global\showpagenumflag=0 %
\def\nextoddpage%
{\newpage\ifodd\pageno%
\else\global\showpagenumflag=0 %
\null\vfil\eject%
\global\showpagenumflag=1 %
\fi}%
%
%
\font\headfont=cmb12%
\def\newhead#1[#2]%
{\ifhmode\mypar\fi%
\ifnum\headno=0 \else\goodbreak\bigskip\fi%
{\headfont\noindent\nextheadno\ - #1.}\headlabel{#2}%
\nobreak\medskip}%
%
%
\def\newsubhead#1[#2]%
{\ifhmode\mypar\fi%
\ifnum\subheadno=0 \else\goodbreak\medskip\fi%
{\bf\noindent\nextsubheadno\ - #1.\ }\subheadlabel{#2}}%
%
%
\newif\ifinproclaim%
\global\inproclaimfalse%
\def\proclaim#1{%
\goodbreak\medskip
\bgroup\inproclaimtrue%
\noindent{\bf #1}%
\nobreak\medskip\sl}%
\def\noskipproclaim#1{%
\goodbreak\medskip%
\bgroup\inproclaimtrue%
\noindent{\bf #1}\nobreak\sl}%
\def\endproclaim{\mypar\egroup\nobreak\medskip\ignorespaces}%
%
%
%
\newcount\xpos\newcount\ypos
\def\makelabelgrid{%
\xpos=-5 \ypos=-5 %
\loop\ifnum\xpos<6 %
{\loop\ifnum\ypos<6 %
\def\labeltext{x}%
\ifnum\xpos=0\def\labeltext{+}\fi%
\ifnum\ypos=0\def\labeltext{+}\fi%
\placelabel[\xpos][\ypos]{\labeltext}%
\advance\ypos by 1 %
\repeat}%
\advance\xpos by 1 %
\repeat}%
\def\placelabel[#1][#2]#3{{%
\setbox10=\hbox{\raise #2cm \hbox{\hskip #1cm #3}}%
\ht10=0pt \dp10=0pt \wd10=0pt \box10}}%
%
%
%
%
\def\myitem#1{\noindent\hbox to .5cm{\hfill#1\hss}}%
%
%
%
%
%
%
%
%
%
\font\sansseriften=cmss10%
\font\sansserifseven=cmss7%
\font\sansseriffive=cmss5%
\newfam\sansseriffam%
\textfont\sansseriffam=\sansseriften%
\scriptfont\sansseriffam=\sansserifseven%
\scriptscriptfont\sansseriffam=\sansseriffive%
\def\mathsf{\fam\sansseriffam}%
%
%
%
\font\boldten=cmb10%
\font\boldseven=cmb7%
\font\boldfive=cmb5%
\newfam\mathboldfam%
\textfont\mathboldfam=\boldten%
\scriptfont\mathboldfam=\boldseven%
\scriptscriptfont\mathboldfam=\boldfive%
\def\mathbf{\fam\mathboldfam}%
%
%
%
\font\mycmmiten=cmmi10%
\font\mycmmiseven=cmmi7%
\font\mycmmifive=cmmi5%
\newfam\mycmmifam%
\textfont\mycmmifam=\mycmmiten%
\scriptfont\mycmmifam=\mycmmiseven%
\scriptscriptfont\mycmmifam=\mycmmifive%
\def\hexa#1{\ifcase #1 0\or 1\or 2\or 3\or 4\or 5\or 6\or 7\or 8\or 9\or A\or B\or C\or D\or E\or F\fi}%
\mathchardef\mathi="7\hexa\mycmmifam7B%
\mathchardef\mathj="7\hexa\mycmmifam7C%
%
%
\font\mymsbmten=msbm10 at 8pt%
\font\mymsbmseven=msbm7 at 5.6pt
\font\mymsbmfive=msbm5 at 4pt%
\newfam\mymsbmfam%
\textfont\mymsbmfam=\mymsbmten%
\scriptfont\mymsbmfam=\mymsbmseven%
\scriptscriptfont\mymsbmfam=\mymsbmfive%
\mathchardef\mybeth="7\hexa\mymsbmfam69%
\mathchardef\mygimmel="7\hexa\mymsbmfam6A%
\mathchardef\mydaleth="7\hexa\mymsbmfam6B%
%
%
%
%
\def\proof{{\noindent\bf Proof:\ }}%
\def\remark{{\noindent\bf Remark:\ }}%
\def\qed{~$\square$}%
\def\makeop#1{\global\expandafter\def\csname op#1\endcsname{{\text{#1}}}}%
\def\makeopsmall#1{\global\expandafter\def\csname op#1\endcsname{{\text{\lowercase{#1}}}}}%
%
%
\def\munion{\mathop{\cup}}%
\def\minter{\mathop{\cap}}%
%
%
\makeop{Ext}%
\makeop{Int}%
\makeop{Dist}%
\makeop{Diam}%
\makeop{Length}%
%
%
%
%
%
%
%
%
\def\msup{\mathop{{\text{Sup}}}}%
%
%
%
\makeop{Dim}%
\makeop{Ker}%
\makeop{Coker}%
\makeop{Tr}%
\makeop{Adj}%
\makeop{Det}%
\makeop{End}%
\makeop{Lin}%
\makeop{Symm}%
\makeop{Mult}%
%
%
\makeop{dx}%
\makeop{dy}%
\makeop{dz}%
\makeop{dt}%
\makeop{dVol}%
\makeop{dArea}%
\makeop{Supp}%
\makeop{Hess}%
\makeop{Lip}%
%
%
\makeop{Re}%
\makeop{Im}%
\makeop{Arg}%
\makeop{Log}%
\makeop{Exp}%
%
%
\makeopsmall{Cos}%
\makeopsmall{Sin}%
\makeopsmall{Tan}%
\makeopsmall{Sec}%
\makeopsmall{Cosec}%
\makeopsmall{Cot}%
\makeopsmall{ArcCos}%
\makeopsmall{ArcSin}%
\makeopsmall{ArcTan}%
\makeopsmall{ArcSec}%
\makeopsmall{ArcCosec}%
\makeopsmall{ArcCot}%
%
%
\makeopsmall{Cosh}%
\makeopsmall{Sinh}%
\makeopsmall{Tanh}%
\makeopsmall{ArcCosh}%
\makeopsmall{ArcSinh}%
\makeopsmall{ArcTanh}%
%
%
\makeop{Vol}%
\makeop{Area}%
\makeop{Riem}%
\makeop{Ric}%
\makeop{Scal}%
\makeop{Euc}%
\makeop{Imm}%
\makeop{Emb}%
%
%
\makeop{Id}%
\makeop{Ad}%
\makeop{O}%
\makeop{SO}%
\makeop{SL}%
\makeop{GL}%
\makeop{Conf}%
\makeop{Homeo}%
\makeop{Diff}%
\makeop{Isom}%
%
%
\makeop{Ind}%
\makeop{Sig}%
\makeop{Spec}%
%
%
\makeop{Conv}%
\makeop{Max}%
\makeop{Min}%
\makeop{Mod}%
\makeop{Deg}%
\makeop{loc}%
%
%
%
%
%
%
%
%
%
%
%
%
%
 %
%
%
%
%
%
\newref{AltschulerWu}{Altschuler S., Wu L. F., Translating surfaces of the non-parametric mean curvature flow with prescribed contact angle, {\sl Calc. Var. Partial Differential Equations}, {\bf 2}, (1994), 101--111}
\newref{ClutterbuckSchnurerSchulze}{Clutterbuck J., Schn\"urer O., Schulze F., Stability of translating solutions
to mean curvature flow, {\sl Calc. Var. Partial Differential Equations}, {\bf 29}, (2007), 281--293}
\newref{ChengZhou}{Cheng X., Zhou D., Stability properties and gap theorem for complete $f$-minimal hypersurfaces, arXiv 1307.5099}
\newref{ChengMejiaZhouI}{Cheng X., Mejia T., Zhou D., Simons' type equation for $f$-minimal hypersurfaces and applications, arXiv:1305.2379}
\newref{ChengMejiaZhouII}{Cheng X., Mejia T., Zhou D., Eigenvalue estimate and compactness for closed $f$-minimal surfaces, arXiv:1210.8448}
\newref{ChengMejiaZhouIII}{Cheng X., Mejia T., Zhou D., Stability and compactness for complete $f$-minimal surfaces, to appear in {\sl Transactions of Amer. Math. Soc.}, arXiv:1210.8076}
\newref{DavilaDelPinoNguyen}{D\'avila J., Del Pino M., Ngueyn X. H., Finite topology self-translating surfaces for the mean curvature flow in $\Bbb{R}^3$, {\sl Adv. Math.}, {\bf 320}, 674--729}
\newref{GilbTrud}{Gilbarg D., Trudinger N. S., {\sl Elliptic partial differential equations of second order},
Classics in Mathematics, Springer-Verlag, Berlin, 2001}
\newref{Halldorsson}{Halldorsson H. P., Helicoidal surfaces rotating/translating under the mean curvature flow, {\sl Geom. Dedicata}, {\bf 162}, (2013), 45--65}
\newref{HauswirthPacard}{Hauswirth L., Pacard F., Minimal surfaces of finite genus with two limit ends, {\sl Invent. Math.}, {\bf 169}, no. 3, (2007), 569--620}
\newref{HoffmanMeeks}{Hoffman D., Meeks W.H., Embedded minimal surfaces of finite topology, {\sl Ann. of Math.}, {\bf 131}, (1990), 1--34}
\newref{KapouleasI}{Kapouleas N., Complete constant mean curvature surfaces in Euclidean three space, {\sl Ann. of Math.}, {\bf 121}, (1990), 239--330}
\newref{KapouleasII}{Kapouleas N., Constant mean curvature surfaces constructed by fusing Wente tori, {\sl Invent. Math.}, {\bf 119}, (1995), 443--518}
\newref{KapouleasIII}{Kapouleas N., Complete embedded minimal surfaces of finite total curvature, {\sl J. Diff Geom.}, {\bf 47}, (1997), 95--169}
\newref{MartinSavasSmoczyk}{Martin F., Savas-Halilaj A., Smoczyk K., On the topology of translating solitons of the mean curvature flow, arXiv:1404.6703}
\newref{MazzeoPacard}{Mazzeo R., Pacard F., Constant mean curvature surfaces with Delaunay ends, {\sl Comm. Anal. Geom.}, {\bf 9}, no. 1, (2001), 169--237}
\newref{MazzeoPollack}{Mazzeo R., Pollack D., Gluing and moduli for noncompact geometric problems, in {\sl Geometric theory of singular phenomena in partial differential equations}, Symposia Mathematica Vol. XXXII, Cambridge Univ. Press, (1998), pp. 17--51}
\newref{Morabito}{Morabito F., Index and nullity of the Gauss map of the Costa-Hoffman-Meeks surfaces, {\sl  Indiana university mathematics journal}, {\bf 58}, no. 2, (2009), 677--708}
\newref{Nayatani}{Nayatani S., Morse Index and Gauss maps of complete minimal surfaces in Euclidean 3-space, {\sl Comment. Math. Helv.}, {\bf 68}, no. 4, (1993), 511--537}
\newref{NguyenI}{Nguyen X. H., Complete embedded self-translating surfaces under mean curvature flow, {\sl J. Geom. Anal.}, {\bf 23}, (2013), 1379--1426}
\newref{NguyenII}{Nguyen X. H., Translating tridents, {\sl Comm. Partial Differential Equations}, {\bf 34}, (2009), 257--280.}
\newref{Pacard}{}
\newref{Schwarz}{Schwarz M., {\sl Morse homology}, Progress in Mathematics, {\bf 111}, Birkh\"auser Verlag, Basel, (1993)}
\newref{Wang}{Wang X. J., Convex solutions to the mean curvature flow, {\sl Ann. of Math.}, {\bf 173}, (2011), 1185--1239}
\newref{Weber}{Weber M., Classical Minimal Surfaces in Euclidean Space by Examples: geometric and computational aspects of the weierstrass representation, in {\sl Global Theory of Minimal Surfaces} (David Hoffman Ed.), Clay Mathematics Proceedings, Vol. 2, (2004)}
\catcode`\@=11
\def\multiline#1{\null\,\vcenter{\openup1\jot \m@th %
\ialign{\strut$\displaystyle{##}$\hfil\crcr#1\crcr}}\,}
\def\triplealign#1{\null\,\vcenter{\openup1\jot \m@th %
\ialign{\strut\hfil$\displaystyle{##}\quad$&\hfil$\displaystyle{{}##}$&$\displaystyle{{}##}$\hfil\crcr#1\crcr}}\,}
\def\tripleeqalign#1{\null\,\vcenter{\openup1\jot \m@th %
\ialign{\strut$\displaystyle{##}\quad$\hfil&$\displaystyle{{}##}$\hfil&$\displaystyle{{}##}$\hfil\crcr#1\crcr}}\,}
\catcode`\@=12
\makeop{Ln}
\makeop{II}
\def\opCyl{{\text{\rm SF}}}
\makeop{cyl}
\makeopsmall{Arccosh}
\makeopsmall{Log}
\makeopsmall{Sec}
\def\opO{{\text{\rm O}}}
\def\Pagetitle{\ifnum\pageno=1\null\hfill\else\hfil{\rm On complete embedded translating solitons...}\hfil\fi}
\def\Pagefooter{\hfil{\myfontdefault\folio}\hfil}
\font\tablefont=cmr7
\newif\ifshowaddress\showaddresstrue
%
\def\centre{\rightskip=0pt plus 1fil \leftskip=0pt plus 1fil \spaceskip=.3333em \xspaceskip=.5em \parfillskip=0em \parindent=0em}%
\def\textmonth#1{\ifcase#1\or January\or Febuary\or March\or April\or May\or June\or July\or August\or September\or October\or November\or December\fi}
\font\abstracttitlefont=cmr10 at 14pt {\abstracttitlefont\centre On complete embedded translating solitons of the mean curvature flow that are of finite genus.\par}
\bigskip
{\centre 16th January 2015\par}
\bigskip
{\centre Graham Smith\footnote{${}^*$}{{\tablefont Instituto de Matem\'atica, UFRJ, Av. Athos da Silveira Ramos 149, Centro de Tecnologia - Bloco C, Cidade Universit\'aria - Ilha do Fund\~ao, Caixa Postal 68530, 21941-909, Rio de Janeiro, RJ - BRAZIL\hfill}}\par}
\bigskip
\noindent{\bf Abstract:\ }We desingularise the union of $3$ Grim paraboloids along Costa-Hoffman-Meeks surfaces in order to obtain complete embedded translating solitons of the mean curvature flow with $3$ ends and arbitrary finite genus.
\bigskip
\noindent{\bf Key Words:\ }Mean curvature flow, soliton, singular perturbation, Costa-Hoffman-Meeks surface, Grim paraboloid, bowl soliton.
\bigskip
\noindent{\bf AMS Subject Classification:\ } 53C44
\bigskip
\myfontdefault
%
%
\newhead{Introduction}[Introduction]
\newsubhead{Main Result}[BackgroundAndMainResult]
For the purposes of this paper, a {\bf mean curvature flow (MCF) soliton} is a complete surface in $\Bbb{R}^3$ whose evolution under the mean curvature flow is given by translation. In other words, up to rescalings and rigid motions of the ambient spacetime, it is a solution of what we will call the {\bf MCFS equation}
$$
H + \langle N,e_z\rangle = 0,\eqnum{\nexteqnno[IntroMCFSEquation]}
$$
where $H$ here denotes the mean curvature of the surface, $N$ its unit normal vector field, and $e_z$ the unit vector in the direction of the $z$-axis. Since MCF solitons are also minimal surfaces with respect to the conformally rescaled metric
$$
g_{ij} := e^{-z}\delta_{ij},
$$
many techniques used in the construction of minimal surfaces ought also to yield constructions of MCF solitons. We refer the reader to the review \cite{MartinSavasSmoczyk} of Mart\'\i n, Savas-Halilaj and Smoczyk for a good overview of the theory of MCF solitons at the time of writing.
\par
We use surgery to construct embedded MCF solitons with $3$ ends and arbitrary finite genus. Before stating our result, we describe the two components of our construction. First, given a positive integer $g$, the {\bf Costa-Hoffman-Meeks (CHM) surface} of genus $g$, denoted $C_g$, is a properly embedded minimal surface in $\Bbb{R}^3$ with $3$ ends, each of which may be taken to be a graph over an unbounded annulus in $\Bbb{R}^2$ (see \cite{HoffmanMeeks} and \cite{Weber}). For $0\leq k\leq g$, this surface is invariant under reflection in the vertical plane making an angle of $k\pi/(g+1)$ with the $x$-axis at the origin. We call the group of symmetries of $\Bbb{R}^3$ generated by these reflections the group of {\bf horizontal symmetries} of $C_g$.{\centre \footnote{${}^\dagger$}{{\tablefont Hoffman \& Meeks show that the complete symmetry group of ${\scriptstyle C_g}$ is the dihedral group generated by the elements ${\scriptstyle A}$ and ${\scriptstyle B}$, where ${\scriptstyle A}$ is reflection in the ${\scriptstyle (x-z)}$-plane and ${\scriptstyle B}$ is rotation by an angle of ${\scriptstyle k\pi/(g+1)}$ about the ${\scriptstyle z}$-axis followed by reflection in the ${\scriptstyle (x-y)}$-plane.\hfil}}} Next, the {\bf Grim paraboloid} (also known as the {\bf bowl soliton}) is defined to be the unique simply-connected MCF soliton which is symmetric under revolution about the $z$-axis. It is known (see \cite{ClutterbuckSchnurerSchulze}) that this surface is asymptotic at infinity to a vertical translate of the graph of
$$
\frac{1}{2}r^2 - \opLog(r),
$$
where $r$ here denotes the distance in $\Bbb{R}^2$ to the origin.
\proclaim{Theorem A}
\noindent For all $g\in\Bbb{N}$ and for all sufficiently small $\epsilon$, there exists a complete, embedded MCF soliton $\Sigma_{g,\epsilon}$ of genus $g$ with $3$ ends. Furthermore, denoting $R:=\epsilon^{-1/\lambda}$, for some $\lambda\in]4,5[$, we may suppose that
\medskip
\myitem{(1)} $\Sigma_{g,\epsilon}\setminus(B(\epsilon R)\times\Bbb{R})$ consists of three connected components, each of which converges towards the same Grim paraboloid as $\epsilon$ tends to zero;
\medskip
\myitem{(2)} Upon rescaling by a factor of $1/\epsilon$, $\Sigma_{g,\epsilon}\minter(B(\epsilon R)\times\Bbb{R})$ converges to $C_g$ as $\epsilon$ tends to zero; and
\medskip
\myitem{(3)} $\Sigma_{g,\epsilon}$ is invariant under the group of horizontal symmetries of $C_g$.
\endproclaim
\remark Theorem $A$ follows immediately from Theorems \procref{ExistenceTheorem} and \procref{EmbeddednessTheorem}, below.
\medskip
\remark All notation and terminology used in this paper is explained in detail in Appendix \headref{TerminologyConventionsEtc}. Recall, in particular, that, by elliptic regularity, all standard modes of convergence of smooth, embedded solutions of parametric elliptic functionals to smooth, embedded solutions of the same functionals are equivalent over any compact region.
\medskip
\remark The constants that appear in Theorem A have the following geometric significance. The quantity $\epsilon$ determines the scaling factor of the CHM surface. Roughly speaking, it is the ``neck radius'' of $\Sigma_{g,\epsilon}$. The quantity $R$ determines how far along the end of the CHM surface the surgery is carried out. For the construction to work, $\epsilon$ and $R$ must converge in tandem to zero and infinity respectively, hence the condition $R^\lambda\epsilon=1$. Distinct values of $\epsilon$ ought to yield distinct surfaces. Indeed, a refinement of our result ought to yield a continuous family $(\Sigma_{g,\epsilon})_{\epsilon<r_0}$ of distinct embedded MCF solitons with neck radii converging to zero. Since our current argument relies on the Schauder fixed point theorem, it does not yield such fine control over the surfaces constructed.
\newsubhead{Techniques}[Techniques]
The proof of Theorem A follows the standard desingularization construction for minimal surfaces originally described by Kapouleas in \cite{KapouleasI}, \cite{KapouleasII} and \cite{KapouleasIII}. In simple terms, we first use surgery to construct an approximate MCF soliton $\hat{\Sigma}_{g,\epsilon}$ and then apply a fixed-point argument to prove the existence of an actual MCF soliton lying nearby in some suitable function space. As in all singular perturbation constructions, this is much easier said than done, and the main challenge lies in deriving the many non-trivial analytic estimates required to make the perturbation argument work.
\par
The use of CHM surfaces in Kapouleas' construction presents particular difficulties on account of their low orders of rotational symmetry. Indeed, rotational symmetries often serve in such constructions to improve decay rates and thus in turn provide stronger estimates. This phenomenon is well illustrated by the case of bounded solutions $u:S^1\times[0,\infty[\rightarrow\Bbb{R}$ of Laplace's equation $\Delta u=0$. By separation of variables, all such solutions have the form
$$
u(\theta,t) = \sum_{m\in\Bbb{Z}} a_m e^{im\theta}e^{-\left|m\right|t}.
$$
In particular, when $u$ has $n$'th order rotational symmetry, all of its coefficients of order $0<\left|m\right|<n$ vanish, so that the difference $(u-a_0)$ decays like $e^{-nt}$. Since this argument does not apply when CHM surfaces are used, we obtain our estimates by introducing instead, in Section \subheadref{TheRegularComponent}, what we call the {\bf hybrid norm}. This functional norm, which is a combination of H\"older and Sobolev norms, encapsulates the singular nature of our construction as $\epsilon$ tends to zero. Its main properties, which follow readily from the Sobolev embedding theorem and classical interpolation inequalities, are established in Lemma \procref{HybridProperty} and play a significant role in the derivation of key estimates throughout the rest of the paper.
\par
Finally, before reviewing our argument, it is worth highlighting an ingenious variant of the desingularization construction for CHM surfaces used by Hauswirth \& Pacard in \cite{HauswirthPacard}, Mazzeo \& Pacard in \cite{MazzeoPacard} and Morabito in \cite{Morabito}. In each of these papers, it is observed that the Jacobi operator $\hat{J}_{g,\epsilon}$ of the approximate minimal surface $\hat{\Sigma}_{g,\epsilon}$ is, modulo a conformal transformation when necessary, {\sl intrinsically} close to the Jacobi operator of the original CHM surface. A direct perturbation argument then yields a-priori estimates for the norm of its Green's operator, thereby bypassing one of the main technical challenges of the perturbation part of the construction. In addition, in these works, the initial surgery is carried out in a different manner than in \cite{KapouleasI}, \cite{KapouleasII} and \cite{KapouleasIII}, more pleasing to the geometric eye, though it is not clear to us whether this actually leads to simpler estimates. Regardless, their perturbation argument cannot be applied in the present case since the Jacobi operator of the joined surface is not intrinsically of the correct type.
\par
The proof is organised as follows.
\medskip
{\bf\noindent A\ -\ Rotationally Symmetric Grim Surfaces}\ We will desingularise the join of a CHM surface with three rotationally symmetric Grim ends, that is, unbounded, rotationally symmetric MCF solitons in $\Bbb{R}^3$. The geometry of CHM surfaces is well understood (see, for example, \cite{HoffmanMeeks} and \cite{Weber}) and the large-scale geometry of rotationally symmetric Grim ends has been studied by Clutterbuck, Schn\"urer \& Schulze in \cite{ClutterbuckSchnurerSchulze}. In Chapter \headref{GrimSurfaces}, we address the small-scale geometry of rotationally symmetric Grim ends, which has not previously been studied in the literature.
\par
Rotationally symmetric Grim ends exhibit a dual nature over the region of interest to us. Indeed, they are roughly catenoidal towards the lower end of this region, and roughly parabolic towards its upper end. This presents us with our first main challenge, which we address via the following algebraic trick. We introduce two abstract variables which respectively represent the catenoidal part and the parabolic part of the Grim catenoidal end. We then construct formal solutions to the MCFS equation in terms of these variables and the desired formulae are then obtained upon applying the contraction-mapping theorem to their partial sums.
\par
The main results of this chapter are Theorems \procref{AsymptoticSolutionAtInfinity} and \procref{ErrorForExactSolutionWithGivenLogarithmicParameter} which provide asymptotic formulae for the profiles of rotationally symmetric Grim ends over the large and small scales respectively.
\medskip
{\bf\noindent B\ -\ Green's Operators}\ Our perturbation argument requires estimates for the norm of a Green's operator of the MCFS Jacobi operator of the approximate MCF soliton. These are in turn derived from estimates of the norms of the corresponding operators of CHM surfaces and rotationally symmetric Grim ends. Green's operators of Jacobi operators of CHM surfaces are well understood (see, for example, \cite{HauswirthPacard}, \cite{Morabito}, \cite{Nayatani} and \cite{Pacard}). In Chapters \headref{TheGrimParaboloid} and \headref{GrimEnds}, we study the Green's operators of the MCFS Jacobi operators, first of Grim paraboloids, then of rotationally symmetric Grim ends. The former is relatively straightforward, but the latter presents us with our second main challenge, for the simple reason that catenoids become singular as the neck radius tends to zero. Let us call this the {\bf vanishing neck problem}.
\par
We address the vanishing neck problem by making, in Section \subheadref{TheModifiedMCFSJacobiOperator}, two modifications to the Jacobi operator. First, we introduce the {\bf modified MCFS Jacobi operator}, which measures the first-order variation of mean curvature arising from first-order perturbations of the surface in the direction of a suitable modification of the unit normal vector field. We underline that, since different modifications are made on different scales, the precise definition of this operator varies according to context (the description of the general framework, which unifies these definitions, is deferred to Section \subheadref{ModifiedJacobiOperators}). Next, we introduce {\bf canonical extensions} of operators across the region within the neck, which allow the modified MCFS Jacobi operators of distinct rotationally symmetric Grim ends to be compared as if they were all operators defined over $\Bbb{R}^2$.
\par
The main result of these two chapters is Theorem \procref{JIsLinearIsomorphismOverEnds}, which provides the required uniform estimates for the Green's operators of the modified MCFS Jacobi operators of rotationally symmetric Grim ends. We prove this result using a perturbation argument. We first decompose the differences between the modified MCFS Jacobi operators of Grim paraboloids and those of rotationally symmetric Grim ends into regular and singular components. Then, in Section \subheadref{TheRegularComponent}, we show that the operator norms of the regular components tend to zero as $\epsilon$ tends to zero, and in Section \subheadref{TheSingularPart}, making use of the hybrid norm described above, we show the same result for the singular component.
\medskip
{\bf\noindent C\ -\ Surgery and the Deformation Family} In Chapter \headref{SurgeryAndThePerturbationFamily}, we describe the surgery operation used to construct the approximate MCF soliton $\hat{\Sigma}_{g,\epsilon}$ as well as the deformation family around this surface within which the actual MCF soliton $\Sigma_{g,\epsilon}$ will be found. The surgery operation is straightforward and is described in Section \subheadref{TheBasicSurgeryOperation}. Simply put, the ends of the CHM surface are amputated, suitably chosen rescaled rotationally symmetric Grim ends are grafted in their place, and the join is smoothed out using cut-off functions. The construction of the deformation family about $\hat{\Sigma}_{g,\epsilon}$ is more technical and is carried out in Section \subheadref{TheDeformationFamily}. The challenge in understanding (and explaining!) this construction arises from the fact that $4$ different families of deformations must be considered. The first concerns deformations in the direction of a suitable modification of the unit normal vector field. We refer to the resulting first-order perturbations of the surface as {\bf microscopic perturbations}, since they decay at infinity. The remaining three families involve variations of the logarithmic parameters of the ends, starting far inside the locus of surgery, and vertical translations of the ends, starting far inside and far outside the locus of surgery respectively. We refer to the resulting first-order perturbations as {\bf macroscopic perturbations}, since they remain large at infinity.
\par
We associate to each of the four classes of perturbation described above the operator of first-order variation of the MCFS functional about $\hat{\Sigma}_{g,\epsilon}$. We denote these operators by $\hat{J}_{\epsilon}$, $X_\epsilon$, $Y_\epsilon$ and $Z_\epsilon$ respectively. Understanding their analytic properties is key to estimating the norm of the Green's operator of the modified MCFS Jacobi operator of $\hat{\Sigma}_{g,\epsilon}$, and we conclude this chapter by deriving preliminary estimates in Sections \subheadref{MicroscopicAndMacroscopicPerturbations}, \subheadref{ModifiedJacobiOperators} and \subheadref{ControllingMacroscopicPerturbations}.
\medskip
{\bf\noindent D\ -\ Constructing the Green's operator} In Chapter \headref{ConstructingTheGreensOperator}, we construct a Green's operator of the modified MCFS Jacobi operator of $\hat{\Sigma}_{g,\epsilon}$, together with estimates of its operator norm. This chapter constitutes the hardest technical part of the paper. The determination of sufficiently strong estimates is made possible by the correct choice of functional norms over the different components of $\hat{\Sigma}_{g,\epsilon}$ and, once again, the hybrid norm described above will play here a pivotal role.
\par
The estimates for the norm of the Green's operator are obtained in Sections \subheadref{ThePingPongArgumentPartI}, \subheadref{ThePingPongArgumentPartII} and \subheadref{ThePingPongArgumentPartIII} via a classical iteration process which we call the ``ping-pong'' argument. This process, which is common to all singular perturbation constructions, involves passing successive error terms back and forth over the join region. From a conformal perspective, the join region consists of cylinders which become very long as $\epsilon$ tends to zero. It is the exponential decay of harmonic functions over long cylinders that provides the estimates we require whilst guaranteeing the convergence of the iteration process for small $\epsilon$. We believe these ideas are best illustrated by the simplest version of this construction, used in the theory of Morse homology, and described in detail by Schwartz in Section $2.5$ of \cite{Schwarz}.
\par
The first main results of this section are Theorems \procref{NormOfA} and \procref{NormOfB} which provide estimates for the norms of the operators used in the two stages of the iteration process. In addition, Theorems \procref{NormOfW}, \procref{NormOfPInLowerRegion} and \procref{NormOfPInUpperRegion} provide estimates for the norms of the different components of the Green's operator that we construct.
\medskip
{\bf\noindent E\ -\ Existence and Embeddedness} Finally, in Chapter \headref{ExistenceAndEmbeddedness} we prove Theorem A by applying the Schauder fixed point theorem to the MCFS functional about the approximate MCF soliton $\hat{\Sigma}_{g,\epsilon}$. First, we determine estimates for the MCFS functional up to second order about $\hat{\Sigma}_{g,\epsilon}$. Then, using the estimates obtained in Chapter \headref{ConstructingTheGreensOperator}, we prove existence in Theorem \procref{ExistenceTheorem}, and we prove embeddedness in Theorem \procref{EmbeddednessTheorem} using a straightforward geometric argument.
\newsubhead{Notation}[Notation]
So as not to become overwhelmed by a deluge of constants, throughout the paper we use the following notation, which we have found to be of great help. First, given two variable quantities $a$ and $b$, we will write
$$
a\lesssim b.\eqnum{\nexteqnno[AppendixControl]}
$$
to mean that there exists a constant $C$, which for the purposes of this paper will be considered universal, such that
$$
a \leq Cb.
$$
Next, given a function $f$ and a sequence of functions $(g_m)$, we will write
$$
f=\opO(g_m).\eqnum{\nexteqnno[AppendixAsymptoticControl]}
$$
to mean that there exists a sequence $(C_m)$ of constants, which for the purposes of this paper will be considered universal, such that the relation
$$
\left|D^m f\right| \leq C_mg_m
$$
holds pointwise for all $m$. The indexing variable of the sequence $(g_m)$ should be clear from the context. In certain cases, every element of this sequence may be the same. It should also be clear from the context when this occurs. All other notation and terminology is explained in detail in Appendix \headref{TerminologyConventionsEtc}.
\newsubhead{Acknowledgements}[Acknowledgements]
We are grateful to Knut Smoczyk for drawing our attention to this interesting problem. We are likewise grateful to Detang Zhou and Andrew Clarke for helpful conversations and their invaluable insights.
\par
The first version of this paper was completed in January 2015 and was circulated amongst experts in the field before publication on arXiv. At roughly the same time, \cite{DavilaDelPinoNguyen} was also published on arXiv. \cite{DavilaDelPinoNguyen} proves a weaker result using techniques that do not extend to the case studied in this paper.
\newhead{Rotationally Symmetric Grim Surfaces}[GrimSurfaces]
\newsubhead{The Large Scale}[TheLargeScale]
We define a {\bf Grim surface} to be any unit speed MCF soliton which is a graph over some domain in $\Bbb{R}^2$. We define a {\bf Grim end} to be a Grim surface which is defined over some unbounded annulus $A(a,\infty)$. These will be studied in more detail in Chapter \headref{GrimEnds}. In this section, we study rotationally symmetric Grim surfaces defined over some annulus $A(a,b)$. We first recall the general formula for such surfaces. Let $u$ be a twice differentiable function defined over some closed interval $[a,b]$ and let $\Sigma$ be the surface of revolution generated by rotating its graph about the $z$-axis. The principle curvatures of $\Sigma$ in the radial and angular directions are respectively
$$\eqalign{
c_r &= \frac{-u_{rr}}{\sqrt{1+u_r^2}^3},\ \text{and}\cr
c_\theta &= \frac{-u_r}{r\sqrt{1+u_r^2}},\cr}\eqnum{\nexteqnno[PrincipleCurvaturesOfSurfaceOfRevolution]}
$$
where $r$ here denotes the radial distance in $A(a,b)$ from the origin, and the subscript $r$ denotes differentiation with respect to this variable. The vertical component of the upward-pointing, unit normal vector over $\Sigma$ is
$$
\langle N_\Sigma,e_z\rangle = \frac{1}{\sqrt{1 + u_r^2}},\eqnum{\nexteqnno[NormalComponentOfSurfaceOfRevolution]}
$$
so that, by \eqnref{IntroMCFSEquation}, $\Sigma$ is a rotationally symmetric Grim surface whenever
$$
ru_{rr} + (u_r-r)(1+u_r^2) = 0.\eqnum{\nexteqnno[IntegratedFormulaForGrimSurface]}
$$
Solutions of this equation have no straightforward closed form. However it will be sufficient for our purposes to obtain approximations by semi-convergent - that is, asymptotic - series. We first derive an asymptotic formula which is valid as $r$ tends to infinity.
\proclaim{Theorem \nextprocno}
\noindent If $u:]a,\infty[\rightarrow\Bbb{R}$ is a solution of \eqnref{IntegratedFormulaForGrimSurface} then, as $r\rightarrow+\infty$,
$$
u = \frac{1}{2}r^2 - \opLog(r) + a + \opO\left(r^{-(k+2)}\right),\eqnum{\nexteqnno[AsymptoticSolutionAtInfinity]}
$$
for some real constant $a$.\footnote*{{\rm We refer the reader to Section \subheadref{Notation} and Appendix \headref{TerminologyConventionsEtc} for a detailed review of the notation used here and throughout the sequel.}}
\endproclaim
\proclabel{AsymptoticSolutionAtInfinity}
\noindent Theorem \procref{AsymptoticSolutionAtInfinity} follows immediately upon integrating \eqnref{SpecialCaseARLargeScale}, below. A similar result has already been obtained by Clutterbuck, Schn\"urer \& Schulze in \cite{ClutterbuckSchnurerSchulze}. However, we consider it worth deriving \eqnref{AsymptoticSolutionAtInfinity} in full, not only because we use different techniques, but also because we believe it serves as good preparation for the more subtle small-scale asymptotic estimates that will be studied in the following sections.
\par
Define the non-linear operator $\Cal{G}$ by
$$
\Cal{G}v := r v_r + (v-r)(1+v^2),\eqnum{\nexteqnno[GrimSurfaceEquation]}
$$
and observe that $v$ solves $\Cal{G}v=0$ if and only if its integral is the profile of a rotationally symmetric Grim surface. We first derive formal solutions to \eqnref{GrimSurfaceEquation}. To this end, we define a {\bf Laurent series} in the formal variable $R$ to be a formal power series of the form
$$
V := \sum_{m=-\infty}^k V_m R^m,\eqnum{\nexteqnno[AnsatzLaurentPolynomial]}
$$
where, for all $m$, $V_m$ is a real number and $k$ is some finite integer, which we henceforth call the {\bf order} of $V$. Since the operations of formal multiplication and formal differentation are well-defined over the space of Laurent series, the operator $\Cal{G}$ also has a well-defined action over this space.
\proclaim{Lemma \nextprocno}
\noindent There exists a unique Laurent series $V$ such that $\Cal{G}V=0$. Furthermore,
\medskip
\myitem{(1)} $V$ has order $1$;
\medskip
\myitem{(2)} $V_1=1$, $V_{-1}=-1$;
\medskip
\myitem{(3)} if $m$ is even, then $V_m=0$; and
\medskip
\myitem{(4)} if $\hat{V}_n:=\sum_{m=1-2n}^1V_mR^m$ denotes the $n$'th partial sum of $V$, then $\Cal{G}\hat{V}_n$ is a finite Laurent series of order $(1-2n)$.
\endproclaim
\proclabel{LemmaFormalSolution}
\proof Consider the ansatz \eqnref{AnsatzLaurentPolynomial}. If $k\leq -1$, then the highest order term in $\Cal{G}V$ is $(-R)$, if $k=0$, then it is $(-R)(1+V_0^2)$, and if $k\geq 2$, then it is $V_k^3R^{3k}$. Since none of these vanish, it follows that $V$ must be of order $1$. In this case, the highest order term in $\Cal{G}V$ is $V_1^2(V_1-1)R^3$ so that, in order to have non-trivial solutions, we require $V_1=1$. We now have
$$
R\frac{dV}{dR} + (V-R)(1+V^2) =
R + \sum_{m=-\infty}^0(m+1)V_mR^m +
\sum_{m=-\infty}^2\left(\sum_{\scriptstyle p+q+r=m,\atop{\scriptstyle p\leq 0,\atop{\scriptstyle q,r\leq 1}}}V_pV_qV_r\right)R^m.
$$
In particular, setting the respective coefficients of $R^2$ and $R$ equal to $0$ yields
$$
V_0=0,\ V_{-1}=-1.
$$
For all $m\leq -2$, setting the coefficient of $R^{m+2}$ equal to $0$ now yields
$$
V_m + \left(\sum_{\scriptstyle p+q+r=m+2,\atop{\scriptstyle m+1\leq p\leq -1,\atop{\scriptstyle m+2\leq q,r\leq 1}}}V_pV_qV_r\right) + (m+3)V_{m+2} = 0.\eqnum{\nexteqnno[RecurrenceRelationAtInfinity]}
$$
The existence and uniqueness of $V$ now follow from this recurrence relation. Furthermore, if $p+q+r=m+2$, and if $m$ is even, then at least one of $p$, $q$ and $r$ is also even, and since $V_0=0$, it follows by induction that $V_m=0$ for all even $m$. In addition, by \eqnref{RecurrenceRelationAtInfinity}, for all $n$, and for all $m\geq(3-2n)$, the coefficient of $R^m$ in $\Cal{G}\hat{V}_n$ is equal to $0$. However, since $V_{-2n}=0$, by \eqnref{RecurrenceRelationAtInfinity} again, the coefficient of $R^{2-2n}$ in $\Cal{G}\hat{V}_n$ is also equal to $0$, so that $\Cal{G}\hat{V}_n$ is a finite Laurent series of order $(1-2n)$, as desired.\qed
\medskip
For all $n$, define the $n$'th partial sum $v_n:]0,\infty[\rightarrow\Bbb{R}$ by
$$
v_n(r) := \sum_{m=1-2n}^1V_mr^m.\eqnum{\nexteqnno[PartialSumAtInfinity]}
$$
We now show that the sequence $(v_n)$ yields successively better approximations over the large scale of the exact solutions of $\Cal{G}v=0$. We first derive zero'th order bounds.
\proclaim{Lemma \nextprocno}
\noindent If $v:[a,\infty[\rightarrow\Bbb{R}$ solves \eqnref{GrimSurfaceEquation} then, for large $r$,
$$
\left|v_0 - v\right| \lesssim \frac{1}{r}.\eqnum{\nexteqnno[FirstOrderEstimateLargeScale]}
$$
\endproclaim
\proclabel{LemmaFirstOrderEstimate}
\proof Consider the family of polynomials $p_t(y)=(y-1)(t^2 + y^2)$. For all $t>0$, $y=1$ is the unique real root of $p_t$. Since $y=0$ is the unique local maximum of $p_0$, for sufficiently small $t$, the unique local maximum of $p_t$ is also near $0$, and the value of $p_t$ at this point is less than $-t^2/2$. Since $p_0$ is convex over the interval $[1/3,\infty[$, for $1/3<y<1$, $p_0(y)\leqslant (3/2)(1-y)p_0(1/3)=(y-1)/9$ and so, for sufficiently small $t$, over the smaller interval $[1/2,1]$, $p_t(y)\leqslant (y-1)/18$.
\par
Now let $v$ be a solution of $\Cal{G}v=0$. In particular, using a dot to denote differentiation with respect to $r$, we have $\dot{v}=-r^2p_{1/r}(v/r)$. Suppose, furthermore, that $r\gg 1$ so that the estimates of the preceeding paragraph hold for $p_{1/r}$. When $v\geq r$, $\dot{v}-\dot{r}=\dot{v}-1\leq-1$, so that, for sufficiently large $r$, $v(r)\leqslant r$. If $v\leqslant r/2$, then $\dot{v}-\dot{r}/2\geqslant r/2-1/2$, so that, for sufficiently large $r$, $v(r)\geqslant r/2$. Finally, if $r/2\leqslant v\leqslant r$ then, by the preceding discussion, $\dot{v}\geqslant r(r-v)/18$. It follows that if $w:=r(v_0-v)=r(r-v)$, then $w>0$ and $\dot{w} = 2r - v - r\dot{v} \leqslant r + w/r - rw/18$. Since this is negative for $w\geqslant 36$ and $r>6$, the function $w$ is bounded, and the result follows.\qed
\proclaim{Lemma \nextprocno}
\noindent If $v:[a,\infty[\rightarrow\Bbb{R}$ solves \eqnref{GrimSurfaceEquation} then, for all $n$, and for large $r$,
$$
\left|v_n - v\right| \lesssim r^{-(2n+1)}.\eqnum{\nexteqnno[ZeroethOrderAsymptoticRelationLargeScale]}
$$
\endproclaim
\proclabel{ZeroethOrderAsymptoticRelationLargeScale}
\proof For all $n$, let $w_n:=r^{2n-1}(v_n-v)$ be the rescaled error. We prove by induction that $\left|w_n\right|\lesssim r^{-2}$ for all $n$. Indeed, the case $n=0$ follows from \eqnref{FirstOrderEstimateLargeScale}. We suppose therefore that $n\geq 1$. Since $w_n=r^2 w_{n-1} + V_{1-2n}$, it follows by the inductive hypothesis that $w_n$ is bounded. Now let $P(a,b)$ denote any polynomial in the variables $a$ and $b$. Since $\Cal{G}v_{n}$ is a finite Laurent polynomial of order $(1-2n)$, using a dot to denote differentiation with respect to $r$, we have
$$\eqalign{
\dot{w}_n &= \frac{(2n-1)}{r}w_n + r^{2n-2}(r\dot{v}_n - r\dot{v})\cr
&=\frac{1}{r}P\left(\frac{1}{r},w_n\right) - r^{2n-2}\left((v_n-r)(1+v_n^2)-(v-r)(1+v^2)\right)\cr
&=\frac{1}{r}P\left(\frac{1}{r},w_n\right) - \frac{1}{r}w_n\left(1 - r(v_n+v) + (v_n^2 + v_nv + v^2)\right).\cr}
$$
Since $v=v_n-r^{-(2n-1)}w_n$ and since $(v_n-r)$ is also a polynomial in $r^{-1}$ with no constant term, this yields
$$
\dot{w}_n = \frac{1}{r}P\left(\frac{1}{r},w_n\right) - rw_n.
$$
Since $w_n$ is bounded, there therefore exists a constant $B>0$ such that for all $r\geqslant 1$,
$$
\left|\dot{w}_n+rw_n\right|\leqslant Br^{-1}.\eqnum{\nexteqnno[BoundOnDerivativeOfW]}
$$
In particular, for $r\geq 2$ and $r^2 w_n\geq 2B$,
$$
\frac{d}{dr}r^2 w_n = r^2(\dot{w}_n + rw_n) + (2r-r^3)w_n \leq Br - \frac{1}{2}r^3 w_n \leq 0,
$$
so that $r^2w_n$ is bounded from above. Since $(-w_n)$ also satisfies \eqnref{BoundOnDerivativeOfW}, we see that $r^2w_n$ is bounded from below, and this completes the proof.\qed
\proclaim{Lemma \nextprocno}
\noindent If $v:[a,\infty[\rightarrow\Bbb{R}$ solves \eqnref{GrimSurfaceEquation} then, for all $n$,
$$
v_n - v = \opO\left(r^{-(k+2n+1)}\right).\eqnum{\nexteqnno[AsymptoticRelationLargeScale]}
$$
In particular,
$$
v = r - \frac{1}{r} + \opO\left(r^{-(k+3)}\right).\eqnum{\nexteqnno[SpecialCaseARLargeScale]}
$$
\endproclaim
\proclabel{AsymptoticRelationLargeScale}
\proof For all $n$, denote $w_n:=(v_n-v)$. As in the proof of Lemma \procref{ZeroethOrderAsymptoticRelationLargeScale}, we obtain
$$
\dot{w}_n = P_1\left(\frac{1}{r},w_n\right)rw_n + \frac{1}{r}\Cal{G}v_n,
$$
where $P_1$ is some polynomial. Since $\Cal{G}v_n$ is a finite Laurent polynomial of order $(1-2n)$, it follows by induction that, for all $k$,
$$
\frac{d^kw_n}{dr^k} = P_k\left(\frac{1}{r},w_n\right)r^k w_n + Q_k\left(\frac{1}{r},w_n\right)r^{k-(2n+1)},
$$
where $P_k$ and $Q_k$ are polynomials. It follows by \eqnref{ZeroethOrderAsymptoticRelationLargeScale} that, for all $k$,
$$
\left|\frac{d^k}{dr^k}(v_n-v)\right| = \left|\frac{d^kw_n}{dr^k}\right|\lesssim r^{k-(2n+1)}.
$$
However, since $(v_{n+k}-v_n)$ is a finite Laurent series of order $-(2n+1)$, for all $k$,
$$
\left|\frac{d^k}{dr^k}(v_n-v)\right| \leq \left|\frac{d^k}{dr^k}(v_n-v_{n+k})\right| + \left|\frac{d^k}{dr^k}(v_{n+k}-v)\right| \lesssim r^{-(k+2n+1)},
$$
and the result follows.\qed
\newsubhead{The Small Scale - Formal Solutions}[TheSmallScaleFormalSolutions]
We now study solutions to \eqnref{GrimSurfaceEquation} over the small scale. We fix positive constants $K\gg 1$ and $\eta\ll 1$ which we henceforth consider to be universal. Let $\Lambda$ be a large, positive real number, and let $\epsilon, R>0$ and $c\in\Bbb{R}$ be such that
$$
\left(\epsilon R^{4+\eta} + \frac{1}{R^{1-\eta}}\right)\leq \frac{1}{\Lambda},\ \epsilon R^{5-\eta}\geq\Lambda,\ \left|c\right|\leq K.\eqnum{\nexteqnno[ControlOfParameters]}
$$
These conditions will be used repeatedly throughout the paper. Observe, in particular, that \eqnref{ControlOfParameters} implies that $\epsilon$ becomes small and $R$ becomes large as $\Lambda$ tends to infinity. We will prove
\proclaim{Theorem \nextprocno}
\noindent For all sufficiently large $\Lambda$, and for all $R,\epsilon$ satisfying \eqnref{ControlOfParameters}, there exists a smooth function $\sigma[\epsilon,R]:\Bbb{R}\rightarrow\Bbb{R}$ such that, for all $c\in[-K,K]$, if $v:[\epsilon R,\epsilon R^4]\rightarrow\Bbb{R}$ solves $\Cal{G}v=0$ with initial value
$$
v(\epsilon R) = \frac{1}{R}\sigma[\epsilon,R](c) + \frac{\epsilon R}{2},\eqnum{\nexteqnno[DependenceOfICOnLogCoeff]}
$$
then
$$
v(r) = \frac{1}{2}r + \frac{c\epsilon}{r} + \opO\left(\left[1 + \opLog\left(\frac{r}{\epsilon R}\right)\right]\frac{1}{r^k}\left(r + \frac{\epsilon}{r}\right)^3\right),
\eqnum{\nexteqnno[ErrorForExactSolutionWithGivenLogarithmicParameter]}
$$
Furthermore, the function $\sigma[\epsilon,R]$ converges to the identity in the $C^\infty_\oploc$ sense as $\Lambda$ tends to $+\infty$.
\endproclaim
\proclabel{ErrorForExactSolutionWithGivenLogarithmicParameter}
\noindent The function $\sigma[\epsilon,R]$ will be defined explicitly in Section \subheadref{TheSmallScaleExactSolutionsI}, below, and Theorem \procref{ErrorForExactSolutionWithGivenLogarithmicParameter} will follow immediately from Lemma \procref{EstimateOfErrorInLogCoordinates}, below. The constant $c$ will henceforth be referred to as the {\bf logarithmic parameter} of the function $v$. Observe that, up to a small perturbation, it is related to the initial value of $v$ by a linear function. This perturbation is required in order to guarantee good estimates over the whole interval. Indeed, replacing $\sigma[\epsilon,R](c)$ by $c$ in \eqnref{DependenceOfICOnLogCoeff} would increase the error in \eqnref{ErrorForExactSolutionWithGivenLogarithmicParameter}, making it then of order $(\epsilon/r)$.
\par
In order to appreciate Theorem \procref{ErrorForExactSolutionWithGivenLogarithmicParameter} and the argument that follows, we find it helpful to first recall the geometric properties of the function $v$ over the interval $[\epsilon R,\epsilon R^4]$. Indeed, by definition, its integral $u$ is the profile of some rotationally symmetric Grim surface. However, it is known (see \cite{ClutterbuckSchnurerSchulze}) that, near the lower end of this interval, the first term in the MCFS equation \eqnref{IntroMCFSEquation} dominates, so that the graph of $u$ is close to some minimal catenoid in $\Bbb{R}^3$ and the function $u$ is itself approximately logarithmic. On the other hand, near the upper end of this interval, it is the second term in the MCFS equation which dominates, and the function $u$ is approximately quadratic, in accordance with the asympototic formula obtained in the preceding section. These two contrasting behaviours are reflected in \eqnref{ErrorForExactSolutionWithGivenLogarithmicParameter} by the $\epsilon/r$ terms and the $r$ terms respectively.
\par
In order to derive an asymptotic formula for $u$ that simultaneously describes these two behaviours, we introduce two abstract variables $M$ and $N$, where $M$ measures its quadratic behaviour, and $N$ measures its logarithmic behaviour. By expressing the equation $\Cal{G}v=0$ in terms of these new variables, the asymptotic formula for $v$ is then obtained in the same manner as in Section \subheadref{TheLargeScale} namely, by first determining formal solutions which then serve as approximations for exact solutions.
\par
Upon applying the change of variables $r:=\epsilon R e^x$ we obtain
$$
\Cal{G}v = \Cal{D}v - \epsilon R e^x + (v - \epsilon R e^x)v^2,\eqnum{\nexteqnno[ReparametrisedGrimSurfaceEquation]}
$$
where the operator $\Cal{D}$ is defined by
$$
\Cal{D}v := v_x + v,\eqnum{\nexteqnno[DefinitionOfD]}
$$
and the subscript $x$ here denotes differentiation with respect to this variable. Now let $\Bbb{R}[X,M,N]$ be the ring of polynomials with real coefficients in the variables $X$, $M$ and $N$. We consider a general element $V$ of $\Bbb{R}[X,M,N]$ as a sum of the form
$$
V = \sum_{p,q\leq k} V_{p,q}(X)M^pN^q,\eqnum{\nexteqnno[FormalSum]}
$$
where, for all $p,q$, $V_{p,q}$ is a polynomial in the variable $X$ and $k$ is some finite, non-negative integer which we henceforth refer to as the {\bf order} of $V$. There is a natural correspondence sending $\Bbb{R}[X,M,N]$ into the space of continuous functions over $[0,3\opLog(R)]$ given by
$$
V \mapsto v(x) := \sum_{p,q\leq k}V_{p,q}(x)\left(\epsilon R e^x\right)^p\left(\frac{c}{R}e^{-x}\right)^q.\eqnum{\nexteqnno[FormalToFunctionCorrespondence]}
$$
In other words, this correspondence is the unique $\Bbb{R}[X]$-ring homomorphism which sends $M$ to $\epsilon R e^x$ and $N$ to $\frac{c}{R}e^{-x}$. Although this homomorphism is not injective, it keeps track of the parameters $\epsilon$, $R$ and $c$, which is the reason why it serves our purposes. Operators $\Cal{G}$ and $\Cal{D}$ are also defined over $\Bbb{R}[X,M,N]$ by
$$\eqalign{
\Cal{G}V &:= \Cal{D}V - M + (V-M)V^2,\ \text{and}\cr
(\Cal{D}V)_{p,q} &:= \left(\frac{d}{dX} + 1 + (p-q)\right)V_{p,q},\cr}\eqnum{\nexteqnno[FormalOperators]}
$$
where $\frac{d}{dX}$ here denotes the operator of formal differentiation with respect to the variable $X$. In particular, $\Cal{G}$ and $\Cal{D}$ both map through the above correspondence to the operators given in \eqnref{ReparametrisedGrimSurfaceEquation} and \eqnref{DefinitionOfD} respectively, thereby justifying this notation. Observe, furthermore, that $\Cal{D}$ defines a surjective linear map from $\Bbb{R}[X,M,N]$ to itself and that its kernel consists of finite sums of the form
$$
V = \sum_{p\leq k} a_p M^p N^{p+1},
$$
where $a_0,...,a_k$ are real constants.
\par
Let $\Bbb{R}[X][[M,N]]$ be the ring of formal power series over the variables $M$ and $N$ with coefficients that are polynomials in the variable $X$. Observe that the operators $\Cal{G}$ and $\Cal{D}$ naturally extend again to well-defined operators over this space.
\proclaim{Lemma \nextprocno}
\noindent There exists a unique formal power series $V$ in $\Bbb{R}[X][[M,N]]$ such that
\medskip
\myitem{(1)} $V_{0,1}=1$;
\medskip
\myitem{(2)} $V_{p,p+1}(0)=0$ for all $p\geq 1$; and
\medskip
\myitem{(3)} $\Cal{G}V=0$.
\medskip
\noindent Furthermore,
\medskip
\myitem{(4)} $V_{1,0}=\frac{1}{2}$;
\medskip
\myitem{(5)} if $p+q$ is even, then $V_{p,q}=0$; and
\medskip
\myitem{(6)} if $p+q=2k+1$ is odd, then $V_{p,q}$ has order at most $k$ in $X$.
\medskip
\noindent Finally, if we denote
$$
\hat{V}_k := \sum_{p+q\leq 2k+1}V_{p,q}(X)M^pN^q,
$$
\noindent then,
\medskip
\myitem{(7)} if $(p+q)\leq(2k+1)$, then the coefficient of $M^pN^q$ in $\Cal{G}\hat{V}_k$ vanishes;
\medskip
\myitem{(8)} if $(p+q)>(2k+1)$ is even, then the coefficient of $M^pN^q$ in $\Cal{G}\hat{V}_k$ vanishes; and
\medskip
\myitem{(9)} if $(p+q)>(2k+1)$ is odd, then the coefficient of $M^pN^q$ in $\Cal{G}\hat{V}_k$ has order at most $(p+q-3)/2$ in $X$.
\endproclaim
\proclabel{LemmaFormalSolutionOverSmallScale}
\proof Let $V=\sum_{p,q} V_{p,q}(X)M^pN^q$ be an element of $\Bbb{R}[X][[M,N]]$ which solves $\Cal{G}V=0$. For all $(p,q)$, equating the coefficient of $M^pN^q$ in $\Cal{G}V$ to $0$, we obtain
$$\eqalign{
\left(\frac{d}{dX} + (1 + (p-q))\right)V_{p,q}
&=\delta_{p1}\delta_{q0} -
\sum_{\scriptstyle p_1+p_2+p_3=p,\atop{\scriptstyle q_1+q_2+q_3=q.}}
V_{p_1,q_1}V_{p_2,q_2}V_{p_3,q_3}\cr
&\qquad +
\sum_{\scriptstyle p_1+p_2=p-1,\atop{\scriptstyle q_1+q_2=q.}}
V_{p_1,q_1}V_{p_2,q_2}.\cr}\eqnum{\nexteqnno[EqnRecurrenceForSmallScaleFormalSolutions]}
$$
In particular,
$$
\frac{dV_{0,0}}{dX} + V_{0,0}(1 + V_{0,0}^2) = 0,
$$
and since there exists no non-trivial polynomial solution to this equation, it follows that that $V_{0,0}=0$. From this it follows that the two summations on the right hand side of \eqnref{EqnRecurrenceForSmallScaleFormalSolutions} only involve terms of order at most $p+q-2$ in $(M,N)$. In particular, $V_{0,1}$ satisfies
$$
\frac{dV_{0,1}}{dX} = 0.
$$
It is thus constant, and we henceforth set it equal to $1$. It now follows by induction that there exists a unique sequence of polynomials $(V_{p,q})$ satisfying \eqnref{EqnRecurrenceForSmallScaleFormalSolutions} such that $V_{0,1}=1$ and $V_{p,p+1}(0)=0$ for all $p\geq 1$.
\par
To prove $(4)$, observe that $V_{1,0}$ satisfies $\frac{dV_{1,0}}{dX} + 2V_{1,0} = 1$ so that, since it is a polynomial, $V_{1,0}=\frac{1}{2}$, as desired. To prove $(5)$, observe that if $p+q$ is even, then every summand on the right hand side in \eqnref{EqnRecurrenceForSmallScaleFormalSolutions} involves at least one term of the form $V_{p',q'}$, where $p'+q'$ is an even number no greater than $p+q-2$. Since $V_{0,0}=0$, it follows by induction that $V_{p,q}=0$ whenever $p+q$ is even, as desired. To prove $(6)$, suppose that for all $l<k$, and for $p+q=2l+1$, the polynomial $V_{p,q}$ has order at most $l$ in $X$. By \eqnref{EqnRecurrenceForSmallScaleFormalSolutions}, for all $p+q=2k+1$, the polynomial $V_{p,q}$ is obtained by integrating terms of order at most $(k-1)$ in $X$, and it follows by induction that $V_{p,q}$ has order at most $k$ in $X$, as desired.
\par
Finally, observe that, by \eqnref{EqnRecurrenceForSmallScaleFormalSolutions}, the term $V_{p,q}$ is defined by setting the coefficient of $M^pN^q$ equal to zero in $\Cal{G}V$, and $(7)$ follows. Furthermore, for $p+q>(2k+1)$, the coefficient of $M^pN^q$ in $\Cal{G}V$ is equal to the right hand side of \eqnref{EqnRecurrenceForSmallScaleFormalSolutions}. $(8)$ and $(9)$ now follow by similar arguments used to prove $(5)$ and $(6)$, above, and this completes the proof.\qed
\newsubhead{The Small Scale - Exact Solutions I}[TheSmallScaleExactSolutionsI]
Let $V$ be the formal power series constructed in Lemma \procref{LemmaFormalSolutionOverSmallScale}. For $\epsilon, R$ satisfying \eqnref{ControlOfParameters}, for $c\in\Bbb{R}$, and for non-negative, integer $k$, let $v_{k,c}$ be the $k$'th partial sum of $V$ with logarithmic parameter $c$, that is
$$
v_{k,c}(x):=\sum_{p+q\leq 2k+1}V_{p,q}(x)\left(\epsilon R e^x\right)^p\left(\frac{c}{R}e^{-x}\right)^q.\eqnum{\nexteqnno[FunctionOfPartialSum]}
$$
Define the function $\sigma[\epsilon,R,k]:\Bbb{R}\rightarrow\Bbb{R}$ by
$$
\sigma[\epsilon,R,k](c) := Rv_{k,c}(0) - \frac{\epsilon R^2}{2}.\eqnum{\nexteqnno[DefinitionOfSigmaK]}
$$
Trivially, if $v:[0,3\opLog(R)]\rightarrow\Bbb{R}$ satisfies
$$
v(0) = \frac{1}{R}\sigma[\epsilon,R,k](c) + \frac{\epsilon R}{2},
$$
then $v$ has the same initial value as $v_{k,c}$. Observe that $\sigma[\epsilon,R,k]$ is a polynomial in $c$ with coefficients that depend on $\epsilon$, $R$ and $k$ and, for all $k$, $\sigma[\epsilon,R,k]$ converges to the identity in the $C^\infty_\oploc$ sense as $\Lambda$ tends to infinity. We will see presently that the estimates we require follow when $k$ is at least $9$, and we therefore define
$$
\sigma[\epsilon,R](c) := \sigma[\epsilon,R,9](c).\eqnum{\nexteqnno[DefinitionOfSigma]}
$$
This is the function that appears in the statement of Theorem \procref{ErrorForExactSolutionWithGivenLogarithmicParameter}.
\par
As in Section \subheadref{TheLargeScale}, we now determine zero'th order bounds for the difference between $v_{k,c}$ and an exact solution with the same initial value. We achieve this via the contraction mapping theorem. We first introduce the required analytic framework. For $T\in [0,3\opLog(R)]$, let $C^0([0,T])$ be the Banach space of continuous functions over the interval $[0,T]$ furnished with the uniform norm and let $C^1_0([0,T])$ be the Banach space of continuously differentiable functions over this interval with initial value $0$, furnished with the norm
$$
\|w\|_{C_0^1} := \|w_x\|_{C^0},\eqnum{\nexteqnno[FirstOrderNorm]}
$$
where the subscript $x$ here denotes differentiation with respect to this variable. Observe that, for all $w\in C^1_0([0,T])$,
$$
\|w\|_{C^0} \leq T\|w\|_{C_0^1}.\eqnum{\nexteqnno[ComparisonOfBothNorms]}
$$
\proclaim{Lemma \nextprocno}
\noindent The operator $\Cal{D}$ defines a linear isomorphism from $C^1_0([0,T])$ into $C^0([0,T])$. Furthermore, the operator norms of $\Cal{D}$ and its inverse satisfy
$$
\|\Cal{D}\|\leq 1 + T,\ \|\Cal{D}^{-1}\| \leq 2.\eqnum{\nexteqnno[NormOfDAndItsInverse]}
$$
\endproclaim
\proof First, bearing in mind \eqnref{ComparisonOfBothNorms},
$$
\|\Cal{D}w\|_{C^0} \leq \|w_x\|_{C^0} + \|w\|_{C^0} \leq (1+T)\|w\|_{C_0^1},
$$
so that $\|\Cal{D}\|\leq 1+T$. By inspection, for all $w$,
$$
(\Cal{D}^{-1}w)(x) = e^{-x}\int_0^x e^y w(y)dy.
$$
In particular,
$$
\|\Cal{D}^{-1}w\|_{C^0} \leq \|w\|_{C^0}.
$$
Thus,
$$
\|\Cal{D}^{-1}w\|_{C_0^1} = \|(\Cal{D}^{-1}w)_x\|_{C^0} \leq \|\Cal{D}\Cal{D}^{-1}w\|_{C^0} + \|\Cal{D}^{-1}w\|_{C^0} \leq 2\|w\|_{C^0},
$$
so that $\|\Cal{D}^{-1}\|\leq 2$. This completes the proof.\qed
\medskip
\noindent Consider now the functional $\Cal{H}:C^1_0([0,T])\rightarrow C^0([0,T])$ given by
$$
\Cal{H}(w) := \Cal{G}(v_k + w).\eqnum{\nexteqnno[DefinitionOfErrorOperator]}
$$
Its Frechet derivative at $w$ is
$$
D\Cal{H}(w)f := \Cal{D}f + \Cal{E}(w)f,\eqnum{\nexteqnno[DerivativeOfErrorOperator]}
$$
where
$$
\Cal{E}(w)f := 3(v_{k,c} + w)^2 f - 2\epsilon R e^x(v_{k,c} + w)f.\eqnum{\nexteqnno[FormulaForErrorOperator]}
$$
\proclaim{Lemma \nextprocno}
\noindent For all $w\in C^1_0([0,T])$, the operator norm of $\Cal{E}(w)$, considered as a linear map from $C^1_0([0,T])$ into $C^0([0,T])$, satisfies
$$
\|\Cal{E}(w)\| \lesssim T\left((\epsilon R e^T)^2 + \frac{1}{R^2} + T^2\|w\|_{C_0^1}^2\right).\eqnum{\nexteqnno[SizeOfErrorOperator]}
$$
\endproclaim
\proof Indeed, over $[0,T]$,
$$
\|\epsilon R e^x\|_{C^0} \leq \epsilon R e^T,\
\left\|\frac{c}{R}e^{-x}\right\|_{C^0} \leq \frac{c}{R}.
$$
Thus, by Lemma \procref{LemmaFormalSolutionOverSmallScale} and \eqnref{ControlOfParameters},
$$
\|v_k\|_{C^0} \lesssim \sum_{i=0}^k (1+T^i)\left(\epsilon R e^T + \frac{1}{R}\right)^{2i+1} \lesssim \epsilon R e^T + \frac{1}{R},
$$
so that, by \eqnref{ComparisonOfBothNorms} and \eqnref{FormulaForErrorOperator},
$$\eqalign{
\|\Cal{E}(w)f\|_{C^0} &\lesssim \left((\epsilon R e^T)^2 + \frac{1}{R^2} + \|w\|_{C^0}^2\right)\|f\|_{C^0}\cr
&\lesssim T\left((\epsilon R e^T)^2 + \frac{1}{R^2} + T^2\|w\|_{C_0^1}^2\right)\|f\|_{C_0^1},\cr}
$$
as desired.\qed
\medskip
\noindent Define the map $\Phi:C^1_0([0,T])\rightarrow C^1_0([0,T])$ by
$$
\Phi(w) := w - \Cal{D}^{-1}\Cal{H}(w).\eqnum{\nexteqnno[DefinitionOfContractionMapping]}
$$
\proclaim{Lemma \nextprocno}
\noindent For $w,\overline{w}\in C^1_0([0,T])$,
$$
\|\Phi(w) - \Phi(\overline{w})\|_{C_0^1} \lesssim T\left((\epsilon R e^T)^2 + \frac{1}{R^2} + T^2\|w\|_{C_0^1}^2 + T^2\|\overline{w}\|_{C_0^1}^2\right)\|w-\overline{w}\|_{C_0^1}.\eqnum{\nexteqnno[ContractionMappingInequality]}
$$
\endproclaim
\proof Indeed, for $w,\overline{w}\in C^1_0([0,T])$, using \eqnref{DerivativeOfErrorOperator},
$$\eqalign{
\Phi(w) - \Phi(\overline{w}) &= w-\overline{w} - \Cal{D}^{-1}\left(\Cal{H}(w)-\Cal{H}(\overline{w})\right)\cr
&=-\Cal{D}^{-1}\left(\Cal{H}(w)-\Cal{H}(\overline{w}) - \Cal{D}(w-\overline{w})\right)\cr
&=-\Cal{D}^{-1}\left(\int_0^1 \Cal{E}(tw + (1-t)\overline{w})dt\right)(w-\overline{w}).\cr}
$$
Thus, by \eqnref{NormOfDAndItsInverse} and \eqnref{SizeOfErrorOperator},
$$
\|\Phi(w)-\Phi(\overline{w})\|_{C_0^1} \lesssim T\left((\epsilon R e^T)^2 + \frac{1}{R^2} + T^2\|w\|_{C_0^1}^2 + T^2\|\overline{w}\|_{C_0^1}^2\right)\|w-\overline{w}\|_{C_0^1},
$$
as desired.\qed
\medskip
\noindent Applying the contraction-mapping theorem now yields
\proclaim{Lemma \nextprocno}
\noindent For sufficiently large $\Lambda$, if $v_{k,c}$ is the $k$'th partial sum of $V$ with logarithmic parameter $c$, and if $v:[0,3\opLog(R)]\rightarrow\Bbb{R}$ solves $\Cal{G}v=0$ with initial value
$$
v(0) = \frac{1}{R}\sigma[\epsilon,R,k](c) + \frac{\epsilon R}{2},\eqnum{\nexteqnno[InitialConditions]}
$$
then
$$
\|v - v_{k,c}\|_{C^0} \lesssim (1+T^{k+1})\left(\epsilon R e^T + \frac{1}{R}\right)^{2k+3}.\eqnum{\nexteqnno[GlobalEstimateOfError]}
$$
\endproclaim
\proclabel{LemmaGlobalEstimateOfError}
\proof By Lemma \procref{LemmaFormalSolutionOverSmallScale},
$$
\|\Cal{G}v_k\|_{C^0} \lesssim \left(1+T^k\right)\left(\epsilon R e^T + \frac{1}{R}\right)^{2k+3}.
$$
By \eqnref{NormOfDAndItsInverse}, there therefore exists $B>0$, which we may consider to be universal, such that
$$
\|\Phi(0)\|_{C_0^1} = \|\Cal{D}^{-1}\Cal{G}v_k\|_{C_0^1} \leq B(1+T^k)\left(\epsilon R e^T + \frac{1}{R}\right)^{2k+3}.
$$
Let $X$ be the closed ball of radius $2B(1+T^k)\left(\epsilon R e^T + R^{-1}\right)^{2k+3}$ about $0$ in $C^1_0([0,T])$. By \eqnref{ControlOfParameters}, if $w,\overline{w}\in X$ then, in particular,
$$
T\|w\|_{C_0^1},T\|\overline{w}\|_{C_0^1} \lesssim \left(\epsilon R e^T + \frac{1}{R}\right),
$$
so that, by \eqnref{ContractionMappingInequality} and \eqnref{ControlOfParameters} again,
$$
\|\Phi(w) - \Phi(\overline{w})\|_{C_0^1} \lesssim \frac{1}{\Lambda}\|w-\overline{w}\|_{C_0^1}.
$$
The map $\Phi$ thus defines a contraction from $X$ to itself, and there therefore exists $w\in X$ such that $\Phi(w)=w$. In particular $\Cal{H}w=0$, and
$$
\|w\|_{C^0} \leq T\|w\|_{C^1_0} \lesssim (1+T^{k+1})\left(\epsilon R e^T + \frac{1}{R}\right)^{2k+3}.
$$
Finally, by definition of the function $\sigma[\epsilon,R,k]$, $v(0)=v_{k,c}(0)$ so that, by uniqueness of solutions to ODEs with prescribed initial values, $v-v_{k,c}=w$, and the result follows.\qed
\newsubhead{The Small Scale - Exact Solutions II}[TheSmallScaleExactSolutionsII]
The final step in proving Theorem \procref{ErrorForExactSolutionWithGivenLogarithmicParameter} involves extending the estimates obtained in Lemma \procref{LemmaGlobalEstimateOfError} to derivatives of all orders.
\proclaim{Lemma \nextprocno}
\noindent If $v_{k,c}$ and $v$ are as in Lemma \procref{LemmaGlobalEstimateOfError}, then
$$
v = v_{k,c} + \opO\left((1+T^{k+1})\left(\epsilon R e^T + \frac{1}{R}\right)^{2k+3}\right).\eqnum{\nexteqnno[DifferenceBetweenFormalAndExactToAllOrders]}
$$
\endproclaim
\proof Denote $w:=v-v_{k,c}$. Since $v_{k,c}$ is a polynomial in $\epsilon R e^x$ and $\frac{c}{R}e^{-x}$ with coefficients in $\Bbb{R}[X]$, as in the proof of Lemma \procref{ZeroethOrderAsymptoticRelationLargeScale},
$$
w_x = P_1\left(w,\epsilon R e^x,\frac{c}{R}e^{-x}\right)w + \Cal{G}v_{k,c},
$$
for some polynomial $P_1$ with coefficients in $\Bbb{R}[X]$. Since $\Cal{G}v_{k,c}$ is also a polynomial in $\epsilon R e^x$ and $\frac{c}{R}e^{-x}$ with coefficients in $\Bbb{R}[X]$, it follows by induction that for all $l$,
$$
\frac{d^l}{dx^l}w = P_l\left(w,\epsilon R e^x,\frac{c}{R}e^{-x}\right)w + \sum_{p=0}^{l-1}Q_{p,l}\left(w,\epsilon R e^x,\frac{c}{R}\right)\frac{d^p}{dx^p}\Cal{G}v_{k,c},\eqnum{\nexteqnno[HigherDerivativesOfDifference]}
$$
for suitable polynomials $P_l$ and $(Q_{p,l})_{0\leq p\leq l-1}$ also with coefficients in $\Bbb{R}[X]$. However, by \eqnref{GlobalEstimateOfError},
$$
\|w\|_{C^0}\lesssim (1+T^{k+1})\left(\epsilon R e^T + \frac{1}{R}\right)^{2k+3}.
$$
Thus, by \eqnref{ControlOfParameters},
$$
\left\|P_l\left(w,\epsilon R e^x, \frac{c}{R}e^{-x}\right)\right\|_{C^0}, \left\|Q_{p,l}\left(\epsilon R e^x,\frac{c}{R}e^{-x}\right)\right\|_{C^0} \lesssim 1.
$$
Finally, Lemma \procref{LemmaFormalSolutionOverSmallScale} and \eqnref{ControlOfParameters} again,
$$
\left\|\frac{d^{l-1}}{dx^{l-1}}\Cal{G}v_k\right\|_{C^0} \lesssim \left(1+T^k\right)\left(\epsilon R e^T + \frac{1}{R}\right)^{2k+3},
$$
and the result follows upon combining these relations.\qed
\proclaim{Lemma \nextprocno}
\noindent If $v_{k,c}$ is the $k$'th partial sum of $V$ with logarithmic parameter $c$, and if $v:[\epsilon R,\epsilon R^4]\rightarrow\Bbb{R}$ solves $\Cal{G}v=0$ with initial value
$$
v(0) = \frac{1}{R}\sigma[\epsilon,R,4k+9](c) + \frac{\epsilon R}{2},\eqnum{\nexteqnno[BetterInitialConditions]}
$$
then, for sufficiently large $\Lambda$,
$$
v = v_{k,c} + \opO\left(\left(1+x^{k+1}\right)\left(\epsilon R e^x + \frac{1}{R}e^{-x}\right)^{2k+3}\right).\eqnum{\nexteqnno[EstimateOfErrorInLogCoordinates]}
$$
\endproclaim
\proclabel{EstimateOfErrorInLogCoordinates}
\remark Since $r=\epsilon R e^x$, by the chain rule,
$$
\frac{d}{dr}=\frac{1}{r}\frac{d}{dx},
$$
so that Theorem \procref{ErrorForExactSolutionWithGivenLogarithmicParameter} follows immediately from \eqnref{EstimateOfErrorInLogCoordinates} upon setting $k=0$.
\medskip
\proof For non-negative, integer $l$, if $v:[0,3\opLog(R)]\rightarrow\Bbb{R}$ solves $\Cal{G}v=0$ with initial value as in \eqnref{InitialConditions} then, since \eqnref{DifferenceBetweenFormalAndExactToAllOrders} holds for all $T\in[0,3\opLog(R)]$,
$$
v = v_{l,c} + \opO\left(\left(1+x^{l+1}\right)\left(\epsilon R e^x + \frac{1}{R}\right)^{2l+3}\right).
$$
In particular, if $v:[0,3\opLog(R)]\rightarrow\Bbb{R}$ now solves $\Cal{G}v=0$ with initial value given by \eqnref{BetterInitialConditions}, then, bearing in mind \eqnref{ControlOfParameters},
$$
v = v_{4k+9,c} + \opO\left(\left(1+x^{k+1}\right)\left(\epsilon R e^x + \frac{1}{R}e^{-x}\right)^{2k+3}\right).
$$
However, by Lemma \procref{LemmaFormalSolutionOverSmallScale} and \eqnref{ControlOfParameters} again,
$$
v_{4k+9,c} = v_{k,c} + \opO\left(\left(1+x^{k+1}\right)\left(\epsilon R e^x + \frac{1}{R}e^{-x}\right)^{2k+3}\right),
$$
and the result follows.\qed
\newsubhead{The Small Scale - Solutions of the Linearised Equation}[TheSmallScaleJacobiFields]
We conclude this section by studying how solutions of the equation $\Cal{G}v=0$ vary with the logarithmic parameter $c$.
\proclaim{Theorem \nextprocno}
\noindent For sufficiently large $\Lambda$ and for all $R,\epsilon$ satisfying \eqnref{ControlOfParameters}, if, for all $c\in[-K,K]$, the function $v_c:[\epsilon R,\epsilon R^4]\rightarrow\Bbb{R}$ solves $\Cal{G}v_c=0$ with initial value
$$
v_c(\epsilon R) = \frac{1}{R}\sigma[\epsilon,R](c) + \frac{\epsilon R}{2},
$$
then
$$
\frac{dv_c}{dc}(r) = \frac{\epsilon}{r} + \opO\left(\left[1+\opLog\left(\frac{r}{\epsilon R}\right)\right]\frac{1}{r^k}\left(r + \frac{\epsilon}{r}\right)^3\right).\eqnum{\nexteqnno[ErrorForExactJacobiFieldWithGivenLogarithmicParameter]}
$$
\endproclaim
\proclabel{ErrorForExactJacobiFieldWithGivenLogarithmicParameter}
\noindent Theorem \procref{ErrorForExactJacobiFieldWithGivenLogarithmicParameter} follows from \eqnref{SizeOfErrorForJacobiFields}, below, via a similar reasoning to that used in Section \subheadref{TheSmallScaleExactSolutionsII}. It suffices to study solutions of the linearisation of $\Cal{G}$ about $v$, the asymptotic properties of which are readily derived from the analysis of the previous sections. Indeed, let $\Bbb{R}[X][[M,N]]$ be as in Section \subheadref{TheSmallScaleFormalSolutions} and define the operator $\partial_N$ over this space by
$$
(\partial_N V)_{p,q} := (q+1)V_{p,q+1}.\eqnum{\nexteqnno[FormalNDerivative]}
$$
In other words, $\partial_N$ is simply the operator of formal differentiation with respect to $N$. By explicit calculation, $N\partial_N$ commutes with $\Cal{D}$. Now let $V$ be the formal power series constructed in Lemma \procref{LemmaFormalSolutionOverSmallScale} and define
$$
W := N\partial_N V.\eqnum{\nexteqnno[FormalJacobiField]}
$$
Applying $N\partial_N$ to the relation $\Cal{G}V=0$ yields
$$
\Cal{D}W + 3V^2W - 2MVW = 0,\eqnum{\nexteqnno[FormalJacobiEquation]}
$$
so that $W$ is a formal solution to the linearisation of $\Cal{G}$ about the formal series $V$.
\par
Fix a non-negative integer $k$, let $\hat{V}_k$ be as in Lemma \procref{LemmaFormalSolutionOverSmallScale} and denote
$$
\hat{W}_k := \sum_{p+q\leq 2k+1}W_{p,q}(X)M^pN^q.\eqnum{\nexteqnno[PartialSumOfFormalJacobiField]}
$$
By \eqnref{FormalJacobiEquation},
$$
\Cal{D}\hat{W}_k + 3\hat{V}_k^2\hat{W}_k - 2M\hat{V}_k\hat{W}_k = \opO((M + N)^{2k+3}).\eqnum{\nexteqnno[FormalJacobiFieldSolvesFormalJacobiEquation]}
$$
Consider now $\Lambda,K>0$, let $\epsilon,R>0$ and $c\in\Bbb{R}$ satisfy \eqnref{ControlOfParameters}, and let $v_{k,c}$ and $w_{k,c}$ be the functions corresponding to $\hat{V}_k$ and $\hat{W}_k$ respectively. By \eqnref{FormalJacobiFieldSolvesFormalJacobiEquation}, for all $k$,
$$
\Cal{D}w_{k,c} + 3v_{k,c}^2w_{k,c} - 2(\epsilon Re^x)v_{k,c}w_{k,c} = \opO\left(x^{k+1}\left(\epsilon R e^x + \frac{1}{R}e^{-x}\right)^{2k+3}\right).\eqnum{\nexteqnno[FunctionIsApproximateJacobiField]}
$$
\proclaim{Lemma \nextprocno}
\noindent For sufficiently large $\Lambda$ and for all $T\in[0,3\opLog(R)]$, if $v:[0,T]\rightarrow\Bbb{R}$ solves $\Cal{G}v=0$ with initial value
$$
v(0) = \frac{1}{R}\sigma[\epsilon,R,k](c) + \frac{\epsilon R}{2},
$$
and if $w:[0,T]\rightarrow\Bbb{R}$ solves
$$
\Cal{D}w + 3v^2 w - 2\epsilon R e^x vw = 0\eqnum{\nexteqnno[ExactJacobiField]}
$$
with initial value $w(0)=w_{k,c}(0)$, then
$$
\|w-w_{k,c}\|_{C^1_0} \lesssim (1+T)^{k+1}\left(\epsilon R e^T + \frac{1}{R}\right)^{2k+3}.\eqnum{\nexteqnno[SizeOfErrorForJacobiFields]}
$$
\endproclaim
\proof Indeed, by \eqnref{FunctionIsApproximateJacobiField},
$$
\left\|\Cal{D}w_{k,c} + 3v_{k,c}^2w_{k,c} - 2(\epsilon R e^x)v_{k,c}w_{k,c}\right\| \lesssim (1+T)^{k+1}\left(\epsilon R e^T + \frac{1}{R}\right)^{2k+3}.
$$
Observe that
$$
\|v\|_{C^0},\|v_{k,c}\|_{C^0},\|w_{k,c}\|_{C^0}\lesssim 1.
$$
Thus, by \eqnref{GlobalEstimateOfError},
$$\eqalign{
\left\|\left(3v_{k,c}^2 - 3v^2\right)w_{k,c}\right\|_{C^0}
&=3\left\|(v_{k,c} - v)(v_{k,c} + v)w_{k,c}\right\|_{C^0}\cr
&\lesssim (1+T)^{k+1}\left(\epsilon R e^T + \frac{1}{R}\right)^{2k+3}.\cr}
$$
Likewise
$$
\left\|\left(2(\epsilon R e^x)v_{k,c}-2(\epsilon R e^x)v\right)w_{k,c}\right\|_{C^0}
\lesssim (1+T)^{k+1}\left(\epsilon R e^T + \frac{1}{R}\right)^{2k+3}.
$$
Thus
$$\multiline{
\left\|\Cal{D}(w_{k,c}-w) + 3v^2(w_{k,c}-w)-2(\epsilon R e^x)v(w_{k,c}-w)\right\|_{C^0}\cr
\qquad\qquad\qquad=\left\|\Cal{D}w_{k,c} + 3v^2w_{k,c} - 2(\epsilon R e^x)vw_{k,c}\right\|_{C^0}\cr
\qquad\qquad\qquad\lesssim (1+T)^{k+1}\left(\epsilon R e^T + \frac{1}{R}\right)^{2k+3}.\cr}
\eqnum{\nexteqnno[ImageOfDifference]}
$$
Observe now that, for all $\phi:[0,T]\rightarrow\Bbb{R}$,
$$
3v^2\phi - 2\epsilon R e^x \phi = \Cal{E}(v-v_{k,c})\phi,
$$
where $\Cal{E}$ is given by \eqnref{FormulaForErrorOperator}. In particular, by \eqnref{ControlOfParameters}, \eqnref{SizeOfErrorOperator} and \eqnref{GlobalEstimateOfError}, the operator norm of $\Cal{E}(v-v_{k,c})$ consider as a map from $C^1_0([0,T])$ into $C^0([0,T])$ satisfies
$$
\left\|\Cal{E}(v-v_{k,c})\right\| \lesssim T\left(\left(\epsilon R e^T\right)^2+\frac{1}{R^2}\right).
$$
Thus, by \eqnref{NormOfDAndItsInverse}, for sufficiently large $\Lambda$, the operator $\Cal{D} + \Cal{E}(v-v_{k,c})$ defines an invertible map from $C^1_0([0,T])$ into $C^0([0,T])$ and the result now follows by \eqnref{ImageOfDifference}.\qed
\medskip
\noindent Theorem \eqnref{ErrorForExactJacobiFieldWithGivenLogarithmicParameter} now follows as indicated above. In addition, a further iteration of this process also yields
\proclaim{Theorem \nextprocno}
\noindent With the same hypotheses as in Theorem \procref{ErrorForExactJacobiFieldWithGivenLogarithmicParameter},
$$
\frac{d^2v_c}{dc^2}(r) = \opO\left(\left[1+\opLog\left(\frac{r}{\epsilon R}\right)\right]\frac{1}{r^k}\left(r + \frac{\epsilon}{r}\right)^3\right).\eqnum{\nexteqnno[ErrorForSecondDerivativeOfGrimProfile]}
$$
\endproclaim
\newhead{The Grim Paraboloid}[TheGrimParaboloid]
\newsubhead{The MCFS Jacobi Operator}[TheMCFSJacobiOperator]
The {\bf Grim paraboloid}, which we henceforth denote by $G_0$, is defined to be the unique rotationally symmetric MCF soliton which is a graph over the whole of $\Bbb{R}^2$. Put differently, using the notation of Chapter \headref{GrimSurfaces}, there is a unique solution $v$ to the ODE $\Cal{G}v=0$ which is defined over the whole interval $]0,\infty[$. This solution tends to $0$ as $x$ tends to $0$, and the Grim paraboloid is the surface of revolution generated by rotating the graph of its integral about the $z$-axis.
\par
Let $J$ be the MCFS Jacobi operator of the Grim paraboloid as defined in Appendix \subheadref{SurfaceGeometry}. In this section, we show that this operator defines a linear isomorphism over suitably weighted Sobolev and H\"older spaces. We first describe the spaces of interest to us (see Appendix \subheadref{FunctionSpaces} for details). Let $g$ denote the metric induced over $\Bbb{R}^2$ by the graph $G_0$, that is
$$
g := (1+v^2)dr^2 + r^2d\theta^2.\eqnum{\nexteqnno[MetricOverPlane]}
$$
For all non-negative, integer $m$, let $\|\cdot\|_{H^m(G)}$ denote the {\bf Sobolev norm} of order $m$ of functions over $\Bbb{R}^2$ with respect to this metric. Likewise, for all non-negative, integer $m$, and for all $\alpha\in[0,1]$, let $\|\cdot\|_{C^{m,\alpha}(G)}$ denote the {\bf H\"older norm} of order $(m,\alpha)$ of functions over $\Bbb{R}^2$ with respect to this metric. Observe that, by \eqnref{AsymptoticSolutionAtInfinity}, these Sobolev and H\"older norms are uniformly equivalent to the Sobolev and H\"older norms defined with respect to the more straightforward metric
$$
g' := (1+r^2)dr^2 + r^2d\theta^2.\eqnum{\nexteqnno[ApproximateMetricOverPlane]}
$$
For all non-negative, integer $m$, let $H^m(G)$ denote the {\bf Sobolev space} of measurable functions $f$ over $\Bbb{R}^2$ whose distributional derivatives up to and including order $m$ are locally square integrable and which satisfy $\|f\|_{H^m(G)}<\infty$. Likewise, for all non-negative, integer $m$, and for all $\alpha\in[0,1]$, let $C^{m,\alpha}(G)$ denote the {\bf H\"older space} of $m$-times differentiable functions $f$ over $\Bbb{R}^2$ which satisfy $\|f\|_{C^{m,\alpha}(G)}<\infty$. Recall that both $H^m(G)$ and $C^{m,\alpha}(G)$, furnished with the above norms, are Banach spaces.
\par
For all real $\gamma$, define $\phi_\gamma:\Bbb{R}^2\rightarrow\Bbb{R}$ by
$$
\phi_\gamma := e^{\frac{(1+\gamma)u}{2}}.\eqnum{\nexteqnno[LargeScaleWeight]}
$$
where $u$ here denotes the integral of $v$ with initial value $0$. For all non-negative, integer $m$, for all $\alpha\in[0,1]$ and for all real $\gamma$, define the {\bf weighted Sobolev and H\"older norms} of weight $\gamma$ over $\Bbb{R}^2$ by
$$\eqalign{
\|f\|_{H^m_\gamma(G)} &:= \|\phi_\gamma f\|_{H^m(G)},\ \text{and}\cr
\|f\|_{C^{m,\alpha}_\gamma(G)} &:= \|\phi_\gamma f\|_{C^{m,\alpha}(G)}.\cr}
\eqnum{\nexteqnno[DefinitionOfWeightedNormsOverG]}
$$
Observe that, by \eqnref{AsymptoticSolutionAtInfinity} again, these {\bf weighted Sobolev and H\"older norms} are uniformly equivalent to the weighted norms defined using instead of $\phi_\gamma$ the more straightforward weight function
$$
\phi'_\gamma := e^{\frac{(1+\gamma)r^2}{4}}.\eqnum{\nexteqnno[ApproximateLargeScaleWeight]}
$$
For all non-negative, integer $m$, for all $\alpha\in[0,1]$, and for all real $\gamma$, define the {\bf weighted Sobolev} and {\bf H\"older spaces} of {\bf weight} $\gamma$ over $\Bbb{R}^2$ by
$$\eqalign{
H^m_\gamma(G) &:= \left\{ f \ |\ \phi_\gamma f\in H^m(G)\right\},\ \text{and}\cr
C^{m,\alpha}_\gamma(G) &:= \left\{ f\ |\ \phi_\gamma f\in C^{m,\alpha}(G)\right\}.\cr}
\eqnum{\nexteqnno[DefinitionOfWeightedSpacesOverG]}
$$
These spaces, furnished with the weighted Sobolev and H\"older norms are trivially also Banach spaces.
\par
Since $G_0$ is a graph over $\Bbb{R}^2$, its MCFS Jacobi operator may be thought of as an operator acting on functions over $\Bbb{R}^2$. In particular, as we will see presently, for all $\alpha\in[0,1]$, and for all real $\gamma$, $J$ defines bounded linear maps from $H^2_\gamma(G)$ into $H^0_\gamma(G)$ and from $C^{2,\alpha}_\gamma(G)$ into $C^{0,\alpha}_\gamma(G)$. We show
\proclaim{Theorem \nextprocno}
\myitem{(1)} For all sufficiently small $\gamma$, $J$ defines a linear isomorphism from $H^2_\gamma(G)$ into $H^0_\gamma(G)$.
\medskip
\myitem{(2)} For all $\alpha\in]0,1[$ and for all sufficiently small $\gamma$, $J$ defines a linear isomorphism from $C^{2,\alpha}_\gamma(G)$ into $C^{0,\alpha}_\gamma(G)$.
\endproclaim
\proclabel{JIsLinearIsomorphismOverG}
\noindent Theorem \procref{JIsLinearIsomorphismOverG} will follow from Lemmas \procref{LemmaInvertibilityOverSobolev} and \procref{ThmInvertibilityOverHolder} below. Before proceeding, we first observe that, for all $\gamma$, the function $\phi_\gamma$ is strictly positive so that, for all non-negative, integer $m$, and for all $\alpha\in[0,1]$, the operator of multiplication by this function, which we denote by $M_\gamma$, defines linear isomorphisms from $H^m_\gamma(G)$ into $H^m(G)$ and from $C^{m,\alpha}_\gamma(G)$ into $C^{m,\alpha}(G)$. For all real $\gamma$, we therefore define,
$$
J_\gamma := M_\gamma J M_\gamma^{-1}.\eqnum{\nexteqnno[PhiWeightedJacobiOperator]}
$$
This operator is none other than the $\phi_\gamma$-Jacobi operator of the Grim paraboloid, which has been studied in detail in \cite{ChengZhou}, \cite{ChengMejiaZhouI}, \cite{ChengMejiaZhouII} and \cite{ChengMejiaZhouIII}. Trivially, $J$ defines linear isomorphisms from $H^2_\gamma(G)$ into $H^0_\gamma(G)$ and from $C^{2,\alpha}_\gamma(G)$ into $C^{0,\alpha}_\gamma(G)$ if and only if $J_\gamma$ defines linear isomorphisms from $H^2(G)$ into $H^0(G)$ and from $C^{2,\alpha}(G)$ into $C^{0,\alpha}(G)$ respectively.
\proclaim{Lemma \nextprocno}
\noindent For all real $\gamma$,
$$
J_\gamma f = \Delta^{G_0} f - \gamma\langle e_z,\nabla^{G_0} f\rangle + \frac{(\gamma^2-1)}{4}f - \frac{(1+\gamma)^2}{4}\langle e_z,N_{G_0}\rangle^2 f + \opTr(A_{G_0}^2) f.\eqnum{\nexteqnno[FormulaForConjugateOperator]}
$$
\endproclaim
\proof By \eqnref{AppGradientAndHessian},
$$\eqalignno{
\nabla^{G_0}\phi_\gamma^{-1} &= -\frac{(1+\gamma)}{2\phi_\gamma}\pi^{G_0}(e_z),\ \text{and}\cr
\opHess^{G_0}\phi_\gamma^{-1} &= \frac{(1+\gamma)^2}{4\phi_\gamma}dz\otimes dz + \frac{(1+\gamma)}{2\phi_\gamma}\langle e_z,N_{G_0}\rangle\opII^{G_0}.\cr}
$$
However, since ${G_0}$ is a mean curvature flow soliton, $H_{G_0}=-\langle e_z,N_{G_0}\rangle$, and taking the trace therefore yields
$$
\Delta^{G_0}\phi_\gamma^{-1} = \frac{(1+\gamma)^2}{4\phi_\gamma} - \frac{(1+\gamma)(3+\gamma)}{4\phi_\gamma}\langle e_z,N_{G_0}\rangle^2.
$$
Thus, by \eqnref{AppMCFSJacobiOperator},
$$
\phi_\gamma J_0\phi_\gamma^{-1} = \frac{(\gamma^2-1)}{4} - \frac{(1+\gamma)^2}{4}\langle e_z,N_{G_0}\rangle^2 + \opTr(A_{G_0})^2.
$$
The result now follows by \eqnref{AppCommutator}.\qed
\medskip
By \eqnref{AppBasicFormulaInCoordinateCharts} and \eqnref{SpecialCaseARLargeScale},
$$\eqalign{
\langle e_z,N_{G_0}\rangle^2 &= \opO(r^{-(2+k)}),\ \text{and}\cr
\opTr(A_{G_0}^2) &= \opO(r^{-(2+k)}).\cr}\eqnum{\nexteqnno[CoefficientsOfConjugateOperator]}
$$
It follows that, as $\gamma$ tends to $0$, the family $(J_\gamma)$ converges to $J_0$ in every operator norm of relevance to us. Since invertibility is stable under small perturbations, it is therefore sufficient to consider only the case $\gamma=0$ where, in particular, $J_0$ is self-adjoint.
\par
We now derive a formula for $J_0$ which is better adapted to our purposes. First, let $c:]0,\infty[\rightarrow\Bbb{R}$ be such that, for all $r$, $c(r)$ is the geodesic curvature of the circle $C(r)$ with respect to the metric induced by the graph $G_0$ over $\Bbb{R}^2$.
\proclaim{Lemma \nextprocno}
\noindent The function $c$ is given by
$$
c=\frac{1}{r}\langle e_z,N_{G_0}\rangle.\eqnum{\nexteqnno[IntrinsicCurvatureOfCircles]}
$$
In particular, for large values of $r$,
$$
c = \opO(r^{-(2+k)}).\eqnum{\nexteqnno[EqnAsymptoticBoundsOnLittleKappa]}
$$
\endproclaim
\proof Let $D$ denote the Levi-Civita covariant derivative of the Euclidean metric over $\Bbb{R}^3$. Think of $C(r)$ as a horizontal circle in $\Bbb{R}^3$ at height $u(r)$, where $u$ here denotes the integral of $v$ with initial value $0$. In particular, $D_{e_\theta}e_\theta = \frac{1}{r}e_r$, where $e_\theta$ and $e_r$ denote respectively the unit, horizontal vector fields in the angular and radial directions about the $z$-axis. Since the geodesic curvature of $C(r)$ with respect to the induced metric over $G_0$ is equal to the length of the tangential component of this vector, the function $c$ is given by
$$
c = \frac{1}{r}\sqrt{1 - \langle e_r, N_{G_0}\rangle^2} = \frac{1}{r}\langle e_z,N_{G_0}\rangle,
$$
as desired. \eqnref{EqnAsymptoticBoundsOnLittleKappa} now follows from \eqnref{CoefficientsOfConjugateOperator}, and this completes the proof.\qed
\medskip
\noindent Let $\rho:]0,\infty[\rightarrow\Bbb{R}$ be such that, for all $r$, $\rho(r)$ is the intrinsic distance along $G_0$ of any point on the circle $C(r)$ from the origin. Since $\rho$ is obtained by integrating $\sqrt{1+v^2}$, by \eqnref{AsymptoticSolutionAtInfinity} again, for large values of $r$,
$$\eqalign{
\rho_r &= r + \opO(r^{-(k+1)}),\ \text{and}\cr
r_\rho &= \frac{1}{r} + \opO(r^{-(k+3)}),\cr}\eqnum{\nexteqnno[AsymptoticFormulaForRho]}
$$
where the subscripts $r$ and $\rho$ here denote differentiation with respect to the variables $r$ and $\rho$ respectively.
\proclaim{Lemma \nextprocno}
\noindent Away from the $z$-axis,
$$
J_0 f = f_{\rho\rho} + f_{\theta\theta} + c f_\rho - \frac{1}{4}f + \psi f,\eqnum{\nexteqnno[ConjugateJacobiOperatorInPolarCoordinates]}
$$
where the subscripts $\rho$ and $\theta$ denote differentiation along the unit radial and unit angular directions in $G_0$ and, for large values of $\rho$,
$$
\left|\psi\right| \lesssim \rho^{-1}.\eqnum{\nexteqnno[ZeroethOrderTermInConjugateJacobiOperatorInPolarCoordinates]}
$$
\endproclaim
\proof Indeed, away from the $z$-axis,
$$
\Delta^{G_0} f = f_{\rho\rho} + f_{\theta\theta} + c f_\rho,
$$
so that \eqnref{ConjugateJacobiOperatorInPolarCoordinates} follows by \eqnref{FormulaForConjugateOperator} and \eqnref{CoefficientsOfConjugateOperator} with
$$
\left|\psi\right| \lesssim r^{-2}.
$$
Finally, integrating \eqnref{AsymptoticFormulaForRho}, yields $\rho\lesssim r^2$, so that $r^{-2}\lesssim\rho^{-1}$ and the result follows.\qed
\newsubhead{Invertibility over Sobolev Spaces}[InvertibilityOverSobolevSpaces]
We now obtain the invertibility of $J_0$ for Sobolev spaces. The main technical difficulty here arises from the non-compactness of the ambient space. This is compensated for by the following estimate.
\proclaim{Lemma \nextprocno}
\noindent There exist $B,R>0$ such that for all $f$ in $H^2(G)$,
$$
\|f|_{A(R,\infty)}\|_{L^2(G)} \leq B\left(\|f|_{A(R-1,R+1)}\|_{L^2(G)} + \|J_0f|_{A(R-1,\infty)}\|_{L^2(G)}\right).\eqnum{\nexteqnno[EqnSobolevMassConcentration]}
$$
\endproclaim
\proclabel{SobolevControlAtInfinity}
\proof Since $C_0^\infty(G)$ is dense in $H^2(G)$, it suffices to prove the result when $f$ is smooth and has compact support. Denote $g:=J_0 f$ and define $\alpha,\beta:]0,\infty[\rightarrow\Bbb{R}$ by
$$\eqalign{
\alpha(\rho) &:= \int_{C(\rho)}f^2 dl,\ \text{and}\cr
\beta(\rho) &:= \int_{C(\rho)}g^2 dl,\cr}
$$
where $C(\rho)$ here denotes the circle of points lying at intrinsic distance $\rho$ along $G_0$ from the origin. Twice differentiating $\alpha$ yields
$$\eqalign{
\alpha_{\rho} &= \int_{C(\rho)}2f f_\rho + f^2c dl,\ \text{and}\cr
\alpha_{\rho\rho} &= \int_{C(\rho)}2f_\rho^2 + 2f f_{\rho\rho} + 4f f_\rho c + f^2 c_\rho + f^2 c^2 dl,\cr}
$$
where the subscript $\rho$ here denotes differentiation with respect to this variable. By \eqnref{ConjugateJacobiOperatorInPolarCoordinates},
$$
\alpha_{\rho\rho} = \int_{C(\rho)}2f_\rho^2 - 2ff_{\theta\theta} + \frac{1}{2}f^2 - 2\psi f^2 + 2fg + 2f f_\rho c + f^2 c_\rho + f^2 c^2 dl.
$$
Integrating the term $2ff_{\theta\theta}$ by parts and applying the algebraic-geometric mean inequality now yields
$$
\alpha_{\rho\rho} \geq \int_{C(\rho)}\left(\frac{1}{4} - 2\psi + c_\rho - c^2\right)f^2 - 4g^2 dl.
$$
However, by \eqnref{EqnAsymptoticBoundsOnLittleKappa}, \eqnref{AsymptoticFormulaForRho} and \eqnref{ZeroethOrderTermInConjugateJacobiOperatorInPolarCoordinates}, $c$, $c_\rho=c_rr_\rho$ and $\psi$ all tend to $0$ as $\rho$ tends to $+\infty$ so that, for sufficiently large $\rho$
$$
\alpha_{\rho\rho} \geq \frac{1}{8}\alpha - 4\beta.
$$
Since $f$ has compact support, upon integrating this relation we obtain, for sufficiently large $R$,
$$
\|f|_{A(R,\infty)}\|^2_{L^2(G)} = \int_R^\infty\alpha d\rho
\leq 32\int_R^\infty\beta d\rho - 8\alpha_{\rho}(R)
=32\|\hat{J}_0 f|_{A(R,\infty)}\|_{L^2(G)}^2 - 8\alpha_{\rho}(R).
$$
However, by the Sobolev trace formula and classical elliptic estimates,
$$\eqalign{
\alpha_{\rho}(R) &\leq B_1\|f|_{A(R-1/2,R+1/2)}\|^2_{H^2(G)}\cr
&\leq B_2\left(\|f|_{A(R-1,R+1)}\|^2_{L^2(G)} + \|J_0f|_{A(R-1,R+1)}\|_{L^2(G)}^2\right),}
$$
for suitable constants $B_1$ and $B_2$. The result now follows upon combining the last two relations.\qed
\medskip
\noindent Combining Lemma \procref{SobolevControlAtInfinity} with classical elliptic estimates (see \cite{GilbTrud}) yields
\proclaim{Lemma \nextprocno}
\noindent There exist $B,R>0$ such that for all $f$ in $H^2(G)$,
$$
\|f\|_{H^2(G)} \leq B\left(\|f|_{B(R)}\|_{L^2(G)} + \|J_0 f\|_{L^2(G)}\right).\eqnum{\nexteqnno[EqnSobolevEllipticEstimate]}
$$
\endproclaim
\noindent Since $J_0$ is self-adjoint, this in turn yields
\proclaim{Lemma \nextprocno}
\noindent $J_0$ defines a Fredholm map from $H^2(G)$ into $L^2(G)$ of Fredholm index equal to $0$.
\endproclaim
\proclabel{LemmaSobolevFredholm}
\noindent It remains only to prove that $J_0$ has trivial kernel in $H^2(G)$. We obtain a slightly more general result which will serve also for the H\"older space case of the following section.
\proclaim{Lemma \nextprocno}
\noindent There exists no non-trivial, bounded function $f:G_0\rightarrow\Bbb{R}$ such that $J_0f=0$.
\endproclaim
\proclabel{LemmaNoBoundedJacobiFields}
\proof Indeed, suppose that there exists a non-trivial bounded function $f:G_0\rightarrow\Bbb{R}$ such that $J_0f=0$. Upon multiplying by $(-1)$, we may suppose that $f$ is positive at some point. Now, since all vertical translates of $G_0$ are also mean curvature flow solitons, the function $\mu=\langle e_z,N_G\rangle$ is a Jacobi field over this surface, that is,
$$
J_0 \phi_0\mu = \phi_0 J \mu = 0.
$$
Since $G_0$ is a graph, the function $\mu$ is everywhere strictly positive. It follows that $\phi_0\mu$ is also positive, so that the quotient $f/\phi_0\mu$ is smooth. Since $\phi_0\gtrsim e^{r^2/4}$ and $\mu=\opO(r^{-1})$, the function $\phi_0\mu$ tends to infinity as $r$ tends to infinity, and so $f/\phi_0\mu$ attains its maximum value at some point $x$, say, of $G_0$. In particular, upon rescaling, we may suppose that $f/\phi_0\mu\leq 1$ and that $f(x)/\phi_0(x)\mu(x)=1$.
\par
Bearing in mind that $\mu$ is positive, we define the operator $J_\mu:=M_\mu^{-1} J M_\mu$, where $M_\mu$ here denotes the operator of multiplication by $\mu$. Since $J\mu=0$, by \eqnref{AppCommutator}, this operator has no zero'th order term. Thus, since $J_\mu(f/\mu\phi_0)=(1/\mu\phi_0)J_0 f=0$, it follows by the strong maximum principle that $f/\phi_0\mu$ is constant and equal to $1$. However, since $\phi_0\mu$ is unbounded, this is absurd, and the result follows.\qed
\proclaim{Corollary \nextprocno}
\noindent $J_0$ has trivial kernel in $H^2(G)$.
\endproclaim
\proclabel{LemmaSobolevKernelIsEmpty}
\proof Indeed, by the Sobolev embedding theorem, every element of $H^2(G)$ is bounded, and the result now follows by Lemma \procref{LemmaNoBoundedJacobiFields}.\qed
\medskip
\noindent The above results together with a perturbation argument now yield
\proclaim{Lemma \nextprocno}
\noindent For sufficiently small $\gamma$, $J$ defines a linear isomorphism from $H^2_{\gamma}(G)$ into $H^0_\gamma(G)$.
\endproclaim
\proclabel{LemmaInvertibilityOverSobolev}
\newsubhead{Invertibility over H\"older Spaces}[InvertibilityOverHoelderSpaces]
We prove the invertibility of $J_0$ over $C^{2,\alpha}(G)$ in essentially the same manner. We first require the following preliminary result.
\proclaim{Lemma \nextprocno}
\noindent Let $\alpha$ and $\beta$ be positive constants. If $\phi:[0,\infty[\rightarrow]0,\infty[$ is a bounded, positive function such that $\phi''\geq\alpha^2\phi - \beta$ in the viscosity sense, then, for all $t$,
$$
\phi(t) \leq \opMax(\phi(0) - \beta/\alpha^2,0)e^{-\alpha t} + \beta/\alpha^2.\eqnum{\nexteqnno[ViscositySolutionsToODEInequality]}
$$
\endproclaim
\proclabel{LemmaViscositySolutionsToODEInequality}
\proof Let $A=\opMax(\phi(0) - \beta/\alpha^2,0)$ and let $B=\msup_{t\in[0,\infty[}\phi(t)$. Fix $T>0$ and define
$$
f = \frac{Be^{\alpha T} - A}{e^{2\alpha T} - 1}e^{\alpha t} + \frac{A - B e^{-\alpha T}}{1 - e^{-2\alpha T}}e^{-\alpha t} + \beta/\alpha^2.
$$
In other words, $f$ is the unique solution of the ODE problem $f_{tt}=\alpha^2 f - \beta$ with boundary values $f(0)=A+\beta/\alpha^2\geq \phi(0)$ and $f(T)=B+\beta/\alpha^2\geq \phi(T)$. Let $C$ be the minimum value of $f-\phi$ over $[0,T]$ and let $t\in[0,T]$ be the point at which this minimum is attained. If $t$ is a boundary point of this interval, then $C\geq 0$. Otherwise, $f-C\geq \phi$ and $f(t)-C=\phi(t)$. Thus, since $\phi$ is a viscosity solution of $\phi''\geq \alpha^2\phi - \beta$, at this point, we have
$$
\alpha^2 f - \beta = (f-C)_{tt}\geq\alpha^2 (f-C) - \beta
$$
so that, once again, $C\geq 0$. In each case, we therefore obtain
$$
\phi \leq f = \frac{Be^{\alpha T} - A}{e^{2\alpha T} - 1}e^{\alpha t} + \frac{A - B e^{-\alpha T}}{1 - e^{-2\alpha T}}e^{-\alpha t} + \beta/\alpha^2,
$$
and the result follows upon taking the limit as $T$ tends to $+\infty$.\qed
\medskip
\noindent As in the Sobolev case, the non-compactness of the ambient space is compensated for by the following estimate.
\proclaim{Lemma \nextprocno}
\noindent There exist $B,R>0$ such that for all $f$ in $C^{2,\alpha}(G)$,
$$
\|f|_{A(R,\infty)}\|_{C^0(G)} \leq B\left(\|f|_{C(R)}\|_{C^0(G)} + \|J_0f|_{A(R-1,\infty)}\|_{C^0(G)}\right).\eqnum{\nexteqnno[EqnHolderMassConcentration]}
$$
\endproclaim
\proof Define $\alpha:]0,\infty[\rightarrow\Bbb{R}$ by
$$
\alpha(\rho) := \msup_{x\in C(\rho)} f(x)^2,
$$
where $C(\rho)$ here denotes the circle of points lying at intrinsic distance $\rho$ along $G_0$ from the origin. Denote $g:=J_0f$, and define $B\geqslant 0$ by
$$
B := \|g^2|_{A(R,\infty)}\|_{C^0(G)}.
$$
Choose $x\in C(\rho)$ maximising $f^2$, and observe that $ff_{\theta\theta}$ is non-positive at this point. Thus, bearing in mind \eqnref{ConjugateJacobiOperatorInPolarCoordinates},
$$\eqalign{
(f^2)_{\rho\rho} &= 2f_\rho^2 + 2ff_{\rho\rho},\cr
&\geq 2f_\rho^2 + 2fg - 2cf_\rho + \frac{1}{2}f^2 - 2\psi f^2,\cr
&\geq \left(\frac{1}{4} - \frac{1}{2}c^2 - 2\psi\right)f^2 - 4g^2.\cr}
$$
By \eqnref{EqnAsymptoticBoundsOnLittleKappa} and \eqnref{ZeroethOrderTermInConjugateJacobiOperatorInPolarCoordinates}, for sufficiently large $\rho$
$$
(f^2)_{\rho\rho} \geq \frac{1}{8}f^2 - 4g^2.
$$
Since $\alpha$ is the envelope of the restriction of $f(x)^2$ to each radial line, it follows that over $[R,\infty[$,
$$
\alpha_{\rho\rho} \geq \frac{1}{8}\alpha - 4B,
$$
in the viscosity sense. Thus, by Lemma \procref{LemmaViscositySolutionsToODEInequality},
$$
\msup_{x\in A(R,\infty)} f^2(x) = \msup_{\rho\geq R}\alpha(\rho) \leq \opMax(\|f^2|_{C(R)}\|_{C^0} - 32B,0) + 32B,
$$
and the result follows.\qed
\medskip
\noindent Using classical elliptic estimates again, this yields
\proclaim{Lemma \nextprocno}
\noindent There exist $B,R>0$ such that for all $f$ in $C^{2,\alpha}(G)$,
$$
\|f\|_{C^{2,\alpha}(G)} \leq B\left(\|f|_{B(R)}\|_{C^0(G)} + \|J_0f\|_{C^{0,\alpha}(G)}\right).\eqnum{\nexteqnno[EqnHolderEllipticEstimate]}
$$
\endproclaim
\medskip
\noindent Standard arguments of the theory of elliptic operators together with Lemmas \procref{LemmaNoBoundedJacobiFields} and \procref{LemmaInvertibilityOverSobolev} now yield, as before
\proclaim{Lemma \nextprocno}
\noindent For all $\alpha$ and for all sufficiently small $\gamma$, $J$ defines a linear isomorphism from $C^{2,\alpha}_\gamma(G)$ into $C^{0,\alpha}_\gamma(G)$.
\endproclaim
\proclabel{ThmInvertibilityOverHolder}
\newhead{Rotationally Symmetric Grim Ends}[GrimEnds]
\newsubhead{The Modified MCFS Jacobi Operator}[TheModifiedMCFSJacobiOperator]
We now consider the case of rotationally symmetric Grim ends. Let $\Lambda$ be a large, positive real number, let $K>0$ be fixed, and let $\epsilon, R>0$ and $c\in\Bbb{R}$ satisfy \eqnref{ControlOfParameters}. Let $v:[\epsilon R,\infty[\rightarrow\Bbb{R}$ solve \eqnref{GrimSurfaceEquation} with logarithmic parameter $c$ so that, by \eqnref{ErrorForExactSolutionWithGivenLogarithmicParameter}, over the interval $[\epsilon R,\epsilon R^4]$,
$$
v=\frac{1}{2}r + \frac{c\epsilon}{r} + O\left(\left[1+\opLog\left(\frac{r}{\epsilon R}\right)\right]\frac{1}{r^k}\left(r + \frac{\epsilon}{r}\right)^3\right).\eqnum{\nexteqnno[AsymptoticFormulaForGrimEndSmallScale]}
$$
Let $u:[\epsilon R,\infty[\rightarrow\Bbb{R}$ be a primitive of $v$, let $G$ be the Grim end generated by rotating the graph of $u$ about the $z$-axis, and let $J$ be its MCFS Jacobi operator, as defined in Appendix \subheadref{SurfaceGeometry}.
\par
Since $G$ is a graph over $A(\epsilon R,\infty)$, $J$ may again be thought of as an operator acting on functions over this annulus. For all non-negative, integer $m$, for all $\alpha\in[0,1]$, and for all real $\gamma$, we define the norms $\|\cdot\|_{H^m_\gamma(G)}$ and $\|\cdot\|_{C^{m,\alpha}_\gamma(G)}$ as in Section \headref{TheGrimParaboloid}. For all non-negative integer $m$, for all $\alpha\in[0,1]$ and for all real $\gamma$, we define the {\bf hybrid norm} with weight $\gamma$ of functions over $\Bbb{R}^m$ by
$$
\|f\|_{m,\alpha,\gamma} := \|f\|_{C^{m,\alpha}_\gamma(G)} + \frac{1}{(\epsilon R)}\|f\|_{H^m_\gamma(G)}.\eqnum{\nexteqnno[DefinitionHybridNorm]}
$$
We will see presently that the hybrid norm encapsulates the asymptotic behaviour of $J$ as $\Lambda$ tends to infinity. Let $\Cal{L}^{m,\alpha}_\gamma(G)$ denote the Banach space of $m$-times differentiable functions $f$ over $\Bbb{R}^2$ with finite hybrid norm. In this section, we show that, for sufficiently small $\gamma$, and for sufficiently large $\Lambda$, the operator $J$ {\sl more or less} defines linear isomorphisms from $\Cal{L}^{2,\alpha}_\gamma(G)$ into $\Cal{L}^{0,\alpha}_\gamma(G)$ and, furthermore, that the norms of this isomorphism and its inverse are uniformly bounded as $\Lambda$ tends to infinity. However, in order to properly formalise these assertions, it is first necessary to apply the following two modifications.
\par
First, since the zero'th order coefficient of $J$ diverges rapidly over the annulus $A(\epsilon R,\epsilon R^4)$ as $\Lambda$ tends to infinity, we study instead what we call the modified MCFS Jacobi operator. We recall that, since different modifications are applied at different scales, the definition of this operator varies according to context, and that the description of the general framework is deferred to Section \subheadref{ModifiedJacobiOperators}. In the present case, it is defined as follows. Let $\chi_1$ be the cut-off function of the transition region $A(1,2)$ as defined in Appendix \subheadref{GeneralDefinitions} and define $\psi:A(\epsilon R,\infty)\rightarrow\Bbb{R}$ by
$$
\psi(r) = \chi_1\langle e_z,N_G\rangle + (1-\chi_1),\eqnum{\nexteqnno[DefinitionOfWeightSmallScale]}
$$
where $N_G$ here denotes the upward-pointing unit normal vector field over $G$. Bearing in mind that $\psi$ is always positive, the {\bf modified MCFS Jacobi operator} of $G$ is now defined by
$$
\hat{J} := M_\psi^{-1}JM_\psi,\eqnum{\nexteqnno[SmallScaleModifiedJacobiOperator]}
$$
where $M_\psi$ here denotes the operator of multiplication by $\psi$.
\par
Next, observe that $\hat{J}$ is in fact only defined over the annulus $A(\epsilon R,\infty)$. We thus extend it to an operator defined over the whole of $\Bbb{R}^2$ as follows. Given a function $\phi:A(\epsilon R,\infty)\rightarrow\Bbb{R}$, we define its {\bf canonical extension} $\tilde{\phi}:\Bbb{R}^2\rightarrow\Bbb{R}$ such that $\tilde{\phi}(x)=\phi(x)$ over $A(\epsilon R,\infty)$, $\tilde{\phi}(0)$ is equal to the mean value of $\phi$ over the circle $C(\epsilon R)$, and $\tilde{\phi}$ restricts to a linear function over every radial line in $B(\epsilon R)$. In particular, if $\phi$ is Lipschitz, then so too is $\tilde{\phi}$, and
$$
\|\tilde{\phi}\|_{C^{0,1}} \leq \pi\|\phi\|_{C^{0,1}}.
$$
Now, given a linear operator $L$ over $A(\epsilon R,\infty)$, we define its {\bf canonical extension} $\tilde{L}$ to be the operator over $\Bbb{R}^2$ whose coefficients are the canonical extensions of each of the coefficients of $L$. We henceforth identify all operators with their canonical extensions over $\Bbb{R}^2$. Observe, in particular, that if $L$ has any rotational symmetries, then so too does its canonical extension.
\proclaim{Theorem \nextprocno}
\noindent For all sufficiently small $\alpha\in]0,1[$ and for all sufficiently large $\Lambda$, $\hat{J}$ defines a linear isomorphism from $\Cal{L}^{2,\alpha}_\gamma(G)$ into $\Cal{L}^{0,\alpha}_\gamma(G)$. Furthermore, the operator norms of $\hat{J}$ and its inverse are uniformly bounded independent of $\Lambda$.
\endproclaim
\proclabel{JIsLinearIsomorphismOverEnds}
\noindent Theorem \procref{JIsLinearIsomorphismOverEnds} follows from Theorem \procref{JIsLinearIsomorphismOverG} by a perturbation argument and Lemmas \procref{TheOperatorDConvergesToZero} and \procref{TheOperatorEConvergesToZero}, below.
\par
We conclude this section by deriving formulae for $\hat{J}$ over different regions.
\proclaim{Lemma \nextprocno}
\noindent Over $A(\epsilon R,1)$, the modified MCFS Jacobi operator of $G$ is given by
$$
\hat{J}f = g^{ij}f_{ij} - 2\mu g^{ip}g^{jq}u_{pq}u_jf_i.\eqnum{\nexteqnno[ExplicitFormulaForModifiedJacobiOperator]}
$$
\endproclaim
\proof First observe that, for every tangent vector $X$ over $G$,
$$
\langle\nabla^G\psi,X\rangle = X\psi = X\langle N_G,e_z\rangle=\langle D_XN_G,e_z\rangle = \langle A_GX,e_z\rangle = \langle X,A_G\pi^G(e_z)\rangle,
$$
and so,
$$
\nabla^G\psi = A_G\pi^G(e_z).
$$
Since every vertical translate of $G$ is also a rotationally symmetric Grim end, $J\langle e_z, N_G\rangle=0$, and so, by \eqnref{AppCommutator},
$$
\hat{J}f = \Delta^G f + \langle e_z, \nabla^G f\rangle + 2\psi^{-1}\langle A_G\nabla^G f,e_z\rangle.
$$
By \eqnref{AppGradientAndHessian},
$$
\opHess^Gf=\opHess(f)\circ\pi - \langle D(f\circ\pi),N\rangle\opII_G.
$$
Furthermore, since $D(f\circ\pi)$ is horizontal
$$
\langle D(f\circ\pi),N_G\rangle
=-\frac{1}{\langle N_G,e_z\rangle}\langle D(f\circ\pi),e_z - \langle N_G,e_z\rangle N_G\rangle
=-\frac{1}{\langle N_G,e_z\rangle}\langle\nabla^G f,e_z\rangle.
$$
Taking the trace therefore yields
$$
\Delta^G f = g^{ij}f_{ij} + \frac{1}{\langle N_G,e_z\rangle}\langle\nabla^G f,e_z\rangle H_G.
$$
However, since $G$ is a mean curvature flow soliton, $H_G=-\langle N,e_z\rangle$, and so
$$
\Delta^G f = g^{ij}f_{ij} - \langle \nabla^G f,e_z\rangle.
$$
We conclude that
$$
\hat{J}f = g^{ij} f_{ij} + 2\psi^{-1}\langle A_G\nabla^G f,e_z\rangle,
$$
and the result now follows by \eqnref{AppBasicFormulaInCoordinateCharts}.\qed
\proclaim{Lemma \nextprocno}
\noindent Over $A(\epsilon R,2\epsilon R^4)$, the modified MCFS Jacobi operator of $G$ satisfies
$$
\hat{J}f = \Delta f - \left(\frac{1}{2} + \frac{c\epsilon}{r^2}\right)^2x^ix^jf_{ij} - \left(\frac{1}{2} - \frac{2c^2\epsilon^2}{r^4}\right)x^if_i + \Cal{E}_Gf,\eqnum{\nexteqnno[ModifiedJacobiOperatorOfGrimEnd]}
$$
where $\Cal{E}_Gf:=a^{ij}f_{ij} + b^if_i$, and $a$ and $b$ satisfy
$$\eqalign{
a &=\opO\left(\left[1+\opLog\left(\frac{r}{\epsilon R}\right)\right]\frac{1}{r^k}\left(r + \frac{\epsilon}{r}\right)^4\right),\ \text{and}\cr
b &=\opO\left(\left[1+\opLog\left(\frac{r}{\epsilon R}\right)\right]\frac{1}{r^{k+1}}\left(r + \frac{\epsilon}{r}\right)^4\right).\cr}\eqnum{\nexteqnno[ModifiedJacobiOperatorOfGrimEndErrorTerms]}
$$
\endproclaim
\proclabel{ModifiedJacobiOperatorOfGrimEnd}
\proof Indeed, by \eqnref{AsymptoticFormulaForGrimEndSmallScale},
$$
u_i = \frac{1}{2}x_i + \frac{c\epsilon}{r^2}x_i + \opO\left(\left[1+\opLog\left(\frac{r}{\epsilon R}\right)\right]\frac{1}{r^k}\left(r + \frac{\epsilon}{r}\right)^3\right).
$$
Thus, by \eqnref{AppBasicFormulaInCoordinateCharts},
$$\eqalign{
\mu^2 &= 1 - \left(\frac{r}{2} + \frac{c\epsilon}{r}\right)^2 + \opO\left(\left[1+\opLog\left(\frac{r}{\epsilon R}\right)\right]\frac{1}{r^k}\left(r + \frac{\epsilon}{r}\right)^4\right),\cr
g^{ij} &= \delta_{ij} - \left(\frac{1}{2} + \frac{c\epsilon}{r^2}\right)^2x^ix^j + \opO\left(\left[1+\opLog\left(\frac{r}{\epsilon R}\right)\right]\frac{1}{r^{k+1}}\left(r + \frac{\epsilon}{r}\right)^4\right).\cr}
$$
It follows that
$$
g^{ij}f_{ij} = \Delta f - \left(\frac{1}{2} + \frac{c\epsilon}{r^2}\right)^2 x^ix^jf_{ij} + a^{ij}f_{ij},
$$
where $a=\opO\left([1+\opLog(r/\epsilon R)]r^{-k}(r + \epsilon/r)^4\right)$, and since $r^{-1}(r + \epsilon/r)^4$ bounds $(r+\epsilon/r)^3$,
$$
-2\mu g^{ip}g^{jq}u_{pq}u_if_j = - \left(\frac{1}{2} - \frac{2\epsilon^2 c^2}{r^4}\right)x^i f_i + b^i f_i,
$$
where $b=\opO\left([1+\opLog(r/\epsilon R)]r^{-(k+1)}(r + \epsilon/r)^4\right)$. The result follows.\qed
\newsubhead{The Regular Component}[TheRegularComponent]
Theorem \procref{JIsLinearIsomorphismOverEnds} is derived from Theorem \procref{JIsLinearIsomorphismOverG} by a perturbation argument. First, let $v_p:]0,\infty[\rightarrow\Bbb{R}$ denote the unique solution of \eqnref{GrimSurfaceEquation} which is defined over the whole positive half-line, as in Chapter \headref{TheGrimParaboloid}. Let $u_p$ denote its primitive with initial value $0$ so that its graph is a Grim paraboloid. Let $\hat{J}_p$ denote its modified MCFS Jacobi operator, as defined in Section \subheadref{TheModifiedMCFSJacobiOperator}. Over the ball $B(2\epsilon R^4)$,
$$
v_p(r) = \frac{1}{2}r + O(r^{3-k}),\eqnum{\nexteqnno[AsymptoticFormulaForGrimParaboloidSmallScale]}
$$
so that, as in Lemma \procref{ModifiedJacobiOperatorOfGrimEnd}, over $B(0,2\epsilon R^4)$,
$$
\hat{J}_p f = \Delta f - \frac{1}{2}x^ix^j f_{ij} - \frac{1}{2}x^if_i + \Cal{E}_pf,\eqnum{\nexteqnno[ModifiedJacobiOperatorOfGrimParaboloid]}
$$
where $\Cal{E}_pf := a^{ij}f_{ij} + b^if_i$ and
$$\eqalign{
a &= \opO(r^{4-k}),\ \text{and}\cr
b &= \opO(r^{3-k}).\cr}\eqnum{\nexteqnno[ModifiedJacobiOperatorOfGrimParaboloidErrorTerms]}
$$
Define
$$
\hat{J}_\gamma := M^{-1}_\gamma\hat{J} M_{\gamma},\eqnum{\nexteqnno[PhiWeightedModifiedJacoiOperator]}
$$
where $M_\gamma$ here denotes the operator of multiplication by $\chi_2 + (1-\chi_2)\phi_\gamma$, $\phi_\gamma$ is given by \eqnref{LargeScaleWeight}, and $\chi_2$ is the cut-off function of the transition region $A(2,4)$ as defined in Appendix \subheadref{GeneralDefinitions}. Observe that, since $\phi_\gamma$ and $\psi$ only depend on $v$ and its integral $u$, it follows by \eqnref{AppMCFSJacobiOperator} that the coefficients of $\hat{J}_\gamma$ are functions of $u$, $v$ and $v_r$ only. Finally, define
$$
\hat{J}_{p,\gamma} := M^{-1}_\gamma\hat{J}_p M_{\gamma}.\eqnum{\nexteqnno[PhiWeightedModifiedJacobiOperatorOfParaboloid]}
$$
A straightforward modification of Theorem \procref{JIsLinearIsomorphismOverG} shows that, for all $\alpha\in]0,1[$, and for all sufficiently small $\gamma$, $\hat{J}_{p,\gamma}$ defines a linear isomorphism from $\Cal{L}^{2,\alpha}_\gamma(G)$ into $\Cal{L}^{0,\alpha}_\gamma(G)$ whose Green's operator has norm uniformly bounded independent of $\Lambda$.
\par
It will suffice to show that the difference $\hat{J}_{p,0}-\hat{J}_0$ converges to $0$ with respect to the hybrid norm as $\Lambda$ tends to $+\infty$. However, the coefficients of this difference diverge making this argument not entirely trivial. Fortunately, they diverge over a region which itself vanishes, and it is this which allows us to conclude. Formally, we define the operators $D$ and $E$ over $A(\epsilon R,\infty)$ by
$$\eqalign{
Df &:= (\hat{J}_0 - E)f - \hat{J}_{p,0}f,\ \text{and}\cr
Ef &:= \chi\frac{2c^2\epsilon^2}{r^4}x^if_i,\cr}\eqnum{\nexteqnno[DefinitionOfDAndE]}
$$
where $\chi$ here denotes the cut-off function of the transiton region $A(\epsilon R^4,2\epsilon R^4)$. We then extend these operators canonically to operators over the whole of $\Bbb{R}^2$, as in Section \subheadref{TheModifiedMCFSJacobiOperator}. By definition,
$$
\hat{J}_0 := \hat{J}_{p,0} + D + E.\eqnum{\nexteqnno[ExtensionOfModifiedJacobiOperator]}
$$
We call $D$ and $E$ the {\bf regular component} and the {\bf singular component} respectively of the difference. We now show that the coefficients of the regular component tend to zero in all norms that concern us as $\Lambda$ tends to infinity. We will study the singular component in the next section.
\par
By \eqnref{ModifiedJacobiOperatorOfGrimEnd}, \eqnref{ModifiedJacobiOperatorOfGrimEndErrorTerms}, \eqnref{ModifiedJacobiOperatorOfGrimParaboloid}, \eqnref{ModifiedJacobiOperatorOfGrimParaboloidErrorTerms} and \eqnref{DefinitionOfDAndE},
$$
Df = a^{ij}f_{ij} + b^if_i,
$$
where, over $A(\epsilon R,2\epsilon R^4)$,
$$\eqalign{
a^{ij} &= -\frac{c\epsilon}{r^2}x^ix^j - \frac{c^2\epsilon^2}{r^4}x^ix^j + \opO\left(\left[1+\opLog\left(\frac{r}{\epsilon R}\right)\right]\frac{1}{r^k}\left(r + \frac{\epsilon}{r}\right)^4\right),\ \text{and}\cr
b^i &= (1-\chi)\frac{2c^2\epsilon^2}{r^4}x^i + \opO\left(\left[1+\opLog\left(\frac{r}{\epsilon R}\right)\right]\frac{1}{r^{k+1}}\left(r + \frac{\epsilon}{r}\right)^4\right).\cr}\eqnum{\nexteqnno[DifferenceOperatorErrorTerms]}
$$
\proclaim{Lemma \nextprocno}
\noindent For sufficiently small $\alpha$,
$$
\|a|_{B(\epsilon R)}\|_{C^{0,\alpha}}, \|b|_{B(\epsilon R)}\|_{C^{0,\alpha}}\rightarrow 0,\eqnum{\nexteqnno[CoefficientsOfDConvergeToZeroOnSmallestScale]}
$$
as $\Lambda$ tends to infinity.
\endproclaim
\proclabel{CoefficientsOfDConvergeToZeroOnSmallestScale}
\proof Indeed, by \eqnref{DifferenceOperatorErrorTerms}, since $\chi$ equals $1$ near $C(\epsilon R)$, over this circle,
$$\eqalign{
a &= \opO\left(\frac{1}{(\epsilon R)^k}\left(\epsilon + \frac{1}{R^2} + (\epsilon R)^4 + \frac{1}{R^4}\right)\right),\ \text{and}\cr
b &= \opO\left(\frac{1}{(\epsilon R)^{k+1}}\left((\epsilon R)^4 + \frac{1}{R^4}\right)\right).\cr}
$$
Since the Lipschitz seminorms of the canonical extensions of $a$ and $b$ over $B(\epsilon R)$ are controlled by their Lipschitz seminorms over $C(\epsilon R)$, by \eqnref{AppFirstInterpolationInequality}, for all $\alpha\in[0,1]$,
$$\eqalign{
\|a|_{B(\epsilon R)}\|_{C^{0,\alpha}}
&\lesssim \frac{\epsilon^{1-\alpha}}{R^\alpha} + \frac{1}{\epsilon^\alpha R^{2+\alpha}}+ (\epsilon R)^{4-\alpha} + \frac{1}{\epsilon^\alpha R^{4+\alpha}},\ \text{and}\cr
\|b|_{B(\epsilon R)}\|_{C^{0,\alpha}}
&\lesssim (\epsilon R)^{3-\alpha} + \frac{1}{\epsilon^{1+\alpha}R^{5+\alpha}}.\cr}
$$
By \eqnref{ControlOfParameters}, for sufficiently small $\alpha$, these both tend to $0$ as $\Lambda$ tends to infinity, as desired.\qed
\proclaim{Lemma \nextprocno}
\noindent For sufficiently small $\alpha$,
$$
\|a|_{A(\epsilon R,2\epsilon R^4)}\|_{C^{0,\alpha}},\ \|b|_{A(\epsilon R,2\epsilon R^4)}\|_{C^{0,\alpha}}\rightarrow 0,\eqnum{\nexteqnno[CoefficientsOfDConvergeToZeroOnSmallScale]}
$$
as $\Lambda$ tends to infinity.
\endproclaim
\proclabel{CoefficientsOfDConvergeToZeroOnSmallScale}
\proof Indeed, by \eqnref{DifferenceOperatorErrorTerms}, over $A(\epsilon R,2\epsilon R^4)$,
$$
a = \opO\left(\frac{1}{r^k}\left(\epsilon + \frac{\epsilon^2}{r^2}\right)\right) + \opO\left(\left[1+\opLog\left(\frac{r}{\epsilon R}\right)\right]\frac{1}{r^k}\left(r + \frac{\epsilon}{r}\right)^4\right),
$$
and $b=b_1+b_2$, where
$$\eqalign{
b_1 &= \opO\left(\left[1+\opLog\left(\frac{r}{\epsilon R}\right)\right]\frac{1}{r^{k+1}}\left(r^4 + \frac{\epsilon^4}{r^4}\right)\right),\ \text{and}\cr
b_2 &= (1-\chi)\frac{2c^2\epsilon^2}{r^4}x^i.\cr}
$$
Thus, by \eqnref{AppFirstInterpolationInequality} and \eqnref{MaximumOfExponential}, for all $\alpha\in[0,1]$,
$$\eqalign{
\|a|_{A(\epsilon R,2\epsilon R^4)}\|_{C^{0,\alpha}}
&\lesssim \frac{\epsilon^{1-\alpha}}{R^\alpha} + \frac{1}{\epsilon^\alpha R^{2+\alpha}} + \opLog(R)(\epsilon R^4)^{4-\alpha} + \frac{1}{\epsilon^\alpha R^{4+\alpha}}\ \text{and},\cr
\|b_1|_{A(\epsilon R,2\epsilon R^4)}\|_{C^{0,\alpha}}
&\lesssim \opLog(R)(\epsilon R^4)^{3-\alpha} + \frac{1}{\epsilon^{1+\alpha}R^{5+\alpha}}.\cr}
$$
By \eqnref{ControlOfParameters}, for sufficiently small $\alpha$, these both tend to $0$ as $\Lambda$ tends to infinity. Finally, over $A(\epsilon R^4,2\epsilon R^4)$,
$$
b_2 = O(\epsilon^2 r^{-(k+3)}),
$$
so that, by \eqnref{AppFirstInterpolationInequality},
$$
\|b_2|_{A(\epsilon R^4,2\epsilon R^4)}\|_{C^{0,\alpha}} \lesssim \frac{1}{\epsilon^{1+\alpha}R^{12+4\alpha}}.
$$
By \eqnref{ControlOfParameters}, for sufficiently small $\alpha$, this also tends to $0$ as $\Lambda$ tends to infinity, and the result follows.\qed
\proclaim{Lemma \nextprocno}
\noindent If $\epsilon R<s<t<\sqrt{2}$, then
$$
\left|v(t) - v_p(t)\right| \leq \left|v(s) - v_p(s)\right|.\eqnum{\nexteqnno[InequalityPreservedOverShortInterval]}
$$
\endproclaim
\proclabel{InequalityPreservedOverShortInterval}
\proof Indeed, by \eqnref{GrimSurfaceEquation}, using a dot to denote diffferentation with respect to $r$, we have
$$
r(\dot{v} - \dot{v}_p) = -(v-v_p)(1 - r(v+v_p) + (v^2 + vv_p + v_p^2)).
$$
However
$$
1 - r(v + v_p) + (v^2 + vv_p + v_p^2) \geq 1 - \frac{r^2}{2}.
$$
Thus, for $r\leq\sqrt{2}$, $\left|v-v_p\right|$ is decreasing, as desired.\qed
\proclaim{Lemma \nextprocno}
\noindent For all $\alpha\in]0,1]$,
$$
\|a|_{A(\epsilon R^4,1)}\|_{C^1}, \|b|_{A(\epsilon R^4,1)}\|_{C^1}\rightarrow 0,\eqnum{\nexteqnno[CoefficientsOfDConvergeToZeroOnLowerCentralScale]}
$$
as $\Lambda$ tends to infinity.
\endproclaim
\proof By \eqnref{AsymptoticFormulaForGrimEndSmallScale} and \eqnref{AsymptoticFormulaForGrimParaboloidSmallScale}, over $C(2\epsilon R^4)$,
$$
\left|v-v_p\right| \lesssim \frac{1}{R^4} + \opLog(R)(\epsilon R^4)^3 + \opLog(R)\frac{1}{R^{12}}.
$$
By Lemma \procref{InequalityPreservedOverShortInterval}, this inequality continues to hold over the whole of $A(2\epsilon R^4,1)$. Since $v$ and $v_p$ both solve \eqnref{GrimSurfaceEquation}, it follows that, over this annulus,
$$
v-v_p = \opO\left(\frac{1}{(\epsilon R^4)^k}\left(\frac{1}{R^4} + \opLog(R)(\epsilon R^4)^3 + \opLog(R)\frac{1}{R^{12}}\right)\right).
$$
Thus,
$$
\|(v-v_p)|_{[2\epsilon R^4,1]}\|_{C^2} \lesssim \frac{1}{\epsilon^2R^{12}} + \opLog(R)\epsilon R^4 + \opLog(R)\frac{1}{\epsilon^2R^{20}}.
$$
so that, by \eqnref{ControlOfParameters},
$$
\|(v-v_p)|_{[2\epsilon R^4,1]}\|_{C^2}\rightarrow 0
$$
as $\Lambda$ tends to infinity. However, by \eqnref{ExplicitFormulaForModifiedJacobiOperator}, over $A(\epsilon R^4,1)$, the coefficients $a$ and $b$ only depend on the first derivatives of $v$ and $v_p$, so that
$$
\|a|_{A(2\epsilon R^4,1)}\|_{C^1},\ \|b|_{A(2\epsilon R^4,1)}\|_{C^1}\rightarrow 0,
$$
as $\Lambda$ tends to infinity, as desired.\qed
\proclaim{Lemma \nextprocno}
\noindent For all $\epsilon>0$, there exists $R>0$ such that if $\left|v(1)-v_p(1)\right|\leq 1$, then
$$
\|a|_{A(R,\infty)}\|_{C^1(G)}, \|b|_{A(R,\infty)}\|_{C^1(G)} \leq \epsilon.\eqnum{\nexteqnno[CoefficientsOfDConvergeToZeroOnLargeScale]}
$$
\endproclaim
\proof Indeed, over $A(4,\infty)$, both $\hat{J}_0$ and $\hat{J}_{p,0}$ are given by \eqnref{FormulaForConjugateOperator}. The result now follows by local uniform dependence of the esimates in \eqnref{CoefficientsOfConjugateOperator} on the initial value.\qed
\proclaim{Lemma \nextprocno}
\noindent For all $R>1$,
$$
\|a|_{A(1,R)}\|_{C^1},\|b|_{A(1,R)}\|_{C^1}\rightarrow 0,\eqnum{\nexteqnno[CoefficientsOfDConvergeToZeroOnUpperCentralScale]}
$$
as $\Lambda$ tends to infinity.
\endproclaim
\proof By \eqnref{AsymptoticFormulaForGrimEndSmallScale}, \eqnref{AsymptoticFormulaForGrimParaboloidSmallScale} and \eqnref{InequalityPreservedOverShortInterval}, over $C(1)$,
$$
\left|v-v_p\right| \lesssim \frac{1}{R^4} + \opLog(R)\left(\epsilon R^4\right)^3 + \opLog(R)\frac{1}{R^{12}}.
$$
Since solutions of first-order ODEs vary smoothly with their parameters,
$$
\|(v-v_p)|_{[1,R]}\|_{C^2} \rightarrow 0,
$$
as $\Lambda$ tends to $\infty$. However, over $A(1,R)$, $a$ and $b$ only depend on $v$ and $v_p$ and their derivatives up to order $2$, and the result follows.\qed
\medskip
\noindent Combining these results yields,
\proclaim{Lemma \nextprocno}
\myitem{(1)} The operator norm of $D$, considered as a map from $H^2(G)$ into $L^2(G)$ converges to zero as $\Lambda$ tends to infinity; and
\medskip
\myitem{(2)} For sufficiently small $\alpha$, the operator norm of $D$, considered as a map from $C^{2,\alpha}(G)$ into $C^{0,\alpha}(G)$ converges to zero as $\Lambda$ tends to infinity.
\endproclaim
\proclabel{TheOperatorDConvergesToZero}
\proof Indeed, by \eqnref{CoefficientsOfDConvergeToZeroOnSmallestScale}, \eqnref{CoefficientsOfDConvergeToZeroOnSmallScale}, \eqnref{CoefficientsOfDConvergeToZeroOnLowerCentralScale}, \eqnref{CoefficientsOfDConvergeToZeroOnLargeScale} and
\eqnref{CoefficientsOfDConvergeToZeroOnUpperCentralScale}, for sufficiently small $\alpha$, both $\|a\|_{C^{0,\alpha}(G)}$ and $\|b\|_{C^{0,\alpha}(G)}$ converge to $0$ as $\Lambda$ tends to infinity, and the result follows.\qed
\newsubhead{The Singular Component}[TheSingularPart]
\noindent We now write
$$
Ef =: a^if_i.\eqnum{\nexteqnno[DefinitionOfAForOperatorE]}
$$
Since $E$ is defined by canonical extension, over the ball $B(\epsilon R)$,
$$
a^i = \frac{2c^2}{\epsilon^2 R^4}x^i.\eqnum{\nexteqnno[FormulaForAOverCentralBall]}
$$
At this stage we require the following key estimate, which reveals the significance of the hybrid norm.
\proclaim{Lemma \nextprocno}
\noindent For sufficiently small $\alpha$ and for sufficiently small $\gamma$,
$$
\|f\|_{C^{1,\alpha}_\gamma(G)}\lesssim (\epsilon R)^{1-2\alpha}\|f\|_{2,\alpha,\gamma}.\eqnum{\nexteqnno[HybridProperty]}
$$
\endproclaim
\proclabel{HybridProperty}
\remark It will be useful to observe that this relation is also valid for spaces of functions defined over an unbounded annulus.
\medskip
\proof Indeed, by the Sobolev embedding theorem, for all $\beta<1$,
$$
\|f\|_{C^{0,\beta}_\gamma(G)} \lesssim \|f\|_{H^2_\gamma(G)} \lesssim (\epsilon R)\|f\|_{2,\alpha,\gamma}.
$$
Setting $\beta=(1-\alpha)$ and using \eqnref{AppFirstInterpolationInequality} and \eqnref{AppSecondInterpolationInequality}, we obtain
$$
\|f\|_{C^{1,\alpha}_\gamma(G)} \lesssim (\epsilon R)^{\frac{1}{1+2\alpha}}\|f\|_{2,\alpha,\gamma} \lesssim (\epsilon R)^{1-2\alpha}\|f\|_{2,\alpha,\gamma},
$$
as desired.\qed
\proclaim{Lemma \nextprocno}
\noindent For sufficiently small $\alpha\in[0,1]$ and for sufficiently small $\gamma$, the operator norm of $E$, considered as a map from $\Cal{L}^{2,\alpha}_\gamma(G)$ into $C^{0,\alpha}_\gamma(G)$ tends to $0$ as $\Lambda$ tends to infinity.
\endproclaim
\proof Indeed, over $A(\epsilon R,2\epsilon R^4)$,
$$
a^i = \opO\bigg(\frac{\epsilon^2}{r^{3+k}}\bigg),
$$
so that
$$\eqalign{
\|a^i|_{A(\epsilon R,2\epsilon R^4)}\|_{C^0} &\lesssim \frac{1}{\epsilon R^3},\ \text{and}\cr
[a^i|_{A(\epsilon R,2\epsilon R^4)}]_1 &\lesssim \frac{1}{\epsilon^2 R^4}.\cr}
$$
Since $a^i$ is extended canonically over $B(\epsilon R)$, these inequalities also hold over the whole of $B(2\epsilon R^4)$ so that, by \eqnref{AppFirstInterpolationInequality}, for all $\alpha\in[0,1]$,
$$
[a^i]_\alpha \lesssim \frac{1}{(\epsilon R)^\alpha\epsilon R^3}.
$$
It follows by \eqnref{HybridProperty} and \eqnref{AppProductRule} that
$$
\|Ef\|_{C^{0,\alpha}_\gamma(G)} \lesssim \frac{1}{(\epsilon R)^\alpha\epsilon R^3}\|f\|_{C^{1,\alpha}_\gamma(G)},
$$
and the result follows by \eqnref{ControlOfParameters}.\qed
\proclaim{Lemma \nextprocno}
\noindent For sufficiently small $\alpha\in[0,1]$ and for sufficiently small $\gamma$, the operator norm of $(\epsilon R)^{-1}E$ considered as a map from $\Cal{L}^{2,\alpha}_\gamma(G)$ into $H^0_\gamma(G)$ tends to $0$ as $\Lambda$ tends to infinity.
\endproclaim
\proof Indeed, a direct calculation yields
$$
\|a^i\|_{L^2_\gamma(G)}\lesssim \frac{1}{R^2}.
$$
Thus, bearing in mind \eqnref{HybridProperty},
$$\eqalign{
\|(\epsilon R)^{-1}Ef\|_{L^2_\gamma(G)}
&\lesssim (\epsilon R)^{-1}\|a^i\|_{L^2_\gamma(G)}\|Df\|_{L^\infty(G)}\cr
&\lesssim (\epsilon R)^{-1}\|a^i\|_{L^2_\gamma(G)}\|f\|_{C^{1,\alpha}_\gamma(G)}\cr
&\lesssim \frac{1}{(\epsilon R)^{2\alpha}R^2}\|f\|_{2,\alpha,\gamma},\cr}
$$
and the result follows by \eqnref{ControlOfParameters}.\qed
\medskip
\noindent Combining these results yields
\proclaim{Lemma \nextprocno}
\noindent For sufficiently small $\alpha\in[0,1]$ and for sufficiently small $\gamma$, the operator norm of $E$ considered as a map from $\Cal{L}^{2,\alpha}_\gamma(G)$ into $\Cal{L}^{0,\alpha}_\gamma(G)$ tends to $0$ as $\Lambda$ tends to infinity.
\endproclaim
\proclabel{TheOperatorEConvergesToZero}
\newhead{Surgery and the Perturbation Family}[SurgeryAndThePerturbationFamily]
\newsubhead{The Basic Surgery Operation}[TheBasicSurgeryOperation]
Recall that our strategy for proving Theorem A consists of two stages. The first involves a surgery operation in which approximate MCF solitons are constructed out of properly embedded minimal surfaces and rotationally symmetric Grim ends. The second involves a fixed-point argument in which these approximate MCF solitons are perturbed into actual MCF solitons. In this section, we describe the surgery operation and in Section \subheadref{TheDeformationFamily}, we describe the family of deformations of the approximate MCF soliton in which the actual MCF soliton will be found. Though conceptually simple, our construction is inevitably rather technical. However, we believe that a careful reading of the following two sections will be rewarded by a clear understanding of the essence of this paper.
\par
Consider first a properly embedded surface $C$ in $\Bbb{R}^3$, minimal outside of some compact set, and with finitely many ends, all of which are horizontal. Let $R_0>0$ be such that every component of $C\minter(A(R_0,\infty)\times\Bbb{R})$ is a minimal graph over $A(R_0,\infty)$. Let $F:A(R_0,\infty)\rightarrow\Bbb{R}$ be the profile of one of these minimal ends, so that
$$
F = a + c\opLog(r) + \opO\left(r^{-(1+k)}\right),
$$
for some real constants $a$ and $c$, which will henceforth be referred to respectively as the {\bf constant term} and the {\bf logarithmic parameter} of the minimal end. In particular, planar ends are simply catenoidal ends with vanishing logarithmic parameters. We will only be concerned with minimal ends invariant under reflection in at least two distinct vertical planes. We therefore assume that $F$ contains no terms of order $(-1)$, so that
$$
F = a + c\opLog(r) + \opO\left(r^{-(2+k)}\right).\eqnum{\nexteqnno[ProfileOfCatenoidalEnd]}
$$
This asymptotic formula will be used repeatedly throughout the sequel.
\par
Let $\Lambda$ be a large, positive number, let $K>0$ be a fixed constant, and choose $\epsilon,R>0$ and $\left|c\right|<K$ as in \eqnref{ControlOfParameters}. Let $G:A(R/4,\infty)\rightarrow\Bbb{R}$ be the profile of a rotationally symmetric Grim end with constant term $a$, logarithmic parameter $c$ and speed $\epsilon$. Rescaling and integrating \eqnref{ErrorForExactSolutionWithGivenLogarithmicParameter} we obtain, over the annulus $A(R/4,2R^4)$,
$$
G = a + c\opLog(r) + \frac{1}{4}\epsilon r^2 + \opO\left(\left[1+\opLog\left(\frac{r}{R}\right)\right]r^{1-k}\left(\epsilon r + \frac{1}{r}\right)^3\right).\eqnum{\nexteqnno[ProfileOfGrimEnd]}
$$
Let $\chi_c$ be the cut-off function of the central transition region $A(R,2R)$, as defined in Appendix \subheadref{GeneralDefinitions}, and define the function $H$ over $A(R_0,\infty)$ by
$$
H := \chi_c F + (1-\chi_c)G.\eqnum{\nexteqnno[DefinitionOfJoinedSurface]}
$$
Its graph will be called the {\bf joined end}. Observe that $H$ is entirely determined by $F$ and the parameters $\epsilon$ and $R$. Furthermore, over the annuli $A(R_0,R)$ and $A(2R,\infty)$, $H$ simply coincides with $F$ and $G$ respectively whilst, over the annulus $A(R,2R)$, by \eqnref{ProfileOfCatenoidalEnd}, \eqnref{ProfileOfGrimEnd} and the fact that $\chi_c=O(r^{-k})$
$$
H = a + c\opLog(r) + \frac{1}{4}\epsilon(1-\chi_c)r^2 + O\left(r^{-(2+k)}\right).\eqnum{\nexteqnno[FormulaForJoinedSurface]}.
$$
\newsubhead{The Deformation Family}[TheDeformationFamily]
Continuing to use the notation of Section \subheadref{TheBasicSurgeryOperation}, let $S$ denote the surface obtained by replacing each of the ends of $C$ with their respective joined ends. We now construct a family of deformations of $S$ out of which the actual MCF soliton will be selected when $\Lambda$ is large. We first describe how the logarithmic parameters of $C$ and $S$ are varied. Let $n$ denote the number of ends of $C$, and for each $1\leq i\leq n$, let $a_{0,i}$ and $c_{0,i}$ be respectively the constant term and the logarithmic parameter of the $i$'th end. Let $U$ be a neighbourhood of $(c_{0,1},...,c_{0,n})$ in $\Bbb{R}^n$ and let $(C_c)_{c\in U}$ be a smoothly varying family of immersed surfaces in $\Bbb{R}^3$ such that $C_{c_0}=C$ and, for all $c\in U$ and for all $1\leq i\leq n$, the $i$'th component of $C_c\minter(A(R_0,\infty)\times\Bbb{R})$ is a horizontal, minimal end with constant term $a_{0,i}$ and logarithmic parameter $c_i$. Finally, for all $c\in U$, let $S_c$ denote the surface obtained by replacing each end of $C_c$ with its corresponding joined end, as described in Section \subheadref{TheBasicSurgeryOperation}.
\par
Let $E:U\times S\rightarrow\Bbb{R}^3$ be a smooth function such that
\medskip
\myitem{(1)} for all $c\in U$, $E_c:=E(c,\cdot)$ parametrises $S_c$; and
\medskip
\myitem{(2)} for all $c\in U$, and for all $p\in S\minter(A(R_0,+\infty)\times\Bbb{R})$, the point $E_c(p)$ lies vertically above or below the point $p$.
\medskip
\noindent Let $\chi_0$, $\chi'_0$, $\chi'_\epsilon$ and $\chi_\epsilon$ be the cut-off functions of the transition regions $A(R_0,2R_0)$, $A(2R_0,4R_0)$, $A(1/2\epsilon,1/\epsilon)$ and $A(1/\epsilon,2/\epsilon)$ respectively, as defined in Appendix \subheadref{GeneralDefinitions}. By composing with vertical projections onto $\Bbb{R}^2$, we think of these functions also as functions defined over $S$. For all $c\in  U$, let $N_c$ denote the unit normal vector field over $S_c$. For all $1\leq i\leq n$, let $\Bbb{I}_i:S\rightarrow\left\{0,1\right\}$ denote the indicator function of the $i$'th component of $S_c\minter(A(R_0,\infty)\times\Bbb{R})$. Observe that, since this intersection is a union of graphs, every component is transverse to the unit vertical vector $e_z$. For all $1\leq i\leq n$, let $\epsilon_i\in\left\{\pm 1\right\}$ be such that $\epsilon_ie_z$ lies on the same side of the $i$'th component as $N_c$. For all $c\in U$, define the {\bf modified normal vector field} over $S_c$ by
$$
\hat{N}_c := (\chi_\epsilon-\chi_0)\epsilon_i e_z + (1-(\chi_\epsilon-\chi_0))N_c.\eqnum{\nexteqnno[ModifiedNormalVectorField]}
$$
Observe that, over the regions $S_c\minter(B(R_0)\times\Bbb{R})$ and $S_c\minter(A(2/\epsilon,\infty)\times\Bbb{R})$, this vector field coincides with $N_c$ whilst, over the region $S_c\minter(A(2R_0,1/\epsilon)\times\Bbb{R})$, it coincides with $\pm e_z$. Now let $V$ and $W$ be neighbourhoods of $0$ in $\Bbb{R}^n$ and define $\tilde{E}:U\times V\times W\times C^\infty(S)\rightarrow C^\infty(S,\Bbb{R}^3)$ by
$$
\tilde{E}_{c,a,b,f}(p) := E_c(p) + f(p)\hat{N}_c(p) + \sum_{i=1}^n\epsilon_i\Bbb{I}_i(p)\left(a_i\left(1-\chi_0'(p)\right) + b_i\left(1-\chi_\epsilon'(p)\right)\right)e_z.\eqnum{\nexteqnno[PerturbationFamily]}
$$
Upon reducing $U$, $V$ and $W$ if necessary, there trivially exists $\delta>0$, which is independent of $\Lambda$, $\epsilon$ and $R$, such that, for all $(c,a,b)\in U\times V\times W$, and for all $\|f\|_{C^0}<\delta$, the function $\tilde{E}_{c,a,b,f}$ defines an immersion of $S$ into $\Bbb{R}^3$. This concludes the description of the deformation family in which the actual MCF solution will be found.
\newsubhead{Microscopic and Macroscopic Perturbations}[MicroscopicAndMacroscopicPerturbations]
Continuing to use the notation of Sections \subheadref{TheBasicSurgeryOperation} and \subheadref{TheDeformationFamily}, we consider now the first-order perturbations of $S$ defined by the above deformation family. We classify these perturbations into two main types. Those in the direction of $C^\infty(S)$ will be called {\bf microscopic perturbations}, and those in the directions of $U$, $V$ and $W$ will be called {\bf macroscopic perturbations}. We now describe the first-order variations of the MCFS functional resulting from macroscopic perturbations. The first-order variations resulting from microscopic perturbations will be studied in the next section.
\par
Recall that, as in Appendix \subheadref{SurfaceGeometry}, the MCFS functional with speed $\epsilon$ of an immersion $E:S\rightarrow\Bbb{R}^3$ is given by
$$
M_E := H_E + \epsilon\langle N_E,e_z\rangle,\eqnum{\nexteqnno[MCFSFunctionalAgain]}
$$
where $H_E$ here denotes the mean curvature function of $E$, and $N_E$ here denotes its unit normal vector field. We define $M_\epsilon:U\times V\times W\rightarrow C_0^\infty(S)$ such that, for all $(c,a,b)\in U\times V\times W$, and for all $p\in S$, $M_{\epsilon,c,a,b}(p)$ is the value of this functional for the immersion $E_{c,a,b}$ at the point $p$. We define the operators $X_\epsilon,Y_\epsilon,Z_\epsilon:\Bbb{R}^d\rightarrow C^\infty_0(S)$ by
$$\eqalign{
(X_\epsilon u)(p) &:= \frac{1}{\langle\hat{N}_S,N_S\rangle}\frac{d}{dt}M_{\epsilon,c_0+tu,0,0}(p)\bigg|_{t=0},\cr
(Y_\epsilon v)(p) &:= \frac{1}{\langle\hat{N}_S,N_S\rangle}\frac{d}{dt}M_{\epsilon,c_0,tv,0}(p)\bigg|_{t=0},\ \text{and}\cr
(Z_\epsilon w)(p) &:= \frac{1}{\langle\hat{N}_S,N_S\rangle}\frac{d}{dt}M_{\epsilon,c_0,0,tw}(p)\bigg|_{t=0}.\cr}\eqnum{\nexteqnno[PerturbationOperators]}
$$
These are the first-order variations of the MCFS functional arising from the $3$ types of macroscopic perturbation. In particular, since $M_{\epsilon,c,0,0}$ vanishes over $S\minter(A(2R,+\infty)\times\Bbb{R})$ for all $c\in V$, for all $u\in\Bbb{R}^d$, $Xu$ is supported over $S\minter (B(2R)\times\Bbb{R})$. Likewise, for all $v,w\in\Bbb{R}^n$, $Yv$ and $Zw$ are supported over $S\minter(A(2R_0,4R_0)\times\Bbb{R})$ and $S\minter(A(1/2\epsilon,1/\epsilon)\times\Bbb{R})$ respectively. In later sections, when no ambiguity arises, the subscript $\epsilon$ will be suppressed, and these operators will be denoted simply by $X$, $Y$ and $Z$ respectively.
\newsubhead{Modified Jacobi Operators}[ModifiedJacobiOperators]
The operator of first-order variation of the MCFS functional resulting from microscopic perturbations is none other than the modified MCFS Jacobi operator. In this section, we determine asymptotic formulae for its coefficients over different regions. We recall that, since different modifications are made on different scales, the precise definition of the modified MCFS Jacobi operator varies with context. We now describe the framework which unifies these different definitions. We will then study three different cases corresponding to, in order, the scale of joined surfaces, the scale of CHM surfaces, and the scale of rotationally symmetric Grim ends.
\par
Consider first a general immersed surface $\Sigma$ in $\Bbb{R}^3$ such that, for some $R_0>0$, every component of $\Sigma\minter(A(R_0,\infty)\times\Bbb{R})$ is a graph over $A(R_0,\infty)$. Let $\Lambda>0$ be a large, positive number, let $\epsilon,R>0$ be as in \eqnref{ControlOfParameters}, and let $\hat{N}_\Sigma$ be the {\bf modified normal vector field} over $\Sigma$ as defined in \eqnref{ModifiedNormalVectorField}. We define $E:C_0^\infty(\Sigma)\rightarrow C^\infty(\Sigma,\Bbb{R}^3)$ by
$$
E_f(p) := p + f(p)\hat{N}_\Sigma(p).
$$
Observe that if $f$ is sufficiently small, then $E_f$ is an immersion. Define $M:C_0^\infty(\Sigma)\rightarrow C^\infty(\Sigma)$ such that, for all such $f$, and for all $p\in\Sigma$, $M_f(p)$ is the value of the MCFS functional \eqnref{MCFSFunctionalAgain} with speed $\epsilon$ for the immersion $E_f$ at the point $p$. The {\bf modified MCFS Jacobi operator} of $\Sigma$ with speed $\epsilon$ is now defined by
$$
(\hat{J}_{\Sigma,\epsilon} f)(p) := \frac{1}{\langle\hat{N}_\Sigma,N_\Sigma\rangle}\frac{d}{dt}M_{tf}(p)\bigg|_{t=0}.\eqnum{\nexteqnno[GeneralModifiedJacobiOperator]}
$$
In later sections, when no ambiguity arises, the subscript $\epsilon$ will be suppressed, and this operator will be denoted simply by $\hat{J}_\Sigma$.
\par
Over the annulus $A(R/4,1/\epsilon)$, since $\hat{N}_\Sigma$ here coincides with $e_z$, the operator $\hat{J}_{\Sigma,\epsilon}$ is simply $\langle N_\Sigma,e_z\rangle^{-1}$ times the linearisation of the MCFS functional for graphs. Consequently, if $F:A(R/4,1/\epsilon)\rightarrow\Bbb{R}$ is the profile of a component of $\Sigma\minter(A(R/4,1/\epsilon)\times\Bbb{R})$ then, upon differentiating \eqnref{AppMCFSFunctionalForGraphs} we obtain, over this annulus,
$$
\hat{J}_{\Sigma,\epsilon}f = g^{ij}f_{ij} - \mu^2g^{ij}F_{ij}F_kf_k + 2\mu^4F_iF_jF_kF_{ij}f_k - 2\mu^2F_{ij}F_if_j - \epsilon\mu^2F_i f_i.\eqnum{\nexteqnno[ExplicitFormulaOfJacobiOperatorForGraphs]}
$$
In particular, for all $v,w\in\Bbb{R}^n$, and for all $p\in S$,
$$\eqalign{
(Yv)(p) &= -\sum_{i=1}^n\Bbb{I}_i(p)v_i(\hat{J}_{\Sigma,\epsilon}\chi_0')(p),\ \text{and}\cr
(Zw)(p) &= -\sum_{i=1}^n\Bbb{I}_i(p)w_i(\hat{J}_{\Sigma,\epsilon}\chi_\epsilon')(p).}\eqnum{\nexteqnno[ValuesOfYVAndZW]}
$$
\par
Now let $C$ be a minimal end over the annulus $A(R_0,\infty)$ satisfying \eqnref{ProfileOfCatenoidalEnd} and let $\hat{J}_{C,\epsilon}$ be its modified MCFS Jacobi operator with speed $\epsilon$.
\proclaim{Lemma \nextprocno}
\noindent Over $A(R/4,2R^4)$,
$$
\hat{J}_{C,\epsilon} f = \Delta f - \frac{c^2}{r^4}x^ix^jf_{ij} - \frac{\epsilon c}{r^2}x^if_i + \frac{2c^2}{r^4}x^if_i + \Cal{E}_{C,\epsilon} f,\eqnum{\nexteqnno[JacobiOperatorOfCHM]}
$$
where $\Cal{E}_{C,\epsilon}f:=a^{ij}f_{ij} + b^if_i$ and $a$ and $b$ satisfy
$$\eqalign{
a &= \opO\left(r^{-(k+4)}\right),\ \text{and}\cr
b &= \opO\left(r^{-(k+4)}\left(\epsilon r + \frac{1}{r}\right)\right).\cr}\eqnum{\nexteqnno[ErrorTermsInJacobiOperatorOfCHM]}
$$
\endproclaim
\proof By \eqnref{ProfileOfCatenoidalEnd},
$$
F_i = \frac{c}{r^2}x^i + \opO\left(r^{-(k+3)}\right).
$$
Thus, by \eqnref{AppBasicFormulaInCoordinateCharts},
$$\eqalign{
\mu^2 &= 1 - \frac{c^2}{r^2} + \opO\left(r^{-(k+4)}\right),\ \text{and}\cr
g^{ij} &= \delta_{ij} - \frac{c^2}{r^4}x^ix^j + \opO\left(r^{-(k+4)}\right).\cr}
$$
Thus,
$$
g^{ij}f_{ij} = \Delta f - \frac{c^2}{r^4}x^ix^j f_{ij} + a^{ij}f_{ij},
$$
where $a=\opO\left(r^{-(k+4)}\right)$. Likewise,
$$\eqalign{
\mu^2 g^{ij} F_{ij}F_kf_k &= b_1^if_i,\cr
2\mu^4 F_iF_jF_kF_{ij}f_k &= b_2^if_i,\ \text{and}\cr
-2\mu^2F_{ij}F_if_j &= \frac{2c^2}{r^4}x^if_i + b_3^if_i,\cr}
$$
where $b_1^i, b_2^i, b_3^i = \opO\left(r^{-(k+5)}\right)$. Finally,
$$
\epsilon\mu^2 F_i f_i = \frac{\epsilon c}{r^2}x^if_i + b_4^i,
$$
where $b_4^i = \opO\left(\epsilon r^{-(k+3)}\right)$. The result follows.\qed
\medskip
Next let $G$ be a rotationally symmetric Grim end of speed $\epsilon$ over the annulus $A(R/4,\infty)$ and let $\hat{J}_{G,\epsilon}$ be its modified MCFS Jacobi operator with speed $\epsilon$. Define $\psi:G\rightarrow\Bbb{R}$ by
$$
\psi := \langle\hat{N}_G,N_G\rangle = \chi_\epsilon \langle e_z,N_G\rangle + (1-\chi_\epsilon),\eqnum{\nexteqnno[DefinitionOfWeightPsi]}
$$
and denote by $M_\psi$ the operator of multiplication by $\psi$.
\proclaim{Lemma \nextprocno}
\noindent Over $A(R/4,\infty)$,
$$
\hat{J}_{G,\epsilon} := M_\psi^{-1}J_{G,\epsilon}M_\psi,\eqnum{\nexteqnno[ModifiedMCFSJacobiOperatorOfG]}
$$
where $J_{G,\epsilon}$ denotes the MCFS Jacobi operator with speed $\epsilon$ of $G$, as defined in Appendix \subheadref{SurfaceGeometry}.
\endproclaim
\proclabel{DefnsOfModifiedMCFSJacobiOpMatch}
\remark In particular, in the case of rotationally symmetric Grim ends, the modified MCFS Jacobi operator as defined above coincides, up to rescaling, with the modified MCFS Jacobi operator as defined in Section \subheadref{TheModifiedMCFSJacobiOperator}.
\medskip
\proof Indeed, more generally, with $M:=M_0$ defined as at the beginning of this section, for all $f\in C_0^\infty(\Sigma)$,
$$
\hat{J}_{\Sigma,\epsilon}f = M_\psi^{-1}J_{\Sigma,\epsilon} M_\psi f + M_\psi^{-1}\langle X,\nabla M\rangle f,
$$
where $X$ here denotes the tangential component of the vector field $\hat{N}_\Sigma$. The result now follows since $M$ vanishes identically over $G$.\qed
\medskip
\noindent In particular, rescaling \eqnref{ModifiedJacobiOperatorOfGrimEnd} and \eqnref{ModifiedJacobiOperatorOfGrimEndErrorTerms} immediately yields
\proclaim{Lemma \nextprocno}
\noindent Over $A(R/4,\infty)$,
$$
\hat{J}_{G,\epsilon}f = \Delta f - \left(\frac{\epsilon}{2}+\frac{c}{r^2}\right)^2x^ix^jf_{ij} - \left(\frac{\epsilon^2}{2}-\frac{2c^2}{r^4}\right)x^i f_i + \Cal{E}_G f.\eqnum{\nexteqnno[JacobiOperatorOfRescaledGrimEnd]}
$$
where $\Cal{E}_{G,\epsilon} f := a^{ij}f_{ij} + b^if_i$, and $a$ and $b$ satisfy,
$$\eqalign{
a &= \opO\left(\left[1+\opLog\left(\frac{r}{R}\right)\right]\frac{1}{r^k}\left(\epsilon r + \frac{1}{r}\right)^4\right),\ \text{and}\cr
b &= \opO\left(\left[1+\opLog\left(\frac{r}{R}\right)\right]\frac{1}{r^{k+1}}\left(\epsilon r + \frac{1}{r}\right)^4\right).\cr}\eqnum{\nexteqnno[ErrorTermsInJacobiOperatorOfRescaledGrimEnd]}
$$
\endproclaim
Finally, let $S$ be a joined end, as constructed in Section \subheadref{TheBasicSurgeryOperation}, and let $\hat{J}_{S,\epsilon}$ denote its modified MCFS Jacobi operator with speed $\epsilon$.
\proclaim{Lemma \nextprocno}
\noindent Over $A(R,2R)$,
$$\eqalign{
\left(\hat{J}_{S,\epsilon} - \hat{J}_{C,\epsilon}\right)f &= a_1^{ij}f_{ij} + b_1^if_i,\ \text{and}\cr
\left(\hat{J}_{S,\epsilon} - \hat{J}_{G,\epsilon}\right)f &= a_2^{ij}f_{ij} + b_2^if_i.\cr}
$$
where $a_1$, $a_2$, $b_1$ and $b_2$ satisfy,
$$\eqalign{
a_1,a_2 &= \opO\left(r^{-(4+k)}\right),\ \text{and}\cr
b_1,b_2 &= \opO\left(r^{-(5+k)}\right).\cr}\eqnum{\nexteqnno[ErrorFromJacobiOperatorOfJoinedSurface]}
$$
\endproclaim
\proof By \eqnref{ProfileOfCatenoidalEnd}, \eqnref{FormulaForJoinedSurface} and \eqnref{ControlOfParameters}, over $A(R,2R)$,
$$\eqalign{
H_i - F_i &= \opO\left(r^{-(3+k)}\right),\ \text{and}\cr
F_i, H_i &= \opO\left(r^{-(1+k)}\right).\cr}
$$
Thus, by \eqnref{AppBasicFormulaInCoordinateCharts},
$$\eqalign{
\mu_H - \mu_F &= \opO\left(r^{-(4+k)}\right),\ \text{and}\cr
g^{ij}_H - g^{ij}_F &= \opO\left(r^{-(4+k)}\right),\cr}
$$
The result follows for $(\hat{J}_{S,\epsilon}-\hat{J}_{C,\epsilon})$ by \eqnref{ExplicitFormulaOfJacobiOperatorForGraphs}. The result for $(\hat{J}_{S,\epsilon}-\hat{J}_{G,\epsilon})$ follows in a similar manner, and this completes the proof.\qed
\medskip
We conclude this section by studying commutators of modified Jacobi operators with certain multiplication operators. Indeed, let $[\hat{J}_{C,\epsilon},\chi_l]$ denote the commutator of $\hat{J}_{C,\epsilon}$ with the operator of multiplication by the cut-off function $\chi_l$ of the lower transition region $A(R/4,R/2)$. Likewise, let $[\hat{J}_{G,\epsilon},\chi_u]$ denote the commutator of $\hat{J}_{G,\epsilon}$ with the operator of multiplication by the cut-off function $\chi_u$ of the upper transition region $A(R^4,2R^4)$. Observe that these operators are supported over the annuli $A(R/4,R/2)$ and $A(R^4,2R^4)$ respectively.
\noskipproclaim{Lemma \nextprocno}
$$\eqalign{
\left[\hat{J}_{C,\epsilon},\chi_l\right]f &= a_1^if_i + b_1f,\ \text{and}\cr
\left[\hat{J}_{G,\epsilon},\chi_u\right]f &= a_2^if_i + b_2f,\cr}
$$
where $a_1$, $a_2$, $b_1$ and $b_2$ satisfy,
$$\eqalign{
a_1, a_2 &= \opO\left(r^{-(k+1)}\right),\ \text{and}\cr
b_1, b_2 &= \opO\left(r^{-(k+2)}\right).\cr}\eqnum{\nexteqnno[ErrorTermsInCommutators]}
$$
\endproclaim
\proclabel{ErrorTermsInCommutators}
\proof Indeed, observe that $\chi_l,\chi_u=\opO(r^{-k})$. The result now follows by \eqnref{JacobiOperatorOfCHM}, \eqnref{ErrorTermsInJacobiOperatorOfCHM}, \eqnref{JacobiOperatorOfRescaledGrimEnd} and \eqnref{ErrorTermsInJacobiOperatorOfRescaledGrimEnd}.\qed
\newsubhead{Controlling Macroscopic Perturbations}[ControllingMacroscopicPerturbations]
We complete this chapter by studying the first-order variation of the MCFS functional resulting from the first macroscopic perturbation. Recall that, for all $u\in\Bbb{R}^n$, $Xu$ vanishes outside $B(2R)$. Inside this ball, we have
\proclaim{Lemma \nextprocno}
\noindent For $u\in\Bbb{R}^d$ such that $\|u\|=1$, over $A(2R_0,R)$,
$$
Xu = \opO\left(\epsilon r^{-(2+k)}\right),\eqnum{\nexteqnno[NormOfDeficiencyTermInLowerRegion]}
$$
and over $A(R,2R)$,
$$
Xu = \opO\left(r^{-(4+k)}\right).\eqnum{\nexteqnno[NormOfDeficiencyTermOfJoinedSurface]}
$$
\endproclaim
\proof For notational convenience, we suppose that $C$ and $S$ each only have one end and, in particular, that $u=1$. Let $C_c$ and $S_c$ be smooth families of immersed surfaces as in Section \subheadref{TheBasicSurgeryOperation}. For all $t$, let $F_t:A(2R_0,\infty)\rightarrow\Bbb{R}$ and $H_t:A(2R_0,\infty)\rightarrow\Bbb{R}$ denote the profiles of $C_{c_0+t}\minter(A(2R_0,\infty)\times\Bbb{R})$ and $S_{c_0+t}\minter(A(2R_0,\infty)\times\Bbb{R})$ respectively. Denote
$$\eqalign{
Z &:= \frac{d}{dt}F_t|_{t=0},\ \text{and}\cr
W &:= \frac{d}{dt}H_t|_{t=0},\cr}
$$
and observe that, over $A(2R_0,2R)$,
$$
Xu = \hat{J}_{S,\epsilon} W.
$$
Now, by \eqnref{ProfileOfCatenoidalEnd},
$$
Z=\opLog(r) + \opO\left(r^{-(2+k)}\right).
$$
Next, by \eqnref{ErrorForExactJacobiFieldWithGivenLogarithmicParameter} and \eqnref{DefinitionOfJoinedSurface}, and bearing in mind that $\chi_c=\opO(r^{-k})$, over $A(R,2R)$, we have
$$
W = \opLog(r) + \opO\left(r^{-(2+k)}\right) = Z + \opO\left(r^{-(2+k)}\right),\eqnum{\nexteqnno[DifferenceBetweenZAndW]}
$$
and since $Z=W$ over $A(2R_0,R)$, \eqnref{DifferenceBetweenZAndW} in fact holds over the whole of $A(2R_0,2R)$. We now write
$$
Xu = \hat{J}_{C,\epsilon}Z + \left(\hat{J}_{S,\epsilon} - \hat{J}_{C,\epsilon}\right)Z + \hat{J}_{S,\epsilon}(W-Z).
$$
The second and third terms are supported over $A(R,2R)$, and by \eqnref{JacobiOperatorOfCHM} and \eqnref{ErrorFromJacobiOperatorOfJoinedSurface},
$$\eqalign{
\left(\hat{J}_{S,\epsilon} - \hat{J}_{C,\epsilon}\right)Z &= \opO\left(r^{-(6+k)}\right),\ \text{and}\cr
\hat{J}_{S,\epsilon}(W-Z) &= \opO\left(r^{-(4+k)}\right).\cr}
$$
Finally, since the graph of $F_t$ is minimal for all $t$, by \eqnref{AppBasicFormulaInCoordinateCharts} and \eqnref{AppMCFSFunctionalForGraphs},
$$
\hat{J}_{C,\epsilon}Z = -\epsilon\mu^2F_{0,i}Z_i = \opO\left(\epsilon r^{-(2+k)}\right),
$$
and the result follows by \eqnref{ControlOfParameters}.\qed
\newhead{Constructing the Green's Operator}[ConstructingTheGreensOperator]
\newsubhead{The Cylindrical, Grim and Hybrid Norms}[TheCylindricalAndGrimNorms]
We now prepare the ground for the perturbation argument that will be used to construct actual MCF solitons out of the approximate MCF solitons constructed in Section \subheadref{TheBasicSurgeryOperation}. In this chapter, we construct the Green's operator of the modified MCFS Jacobi operator of the approximate MCF soliton together with estimates of its operator norm. It is the determination of suitable estimates, requiring a careful and lengthy analysis, which constitutes the hardest part of this paper. We will see presently that sufficiently strong estimates are made possible by the correct choice of functional norms over the different components of the approximate MCF soliton, as well as the use of the hybrid norm, already mentioned in Chapter \headref{GrimEnds}. Throughout this section, we will make use of \eqnref{ControlOfParameters} without comment.
\par
We first study the analytic properties of Green's operators over CHM surfaces. Thus, for $g$ a positive integer, let $C:=C_g$ be the CHM surface of genus $g$. Observe that functions over $C\minter(A(R_0,\infty)\times\Bbb{R})$ may be considered as functions over three copies of $A(R_0,\infty)$. In defining norms over spaces of functions, we will pass between these two perspectives without comment. Consider now the triplet $(X,Y,\hat{J}_C)$ where $X$ and $Y$ are the operators constructed in Section \subheadref{TheBasicSurgeryOperation} and $\hat{J}_C$ is the modified MCFS Jacobi operator of $C$ as constructed in Section \subheadref{ModifiedJacobiOperators}. We now construct a right inverse for this operator when $\Lambda$ is large. We first gather various basic results that will be of use to us. Let $D$ denote the total differentiation operator over $\Bbb{R}^2$ and denote
$$
D_\opCyl := rD,\eqnum{\nexteqnno[DefinitionOfCylindricalDerivative]}
$$
where $r$ here denotes the radial distance from the origin. Likewise, for $\alpha\in[0,1]$ and for $f:\Bbb{R}^2\rightarrow\Bbb{R}$, denote
$$
\delta^\alpha_\opCyl f(r) := r^\alpha\left[f|_{A(r/2,2r)}\right]_\alpha.\eqnum{\nexteqnno[DefinitionOfFractionalCylindricalDerivative]}
$$
For all non-negative integer $m$, for all $\alpha\in[0,1]$ and for all real $\delta$, define the {\bf scale free weighted H\"older norm} of any $m$-times differentiable function $f:A(R_0,\infty)\rightarrow\Bbb{R}$ by
$$
\|f\|_{C^{m,\alpha}_{\delta,\opCyl}(A(R_0,\infty))}
:=\sum_{i=0}^m\left\|r^\delta D^i_\opCyl f\right\|_{C^0(A(R_0,\infty))}
+\left\|r^\delta\delta^\alpha_\opCyl D^m_\opCyl f\right\|_{C^0(]2R_0,\infty[)}.\eqnum{\nexteqnno[DefinitionOfCylindricalNorms]}
$$
For non-neg\-ative, integer $m$, for all $\alpha\in[0,1]$, for all real $\delta$ and for any $m$-times differentiable function $f:C\rightarrow\Bbb{R}$, define
$$
\left\|f\right\|_{C^{m,\alpha}_{\delta,\opCyl}(C)} :=
\left\|f|_{C\minter(B(2R_0)\times\Bbb{R})}\right\|_{C^{m,\alpha}} + \left\|f|_{C\minter(A(R_0,\infty)\times\Bbb{R})}\right\|_{C^{m,\alpha}_{\delta,\opCyl}(A(R_0,\infty))}.
\eqnum{\nexteqnno[ScaleFreeNormOverCosta]}
$$
For all such $m$, $\alpha$ and $\delta$, let $C^{m,\alpha}_{\delta,\opCyl,g}(C)$ denote the space of $m$-times differentiable functions $f$ over $C$ which satisfy $\|f\|_{C^{m,\alpha}_{\delta,\opCyl}(C)}<\infty$ and which also satisfy $f\circ\sigma=f$ for every horizontal symmetry $\sigma$ of $C$. Observe in particular that, since each of $X$ and $Y$ has compact support, we may also think of them as taking values in $C^{0,\alpha}_{\delta+2,\opCyl,g}(C)$.
\par
Recall that, with the above symmetries imposed, for all $\delta\in]1,2[$, and for all $\alpha\in]0,1[$, the {\sl Jacobi operator} $J_C$ of $C$ defines an injective Fredholm map of Fredholm index $(-3)$ from $C^{2,\alpha}_{\delta,\opCyl,g}(C)$ into $C^{0,\alpha}_{\delta+2,\opCyl,g}(C)$ (c.f. \cite{HauswirthPacard}, \cite{Morabito}, \cite{Nayatani} and \cite{Pacard}).\footnote*{We aim to include an overview of the perturbation theory of the Costa-Hoffman-Meeks surfaces in forthcoming work, as we are not aware of any readily accessible account in the literature.}
\proclaim{Lemma \nextprocno}
\noindent For all $\alpha\in]0,1[$, for all $\delta\in]1,2[$, for all $R_0>0$ sufficiently large, and for all $\Lambda>0$ sufficiently large, the triplet $(X,Y,\hat{J}_{C})$ defines a surjective Fredholm map from $\Bbb{R}^3\oplus\Bbb{R}^3\oplus C^{2,\alpha}_{\delta,\opCyl,g}(C)$ into $C^{0,\alpha}_{2+\delta,\opCyl,g}(C)$ of Fredholm index $3$. Furthermore, the right inverse $(U,V,\Phi)$ can be chosen in such a manner that its norm is uniformly bounded, independent of $\Lambda$.
\endproclaim
\proclabel{RightInverseOverCosta}
\remark In the sequel, $R_0$ will be chosen large enough for Lemma \procref{RightInverseOverCosta} to hold for all large values of $\Lambda$. It will then be fixed once and for all, and $\Lambda$ will be made to tend to $+\infty$.
\medskip
\proof For all $c\in U$, where $U$ is a suitable open subset of $\Bbb{R}^3$, let $C_c$ be as in Section \subheadref{TheBasicSurgeryOperation} and suppose in addition that $C_c$ is also invariant under all the horizontal symmetries of $C$. Let $E:U\times C\rightarrow\Bbb{R}^3$ be a smooth function such that
\medskip
\myitem{(1)} for all $c\in U$, $E_c$ parametrises $C_c$;
\medskip
\myitem{(2)} for all $c\in U$ and for all $p\in C\minter(A(R_0,+\infty)\times\Bbb{R})$, the point $E_c(p)$ lies vertically above or below the point $p$; and
\medskip
\myitem{(3)} for all $c\in U$, $E_c:=E(c,\cdot)$ is equivariant under all the horizontal symmetries of $C$.
\medskip
\noindent Let $V$ be a neighbourhood of $0$ in $\Bbb{R}^3$ and define $\tilde{E}:U\times V\times C\rightarrow\Bbb{R}^3$ such that, for all $(c,a)\in U\times V$, and for all $p\in C$,
$$
\tilde{E}_{c,a}(p) = E_c(p) + \sum_{i=1}^3\epsilon_i\Bbb{I}_i(p)a_i(1-\chi_0'(p))e_z,
$$
where $(\epsilon_i)_{1\leq i\leq 3}$, $(\Bbb{I}_i)_{1\leq i\leq 3}$ and $\chi_0'$ are defined as in Section \subheadref{TheBasicSurgeryOperation}. Define $H:U\times V\times C\rightarrow\Bbb{R}^3$ such that, for all $(c,a)\in U\times V$, and for all $p\in C$, $H_{c,a}(p)$ is the mean curvature of the immersion $\tilde{E}_{c,a}$ at the point $p$. Define the operators $X_0,Y_0:\Bbb{R}^3\rightarrow C_0^\infty(C)$ by
$$\eqalign{
(X_0 u)(p) &:= \frac{d}{dt}H_{c_0+tu,0}(p)|_{t=0},\ \text{and}\cr
(Y_0 v)(p) &:= \frac{d}{dt}H_{c_0,tv}(p)|_{t=0}.\cr}
$$
By the perturbation theory of Costa-Hoffman-Meeks surfaces (c.f. \cite{HauswirthPacard}), $(X_0,Y_0,J_C)$ defines a surjective Fredholm map of Fredholm index $3$ from $\Bbb{R}^3\oplus\Bbb{R}^3\oplus C^{2,\alpha}_{\delta,\opCyl,g}(C)$ into $C^{0,\alpha}_{\delta+2,\opCyl,g}(C)$.
\par
Let $N$ and $\hat{N}$ be respectively the unit normal vector field and the modified normal vector field over $C$. Observe that, as $\Lambda$ and $R_0$ tend to $+\infty$, the difference $(\hat{N}-N)$ tends to zero in the $C^k$ sense for all $k$ so that the difference $(\hat{J}_{C}-J_C)$ tends to $0$ in the operator norm. Next, it is straightforward to show that, considered as an operator from $\Bbb{R}^3$ into $C^{0,\alpha}_{\delta+2,\opCyl,g}(C)$, $\left\|Y-Y_0\right\|\lesssim\epsilon$. Finally, by \eqnref{ControlOfParameters}, \eqnref{NormOfDeficiencyTermInLowerRegion} and \eqnref{NormOfDeficiencyTermOfJoinedSurface}, considered as another operator between these two spaces, $\left\|X-X_0\right\|\lesssim R^{\delta-2}$. Since these both tend to $0$ as $\Lambda$ tends to $+\infty$, the result follows by the stability of surjectivity of Fredholm maps under small perturbations.\qed
\medskip
We now review the analytic properties of rotationally symmetric Grim ends. Thus, for all non-negative, integer $m$, for all $\alpha\in[0,1]$, for all $\gamma\in\Bbb{R}$ and for all $\epsilon>0$, define the following weighted H\"older and Sobolev norms for functions over $\Bbb{R}^2$,
$$\eqalign{
\|f\|_{C^{m,\alpha}_{\gamma,\epsilon}(G)} &:= \|f(\cdot/\epsilon)\|_{C^{m,\alpha}_\gamma(G)},\ \text{and}\cr
\|f\|_{H^m_{\gamma,\epsilon}(G)} &:= \|f(\cdot/\epsilon)\|_{H^m_\gamma(G)}.\cr}\eqnum{\nexteqnno[DefinitionOfGrimNorms]}
$$
and define the {\bf hybrid norm} by
$$
\|f\|_{m,\alpha,\gamma,\epsilon} := \|f\|_{C^{m,\alpha}_{\gamma,\epsilon}(G)} + \frac{1}{\epsilon R}\|f\|_{H^m_{\gamma,\epsilon}(G)}.\eqnum{\nexteqnno[HybridNormII]}
$$
For all such $m,\alpha,\gamma$, let $\Cal{L}^{m,\alpha}_{\gamma,\epsilon,g}(G)$ denote the space of $m$-times differentiable functions with finite hybrid norm.
\par
Let $G$ be a rotationally symmetric Grim end of speed $\epsilon$ over the annulus $A(R/4,+\infty)$, and let $\hat{J}_G$ be its modified MCFS Jacobi operator, as defined in Sections \subheadref{TheModifiedMCFSJacobiOperator} and \subheadref{ModifiedJacobiOperators}. Upon rescaling, Theorem \procref{JIsLinearIsomorphismOverEnds} immediately yields
\proclaim{Lemma \nextprocno}
\noindent For all $\alpha\in]0,1[$, for all sufficiently small $\gamma$, and for sufficiently large $\Lambda$, the operator $\epsilon^2\hat{J}_G$ defines a linear isomorphism from $\Cal{L}^{2,\alpha}_{\gamma,\epsilon,g}(G)$ into $\Cal{L}^{0,\alpha}_{\gamma,\epsilon,g}(G)$. Furthermore, we may suppose that the operator norm of its inverse is uniformly bounded independent of $\Lambda$.
\endproclaim
\proclabel{RightInverseOverGrim}
\noindent We also find it useful to introduce the following notation. Define the operators $D_G$ and $\delta^\alpha_G$ by
$$\eqalignno{
D_G&:=\frac{1}{\epsilon}D,\ \text{and}&\nexteqnno[DefinitionOfGrimDerivative]\cr
\delta^\alpha_Gf(x) &:= \frac{1}{\epsilon^\alpha}\left[f|_{B(x,1/\epsilon)}\right]_\alpha.&\nexteqnno[DefinitionOfFractionalGrimDerivative]\cr}
$$
In particular, up to uniform equivalence, for any function $f$ supported in $A(R/4,2R^4)$,
$$
\|f\|_{C^{m,\alpha}_{\gamma,\epsilon}(G)} = \sum_{i=0}^m\|D_G^i f\|_{C^0} + \|\delta^\alpha_GD^m_G f\|_{C^0}.
\eqnum{\nexteqnno[OtherFormulaForGrimNorm]}$$
Likewise, let $\opdVol$ denote the canonical volume form of $\Bbb{R}^2$ and, in analogy to \eqnref{DefinitionOfCylindricalDerivative}, \eqnref{DefinitionOfFractionalCylindricalDerivative}, \eqnref{DefinitionOfGrimDerivative} and \eqnref{DefinitionOfFractionalGrimDerivative}, define
$$\eqalign{
\opdVol_\opCyl &:= \frac{1}{r^2}\opdVol,\ \text{and}\cr
\opdVol_G &:= \epsilon^2\opdVol.\cr}\eqnum{\nexteqnno[DefinitionOfCylindricalAndGrimVolumeForms]}
$$
In particular, a formula similar to \eqnref{OtherFormulaForGrimNorm} also holds for $\|f\|_{H^{m}_{\gamma,\epsilon}(G)}$ when $f$ is supported over the annulus $A(R/4,2R^4)$.
\newsubhead{Ping-Pong - Batting Up}[ThePingPongArgumentPartI]
We now describe the iteration process used to construct the Green's operator of the approximate MCF soliton. As before, for $g$ a positive integer, let $C:=C_g$ denote the CHM surface of genus $g$ and let $S:=S_g$ denote the surface obtained by replacing each of its ends with their respective joined ends, as described in Section \subheadref{TheBasicSurgeryOperation}. Since there is a natural diffeomorphism from $C$ to $S$ which maps points in the ends of $C$ vertically upwards or downwards, functions over $C$ may equally well be considered as functions over $S$ and vice-versa. As before, we will pass between these two perspectives without comment.
\par
For notational convenience, we will henceforth work as if $C$ and $S$ had only one end. Consider now the following seminorms for functions over $S$.
$$\eqalign{
\|f\|_{m,C} &:= \|f|_{B(0,4R)}\|_{C^{m,\alpha}_{(2-m)+\delta,\opCyl}(C)},\cr
\|f\|_{m,G,H} &:= \|f|_{A(R,\infty)}\|_{C^{m,\alpha}_{\gamma,\epsilon}(G)},\cr
\|f\|_{m,G,S} &:= \|f|_{A(R,\infty)}\|_{H^m_{\gamma,\epsilon}(G)},\ \text{and}\cr
\|f\|_{m,G} &:= \|f\|_{m,G,H} + \frac{1}{\epsilon R}\|f\|_{m,G,S}.\cr}\eqnum{\nexteqnno[SeminormsOfJoinedSurface]}
$$
Let $\Cal{E}$ denote the closure with respect to $\|\cdot\|_{0,C}$ of the space of functions supported over $S\minter(B(4R)\times\Bbb{R})$ which are invariant under every horizontal symmetry of the CHM surface $C$. Likewise, let $\Cal{F}$ denote the closure with respect to $\|\cdot\|_{0,G}$ of the space of functions supported over $S\minter(A(R,\infty)\times\Bbb{R})$ that are also invariant under these symmetries.
\par
We define the operator $A:\Cal{E}\rightarrow\Cal{F}$ by
$$
Ae := \hat{J}_S\chi_u\Phi e + XUe + YVe - e,\eqnum{\nexteqnno[DefinitionOfUpwardOperator]}
$$
where $\chi_u$ is the cut-off function of the upper transition region $A(R^4,2R^4)$, and $(U,V,\Phi)$ is defined as in Lemma \procref{RightInverseOverCosta}. This operator measures the extent to which $(U,V,\chi_u\Phi)$ fails to be a Green's operator of $(X,Y,\hat{J}_S)$ for functions in $\Cal{E}$. In particular, since $\hat{J}_S$ coincides with $\hat{J}_C$ over $B(0,R)$, $Ae$ is supported in the interior of $A(R,\infty)$ making it indeed an element of $\Cal{F}$. In addition, by definition of $\hat{J}_S$, and bearing in mind that $X$ and $Y$ are both supported in $B(2R)$,
$$
Ae = [\hat{J}_G,\chi_u]\Phi e + \chi_u(\hat{J}_S - \hat{J}_C)\Phi e.\eqnum{\nexteqnno[UsefulFormalaForUpwardOperator]}
$$
In this section, we prove
\proclaim{Theorem \nextprocno}
\noindent For all $\delta>1$,
$$
\|Ae\|_{0,G} \lesssim \frac{1}{(\epsilon R)^{2\alpha}}\frac{1}{R^{6+\delta}}\|e\|_{0,C}.\eqnum{\nexteqnno[NormOfA]}
$$
\endproclaim
\proclabel{NormOfA}
\noindent Theorem \procref{NormOfA} follows immediately from \eqnref{UsefulFormalaForUpwardOperator} together with \eqnref{HolderNormOfJPhiInUpperRegion}, \eqnref{SobolevNormOfJPhiInUpperRegion}, \eqnref{HolderNormOfCommutatorInUpperRegion} and \eqnref{SobolevNormOfCommutatorInUpperRegion}, below, and the fact that
$$
\|\chi_u\|_{C^{0,\alpha}_{\gamma,\epsilon}(G)}\lesssim\frac{1}{(\epsilon R^4)^\alpha}\lesssim\frac{1}{(\epsilon R)^\alpha}.
$$
For convenience, we now denote $\phi:=\Phi e$.
\noskipproclaim{Lemma \nextprocno}
$$
\left\|\left(\hat{J}_S-\hat{J}_C\right)\Phi e|_{A(R,2R^4)}\right\|_{C^{0,\alpha}_{\gamma,\epsilon}(G)} \lesssim \frac{1}{(\epsilon R)^\alpha}\frac{1}{R^{6+\delta}}\|e\|_{0,C}.\eqnum{\nexteqnno[HolderNormOfJPhiInUpperRegion]}
$$
\endproclaim
\proof Indeed, by \eqnref{DefinitionOfCylindricalDerivative}, for $k\in\left\{0,1,2\right\}$, over $A(R,2R^4)$,
$$
\left|D^k\phi\right| \lesssim \frac{1}{r^{k+\delta}}\|\phi\|_{C^{2,\alpha}_{\delta,\opCyl}(C)} \lesssim \frac{1}{r^{k+\delta}}\|e\|_{0,C}.
$$
Likewise, by \eqnref{DefinitionOfFractionalCylindricalDerivative}, for all $r\in [2R,R^4]$,
$$
\left|\delta^\alpha\left(D^2\phi|_{A(r/2,2r)}\right)\right| \lesssim \frac{1}{r^{k+\alpha+\delta}}\|e\|_{0,C}.
$$
Thus, by \eqnref{JacobiOperatorOfCHM}, \eqnref{ErrorTermsInJacobiOperatorOfCHM}, \eqnref{JacobiOperatorOfRescaledGrimEnd}, \eqnref{ErrorTermsInJacobiOperatorOfRescaledGrimEnd} and \eqnref{ErrorFromJacobiOperatorOfJoinedSurface}, over $A(R,2R^4)$,
$$\eqalign{
\left|\left(\hat{J}_S-\hat{J}_C\right)\phi\right| &\lesssim\left(\frac{\epsilon}{r^{2+\delta}} + \frac{\epsilon^2}{r^\delta} + \left[1+\opLog\left(\frac{r}{R}\right)\right]\epsilon^4r^{2-\delta}\right.\cr
&\qquad\qquad\left.+ \left[1+\opLog\left(\frac{r}{R}\right)\right]\frac{1}{r^{6+\delta}}\right)\|e\|_{0,C},\cr}\eqnum{\nexteqnno[FormulaUsedAlsoInSobolevCase]}
$$
so that, by \eqnref{MaximumOfExponential},
$$
\left|\left(\hat{J}_S-\hat{J}_C\right)\phi|_{A(R,2R^4)}\right| \lesssim \frac{1}{R^{6+\delta}}\|e\|_{0,C}.
$$
Likewise, using also \eqnref{AppFirstInterpolationInequality} and \eqnref{AppProductRule}, for $r\in [2R,R^4]$,
$$\eqalign{
\left|\delta^\alpha\left(\left(\hat{J}_S-\hat{J}_C\right)\phi|_{A(r/2,2r)}\right)\right|
&\lesssim \frac{1}{r^\alpha}\left(\frac{\epsilon}{r^{2+\delta}} + \frac{\epsilon^2}{r^\delta} + \left[1+\opLog\left(\frac{r}{R}\right)\right]\epsilon^4r^{2-\delta}\right.\cr
&\qquad\qquad\left. + \left[1+\opLog\left(\frac{r}{R}\right)\right]\frac{1}{r^{6+\delta}}\right)\|e\|_{0,C},\cr}
$$
so that, by \eqnref{DefinitionOfFractionalGrimDerivative}, for $r\in[2R,R^4]$,
$$\eqalign{
\left|\delta_G^\alpha\left(\left(\hat{J}_S-\hat{J}_C\right)\phi|_{A(r/2,2r)}\right)\right|
&\lesssim \frac{1}{(\epsilon r)^\alpha}\left(\frac{\epsilon}{r^{2+\delta}} + \frac{\epsilon^2}{r^\delta} + \left[1+\opLog\left(\frac{r}{R}\right)\right]\epsilon^4r^{2-\delta}\right.\cr
&\qquad\qquad\left. + \left[1+\opLog\left(\frac{r}{R}\right)\right]\frac{1}{r^{6+\delta}}\right)\|e\|_{0,C}.\cr}
$$
Thus, by \eqnref{AppHolderNormOverUnionII} and \eqnref{MaximumOfExponential},
$$
\left|\delta_G^\alpha\left(\left(\hat{J}_S-\hat{J}_C\right)\phi|_{A(R,2R^4)}\right)\right| \lesssim \frac{1}{(\epsilon R)^\alpha}\frac{1}{R^{6+\delta}}\|e\|_{0,C}.
$$
The result follows upon combining the above relations.\qed
\proclaim{Lemma \nextprocno}
\noindent For all $\delta>1$,
$$
\left\|\left(\hat{J}_S-\hat{J}_C\right)\Phi e|_{A(R,2R^4)}\right\|_{H^0_{\gamma,\epsilon}(G)} \lesssim \frac{(\epsilon R)}{R^{6+\delta}}\|e\|_{0,C}.\eqnum{\nexteqnno[SobolevNormOfJPhiInUpperRegion]}
$$
\endproclaim
\proof By \eqnref{DefinitionOfCylindricalAndGrimVolumeForms} and \eqnref{FormulaUsedAlsoInSobolevCase}, over $A(R,2R^4)$,
$$\eqalign{
\left|\left(\hat{J}_S-\hat{J}_C\right)\phi\right|^2\opdVol_G &\lesssim \left(\frac{\epsilon^4}{r^{2+2\delta}} + \epsilon^6r^{2-2\delta}+\left[1+\opLog\left(\frac{r}{R}\right)^2\right]\epsilon^{10}r^{6-2\delta}\right.\hfill\cr
&\qquad\qquad\left. +\left[1+\opLog\left(\frac{r}{R}\right)^2\right]\frac{\epsilon^2}{r^{10+2\delta}}\right)\|e\|_{0,C}^2\opdVol_\opCyl,\cr}
$$
so that, by \eqnref{IntegralOfExponential},
$$
\int_{A(R,2R^4)}\left|\left(\hat{J}_S-\hat{J}_C\right)\phi\right|^2\opdVol_G \lesssim \frac{(\epsilon R)^2}{R^{12+2\delta}}\|e\|_{0,C}^2,
$$
and the result follows.\qed
\proclaim{Lemma \nextprocno}
\noindent For all $\delta>1$,
$$
\|[\hat{J}_G,\chi_u]\Phi e\|_{C^{0,\alpha}_{\gamma,\epsilon}(G)} \lesssim \frac{1}{(\epsilon R^4)^\alpha}\frac{1}{R^{8+4\delta}}\|e\|_{0,C}.\eqnum{\nexteqnno[HolderNormOfCommutatorInUpperRegion]}
$$
\endproclaim
\proof By \eqnref{DefinitionOfCylindricalDerivative} and \eqnref{DefinitionOfCylindricalNorms} for $k\in\left\{0,1,2\right\}$, over $A(R^4,2R^4)$,
$$
\left|D^k\phi\right| \lesssim \frac{1}{R^{4k+4\delta}}\|\phi\|_{C^{2,\alpha}_{\delta,\opCyl}(C)} \lesssim \frac{1}{R^{4k+4\delta}}\|e\|_{0,C}.
$$
It follows by \eqnref{ErrorTermsInCommutators} that, for $k\in\left\{0,1\right\}$, over this annulus,
$$
\left|D^k[\hat{J}_G,\chi_u]\phi\right| \lesssim \frac{1}{R^{8+4k+4\delta}}\|e\|_{0,C}.\eqnum{\nexteqnno[SecondFormulaAlsoUsedInSobolevCase]}
$$
Thus, by \eqnref{DefinitionOfGrimDerivative}, for $k\in\left\{0,1\right\}$, over this annulus,
$$
\left|D^k_G[\hat{J}_G,\chi_u]\phi\right| \lesssim \frac{1}{(\epsilon R^4)^k}\frac{1}{R^{8+4\delta}}\|e\|_{0,C},
$$
and the result follows by \eqnref{AppFirstInterpolationInequality}.\qed
\noskipproclaim{Lemma \nextprocno}
$$
\|[\hat{J}_G,\chi_u]\Phi e\|_{H^0_{\gamma,\epsilon}(G)} \lesssim \frac{(\epsilon R)}{R^{5+4\delta}}\|e\|_{0,C}.\eqnum{\nexteqnno[SobolevNormOfCommutatorInUpperRegion]}
$$
\endproclaim
\proof By \eqnref{SecondFormulaAlsoUsedInSobolevCase} and \eqnref{DefinitionOfCylindricalAndGrimVolumeForms}, over $A(R^4,2R^4)$,
$$
\left|[\hat{J}_G,\chi_u]\phi\right|^2\opdVol_G \lesssim \frac{\epsilon^2}{R^{8+8\delta}}\|e\|^2_{0,C}\opdVol_\opCyl,
$$
so that, by \eqnref{IntegralOfExponential},
$$
\int_{A(R^4,2R^4)}\left|[\hat{J}_G,\chi_u]\phi\right|^2\opdVol_G \lesssim \frac{\epsilon^2}{R^{8+8\delta}}\|e\|_{0,C},
$$
and the result follows.\qed
\medskip
\noindent These estimates prove Theorem \procref{NormOfA}. In addition, the following estimate will also be of use later.
\proclaim{Lemma \nextprocno}
\noindent For all $\delta>1$,
$$
\|\chi_u\Phi e\|_{2,G} \lesssim \frac{1}{(\epsilon R)^\alpha}\frac{1}{\epsilon^2 R^{2+\delta}}\|e\|_{0,C}.\eqnum{\nexteqnno[NormOfPhiInUpperRegion]}
$$
\endproclaim
\proof Indeed, since $\chi_u=O(r^{-k})$, $\|\chi_u\|_{C^{2,\alpha}_{0,\opCyl}(C)}\lesssim 1$. Thus
$$
\|\chi_u\phi\|_{C^{2,\alpha}_{\delta,\opCyl}(C)} \lesssim \|\phi\|_{C^{2,\alpha}_{\delta,\opCyl}(C)} \lesssim \|e\|_{C^{0,\alpha}_{2+\delta,\opCyl}(C)}=\|e\|_{0,C}.
$$
Thus, by \eqnref{DefinitionOfCylindricalDerivative}, \eqnref{DefinitionOfCylindricalNorms} and \eqnref{DefinitionOfGrimDerivative}, for $k\in\left\{0,1,2\right\}$, over $A(R,2R^4)$,
$$
\left|D^k_G\chi_u\phi\right| \lesssim \frac{1}{(\epsilon r)^k}\frac{1}{r^\delta}\|e\|_{0,C}.
$$
Likewise, by \eqnref{DefinitionOfFractionalCylindricalDerivative}, \eqnref{DefinitionOfCylindricalNorms} and \eqnref{DefinitionOfFractionalGrimDerivative}, for all $r\in [2R,R^4]$,
$$
\left|\delta_G^\alpha\left(D_G^2\chi_u\phi|_{A(r/2,2r)}\right)\right| \lesssim \frac{1}{(\epsilon r)^{2+\alpha}}\frac{1}{r^\delta}\|e\|_{0,C},
$$
so that, by \eqnref{AppHolderNormOverUnionII},
$$
\left|\delta_G^\alpha \left(D_G^2\chi_u\phi|_{A(R,2R^4)}\right)\right| \lesssim \frac{1}{(\epsilon R)^{2+\alpha}}\frac{1}{R^\delta}\|e\|_{0,C}.
$$
Combining the above relations yields
$$
\|\chi_u\phi\|_{2,G,H} \lesssim \frac{1}{(\epsilon R)^\alpha}\frac{1}{\epsilon^2 R^{2+\delta}}\|e\|_{0,C}.
$$
Likewise, by \eqnref{DefinitionOfCylindricalAndGrimVolumeForms}, for all $k$, over $A(R,2R^4)$,
$$
\left|D_G^k\chi_u\phi\right|^2\opdVol_G \lesssim \frac{1}{(\epsilon r)^{2k}}\frac{1}{r^{2\delta}}\|e\|^2_{0,C}(\epsilon r)^2\opdVol_\opCyl,
$$
and since $\delta>1$, it follows by \eqnref{IntegralOfExponential} that
$$
\|\chi_u\phi\|_{2,G,S}\lesssim\frac{1}{R^\delta}\|e\|_{0,C}\lesssim \frac{1}{\epsilon R^{1+\delta}}\|e\|_{0,C}.
$$
The result follows.\qed
\newsubhead{Ping-Pong - Batting Down}[ThePingPongArgumentPartII]
By Lemma \procref{RightInverseOverGrim}, there exists a linear map $\Psi:C^{0,\alpha}_{\gamma,\epsilon,g}(G)\minter H^0_{\gamma,\epsilon,g}(G)\rightarrow C^{2,\alpha,g}_{\gamma,\epsilon,g}(G)\minter H^2_{\gamma,\epsilon,g}(G)$ such that for all $f\in\Cal{F}$,
$$
f = \hat{J}_G\Psi f,
$$
and
$$
\|\Psi f\|_{2,\alpha,\gamma,\epsilon} \lesssim \frac{1}{\epsilon^2}\|f\|_{0,\alpha,\gamma,\epsilon}.\eqnum{\nexteqnno[NormOfPsi]}
$$
Define the operators $B:\Cal{F}\rightarrow\Cal{E}$ and $W:\Cal{F}\rightarrow\Bbb{R}^3$ by,
$$\eqalign{
Bf &:= \hat{J}_S(1-\chi_l)(\Psi f - \chi_\epsilon'(Wf)) - ZWf - f,\ \text{and}\cr
Wf &:= (\Psi f)(0).\cr}\eqnum{\nexteqnno[DefinitionOfDownwardOperator]}
$$
where $\chi_l$ is the cut-off function of the lower transition region $A(R/4,R/2)$, and $\chi_\epsilon'$ is the cut-off function of the transition region $A(1/2\epsilon,1/\epsilon)$, as in Section \subheadref{TheBasicSurgeryOperation}. As before, $B$ measures the extent to which $(-W,(1-\chi_l)(\Psi - \chi_\epsilon' W))$ fails to be a Green's operator of $(Z,\hat{J}_S)$ for functions in $\Cal{F}$. In particular, by \eqnref{ValuesOfYVAndZW} together with the fact that $\hat{J}_S$ coincides with $\hat{J}_G$ over $A(2R,\infty)$, $Bf$ is supported in $B(4R)$, and is thus indeed an element of $\Cal{E}$. In addition, since $\chi_\epsilon'=1$ over $B(4R)$, over this ball, we have
$$
Bf = -[\hat{J}_C,\chi_l](\Psi f - (\Psi f)(0)) + (1-\chi_l)\left(\hat{J}_S-\hat{J}_G\right)\Psi f.\eqnum{\nexteqnno[UsefulFormulaForDownwardOperator]}
$$
In this section, we prove
\proclaim{Theorem \nextprocno}
\noindent For sufficiently small $\alpha$,
$$
\|Bf\|_{0,C} \lesssim \frac{R^2}{(\epsilon R)}\|f\|_{0,G}.\eqnum{\nexteqnno[NormOfB]}
$$
\endproclaim
\proclabel{NormOfB}
\noindent Theorem \procref{NormOfB} follows immediately from \eqnref{UsefulFormulaForDownwardOperator} together with \eqnref{NormOfJPsiInLowerRegion} and \eqnref{NormOfCommutatorInLowerRegion}, below, and the fact that
$$
\|(1-\chi_l)\|_{C^{0,\alpha}_{0,\opCyl}(C)} \lesssim 1.
$$
For convenience, we now denote $\psi:=\Psi f$.
\noskipproclaim{Lemma \nextprocno}
$$
\|Wf\| \lesssim \frac{R^2}{(\epsilon R)}\|f\|_{0,G}.\eqnum{\nexteqnno[NormOfV]}
$$
\endproclaim
\proclabel{NormOfV}
\proof Indeed, by the Sobolev embedding theorem,
$$
\|Wf\| \lesssim \|\Psi f\|_{H^2_{\gamma,\epsilon}(G)} \lesssim (\epsilon R)\|\Psi f\|_{2,\alpha,\gamma,\epsilon}.
$$
Thus, by \eqnref{NormOfPsi},
$$
\|Wf\| \lesssim \frac{R}{\epsilon}\|f\|_{0,\alpha,\gamma,\epsilon} \lesssim \frac{R^2}{(\epsilon R)}\|f\|_{0,G},
$$
as desired.\qed
\noskipproclaim{Lemma \nextprocno}
$$
\|\left(\hat{J}_S-\hat{J}_G\right)\Psi f|_{A(R/4,2R)}\|_{C^{0,\alpha}_{2+\delta,\opCyl}(C)}\lesssim \frac{1}{(\epsilon R)}\frac{1}{R^{2-\delta}}\|f\|_{0,G}.\eqnum{\nexteqnno[NormOfJPsiInLowerRegion]}
$$
\endproclaim
\proof Indeed, by \eqnref{DefinitionOfGrimDerivative}, for $k\in\left\{0,1,2\right\}$, over $A(R/4,2R)$,
$$
\left|D^k\psi\right| \lesssim \epsilon^k\|\psi\|_{C^{2,\alpha}_{\gamma,\epsilon}(G)} \lesssim \frac{1}{\epsilon^{2-k}}\|f\|_{0,G},
$$
and so, by \eqnref{JacobiOperatorOfCHM}, \eqnref{ErrorTermsInJacobiOperatorOfCHM}, \eqnref{JacobiOperatorOfRescaledGrimEnd}, \eqnref{ErrorTermsInJacobiOperatorOfRescaledGrimEnd} and \eqnref{ErrorFromJacobiOperatorOfJoinedSurface}, over $A(R/4,2R)$,
$$
\left|\left(\hat{J}_S-\hat{J}_G\right)\psi\right| \lesssim
\frac{1}{(\epsilon R)}\left(\epsilon + \epsilon^2R^2 + \frac{1}{R^4}\right)\|f\|_{0,G}
\lesssim \frac{1}{(\epsilon R)}\frac{1}{R^4}\|f\|_{0,G}.
$$
Likewise, by \eqnref{DefinitionOfFractionalGrimDerivative},
$$
\left|\delta^\alpha\left(D^2\psi|_{A(R/4,2R)}\right)\right| \lesssim \epsilon^\alpha\|f\|_{0,G}
\lesssim \frac{1}{R^\alpha}\|f\|_{0,G}.
$$
Thus, by \eqnref{DefinitionOfFractionalCylindricalDerivative}, using also \eqnref{AppFirstInterpolationInequality} and \eqnref{AppProductRule},
$$
\left|\delta^\alpha_\opCyl\left(\left(\hat{J}_S-\hat{J}_G\right)\psi|_{A(R/4,2R)}\right)\right| \lesssim \frac{1}{(\epsilon R)}\frac{1}{R^4}\|f\|_{0,G},
$$
and the result follows.\qed
\noskipproclaim{Lemma \nextprocno}
$$
\|\Psi f - (\Psi f)(0)|_{A(R/4,2R)}\|_{2,C} \lesssim \frac{R^{2+\delta}}{(\epsilon R)^{2\alpha}}\|f\|_{0,G}.\eqnum{\nexteqnno[NormOfPsiInLowerRegion]}
$$
\endproclaim
\proof Bearing in mind \eqnref{DefinitionOfFractionalGrimDerivative} and the Sobolev embedding theorem, over $A(R/4,2R)$,
$$
[\psi_0] \lesssim (\epsilon R)^{1-\alpha}\|\psi\|_{C^{0,1-\alpha}_{\gamma,\epsilon}(G)}
\lesssim (\epsilon R)^{1-\alpha}\|\psi\|_{H^2_{\gamma,\epsilon}(G)}
\lesssim (\epsilon R)^{2-\alpha}\|\psi\|_{2,G}.
$$
Consequently, by \eqnref{NormOfPsi},
$$
[\psi_0] \lesssim \frac{R^2}{(\epsilon R)^\alpha}\|f\|_{0,G}.
$$
Likewise, by \eqnref{HybridProperty} and the subsequent remark, over this annulus,
$$
\left|D_G\psi\right| \lesssim (\epsilon R)^{1-2\alpha}\|\psi\|_{2,G}\lesssim \frac{1}{(\epsilon R)^{2\alpha}}\frac{R}{\epsilon}\|f\|_{0,G}.
$$
Finally, over this annulus,
$$
\left|D_G^2\psi\right| \lesssim \|\psi\|_{C^{2,\alpha}_{\gamma,\epsilon}(G)} \lesssim \frac{1}{\epsilon^2}\|f\|_{0,G},
$$
and
$$
\left|\delta_G^\alpha\left(D^2_G\psi|_{A(R/4,2R)}\right)\right|\lesssim\|\phi\|_{C^{2,\alpha}_{\gamma,\epsilon}(G)}\lesssim\frac{1}{\epsilon^2}\|f\|_{0,G}.
$$
The result now follows by \eqnref{DefinitionOfCylindricalDerivative}, \eqnref{DefinitionOfFractionalCylindricalDerivative}, \eqnref{DefinitionOfCylindricalNorms}, \eqnref{DefinitionOfGrimDerivative}, \eqnref{DefinitionOfFractionalGrimDerivative} and \eqnref{ControlOfParameters}.\qed
\proclaim{Lemma \nextprocno}
\noindent For sufficiently small $\alpha$,
$$
\|[\hat{J}_C,\chi_l](\Psi f - (\Psi f)(0))\|_{C^{0,\alpha}_{2+\delta,\opCyl}(C)}\lesssim\frac{1}{(\epsilon R)^{2\alpha}}R^{2+\delta}\|f\|_{0,G}.\eqnum{\nexteqnno[NormOfCommutatorInLowerRegion]}
$$
\endproclaim
\proof This follows from \eqnref{ErrorTermsInCommutators} and \eqnref{NormOfPsiInLowerRegion}.\qed
\newsubhead{Ping-Pong - Iteration}[ThePingPongArgumentPartIII]
By \eqnref{NormOfA} and \eqnref{NormOfB}, for $\delta\in]1,2[$ and for sufficiently small $\alpha$, the operator norms of the products $AB$ and $BA$ satisfy
$$
\|AB\|,\|BA\|\lesssim \frac{1}{(\epsilon R)^{2\alpha}}\frac{1}{\epsilon R^{5+\delta}}\lesssim\frac{1}{\Lambda}.
$$
We therefore define $Q_E:\Cal{E}\rightarrow\Cal{E}$ and $Q_F:\Cal{F}\rightarrow\Cal{F}$ by
$$
Q_E := \sum_{m=0}^\infty(BA)^m,\ Q_F := \sum_{m=0}^\infty(AB)^m.\eqnum{\nexteqnno[DefinitionOfAdjustmentOperators]}
$$
In particular, the operator norms of both $Q_E$ and $Q_F$ are uniformly bounded for large values of $\Lambda$. We now define
$$\eqalign{
U_C e &:= UQ_Ee,\cr
U_G f &:= -UBQ_F f,\cr
V_C e &:= VQ_Ee,\cr
V_G f &:= -VBQ_Ff,\cr
W_C e &:= WAQ_Ee,\cr
W_G f &:= -WQ_Ff,\cr
P_C e &:= \chi_u\Phi Q_E e - (1-\chi_l)\left(\Psi A Q_E e - \chi_\epsilon'(WAQ_Ee)\right),\ \text{and}\cr
P_G f &:= -\chi_u\Phi B Q_F f + (1-\chi_l)\left(\Psi Q_F f - \chi_\epsilon'(WQ_Ff)\right).\cr}\eqnum{\nexteqnno[PreliminariesDefinitionOfRightInverse]}
$$
\noskipproclaim{Lemma \nextprocno}
\noindent For all $e\in\Cal{E}$ and for all $f\in\Cal{F}$,
$$\eqalign{
\hat{J}_S P_C e + XU_C e + YV_C e + ZW_C e&= e,\cr
\hat{J}_S P_G f + XU_G f + YV_G f + ZW_G f&= f.\cr}\eqnum{\nexteqnno[PreliminaryOperatorsAreRightInverses]}
$$
\endproclaim
\proof Indeed, bearing in mind \eqnref{DefinitionOfUpwardOperator} and \eqnref{DefinitionOfDownwardOperator},
$$\multiline{
\hat{J}_SP_C e + XU_C e + YV_C e + ZW_Ce\cr
\qquad\qquad\qquad\qquad=\hat{J}_S\chi_u\Phi Q_E e + X U Q_E e + YVQ_Ee\cr
\qquad\qquad\qquad\qquad\qquad\qquad - \hat{J}_S(1-\chi_l)\left(\Psi A Q_E e - \chi_\epsilon'(W A Q_E e)\right) + ZW A Q_E e\cr
\qquad\qquad\qquad\qquad=AQ_Ee + Q_Ee - BAQ_Ee - AQ_Ee\cr
\qquad\qquad\qquad\qquad=e.\cr}$$
The second relation follows in a similar manner, and this completes the proof.\qed
\medskip
Now let $\chi$ be the cut-off function of the transition region $A(2R,4R)$. Since $\chi=O(r^{-k})$, for all $f$,
$$
\|\chi f\|_{0,C}\lesssim \|f\|_{0,C},\ \|(1-\chi)f\|_{0,G}\lesssim \frac{1}{(\epsilon R)^\alpha}\|f\|_{0,G}.\eqnum{\nexteqnno[NormOfMultiplicationByChi]}
$$
Define
$$\eqalign{
\hat{U} f &:= U_C\chi f + U_G(1-\chi)f,\cr
\hat{V} f &:= V_C\chi f + V_G(1-\chi)f,\cr
\hat{W} f &:= W_C\chi f + W_G(1-\chi)f,\ \text{and}\cr
\hat{P} f &:= P_C\chi f + P_G(1-\chi)f.\cr}\eqnum{\nexteqnno[DefinitionOfRightInverse]}
$$
In particular, by \eqnref{PreliminaryOperatorsAreRightInverses},
$$
\hat{J}_S\hat{P}f + X\hat{U}f + Y\hat{V}f + Z\hat{W}f = f,\eqnum{\nexteqnno[ConstructedOperatorIsRightInverse]}
$$
so that $(\hat{U},\hat{V},\hat{W},\hat{P})$ defines a Green's operator for $(X,Y,Z,\hat{J}_S)$. We conclude this section by determining the norms of its different components. First, since the operator norms of $U$ and $V$ are uniformly bounded, by \eqnref{NormOfB}, \eqnref{PreliminariesDefinitionOfRightInverse} and \eqnref{NormOfMultiplicationByChi},
$$\eqalign{
\|\hat{U}f\| &\lesssim \|f\|_{0,C} + \frac{R^2}{(\epsilon R)^{1+\alpha}}\|f\|_{0,G},\ \text{and}\cr
\|\hat{V}f\| &\lesssim \|f\|_{0,C} + \frac{R^2}{(\epsilon R)^{1+\alpha}}\|f\|_{0,G}.\cr}
\eqnum{\nexteqnno[NormOfJoinedLAndV]}
$$
\proclaim{Theorem \nextprocno}
\noindent For sufficiently small $\alpha$, for all $\delta\in]1,2[$, and for all $f$,
$$
\|\hat{W}f\| \lesssim \|f\|_{0,C} + \frac{R^2}{(\epsilon R)^{1+\alpha}}\|f\|_{0,G}.\eqnum{\nexteqnno[NormOfW]}
$$
\endproclaim
\proclabel{NormOfW}
\proof For $e\in\Cal{E}$, by \eqnref{NormOfA} and \eqnref{NormOfV},
$$
\|W_C e\| = \|WAQ_Ee\| \lesssim \frac{R^2}{(\epsilon R)}\|AQ_E e\|_{0,G} \lesssim \frac{1}{(\epsilon R)^{2\alpha}}\frac{1}{\epsilon R^{5+\delta}}\|e\|_{0,C}\lesssim\|e\|_{0,C}.
$$
For $f\in\Cal{F}$, by \eqnref{NormOfV},
$$
\|W_G f\| = \|WQ_F f\| \lesssim \frac{R^2}{(\epsilon R)}\|f\|_{0,G}.
$$
The result now follows by \eqnref{NormOfMultiplicationByChi}.\qed
\proclaim{Theorem \nextprocno}
\noindent For sufficiently small $\alpha$, for all $\delta\in]1,2[$, and for all $f$,
$$
\|\hat{P}f\|_{2,C} \lesssim \|f\|_{0,C} + \frac{R^2}{(\epsilon R)^{1+\alpha}}\|f\|_{0,G}.\eqnum{\nexteqnno[NormOfPInLowerRegion]}
$$
\endproclaim
\proclabel{NormOfPInLowerRegion}
\proof Consider $e\in\Cal{E}$. Observe that, over $B(4R)$,
$$
P_C e = \Phi Q_E e - (1-\chi_l)(\Psi A Q_E e - \Psi A Q_E e(0)).
$$
Now,
$$
\left\|\Phi Q_E e\right\|_{2,C} \lesssim \|e\|_{0,C},
$$
and by \eqnref{NormOfA} and \eqnref{NormOfPsiInLowerRegion},
$$
\left\|(1-\chi_l)(\Psi A Q_E e - (\Psi A Q_E e)(0))\right\|_{2,C} \lesssim \frac{1}{(\epsilon R)^{4\alpha}R^4}\|e\|_{0,C}\lesssim\|e\|_{0,C},
$$
so that
$$
\|P_C e\|_{2,C} \lesssim \|e\|_{0,C}.
$$
Now consider $f\in\Cal{F}$. Over $B(4R)$,
$$
P_G f = -\Phi B Q_F f - (1-\chi_l)(\Psi Q_F f - \Psi Q_F f(0)).
$$
By \eqnref{NormOfB},
$$
\left\|\Phi B Q_Ff\right\|_{2,C} \lesssim \frac{R^2}{(\epsilon R)}\|f\|_{0,G},
$$
and, by \eqnref{NormOfPsiInLowerRegion},
$$
\left\|(1-\chi_l)(\Psi Q_F f - (\Psi Q_F f)(0))\right\|_{2,C} \lesssim \frac{R^{2+\delta}}{(\epsilon R)^{2\alpha}}\|f\|_{0,G},
$$
so that,
$$
\left\|P_G f\right\|_{2,C} \lesssim \frac{R^2}{(\epsilon R)}\|f\|_{0,G}.
$$
The result now follows by \eqnref{NormOfMultiplicationByChi} and \eqnref{DefinitionOfRightInverse}.\qed
\proclaim{Theorem \nextprocno}
\noindent For sufficiently small $\alpha$, for all $\delta\in]1,2[$, and for all $f$,
$$
\|\hat{P}f\|_{2,G} \lesssim \frac{1}{(\epsilon R)^\alpha}\frac{1}{\epsilon^2 R^{2+\delta}}\left(\|f\|_{0,C} + \frac{R^2}{(\epsilon R)^{1+\alpha}}\|f\|_{0,G}\right).\eqnum{\nexteqnno[NormOfPInUpperRegion]}
$$
\endproclaim
\proclabel{NormOfPInUpperRegion}
\proof Consider $e\in\Cal{E}$. Observe that, over $S\minter(A(R,\infty)\times\Bbb{R})$,
$$
P_C e = \chi_u \Phi Q_E e - \Psi A Q_E e + (WAQ_Ee)\chi'_\epsilon.
$$
By \eqnref{NormOfPhiInUpperRegion},
$$
\left\|\chi_u\Phi Q_E e\right\|_{2,G} \lesssim \frac{1}{(\epsilon R)^\alpha}\frac{1}{\epsilon^2 R^{2+\delta}}\|e\|_{0,C}.
$$
By \eqnref{NormOfA} and \eqnref{NormOfPsi},
$$
\left\|\Psi A Q_E e\right\|_{2,G} \lesssim \frac{1}{(\epsilon R)^{2\alpha}}\frac{1}{\epsilon^2R^{6+\delta}}\|e\|_{0,C}.
$$
By \eqnref{NormOfA} and \eqnref{NormOfV},
$$
\|WAQ_Ee\| \lesssim \frac{R^2}{(\epsilon R)}\|AQ_Ee\|_{0,G} \lesssim \frac{1}{(\epsilon R)^{2\alpha}}\frac{1}{\epsilon R^{5+\delta}}\|e\|_{0,C}.
$$
However,
$$
\|\chi'_\epsilon\|_{2,G} \lesssim \frac{1}{(\epsilon R)},
$$
and so
$$
\|(WAQ_Ee)\chi'_{\epsilon}\|_{2,G} \lesssim \frac{1}{(\epsilon R)^{2\alpha}}\frac{1}{\epsilon^2R^{6+\delta}}\|e\|_{0,C}.
$$
Combining these relations yields,
$$
\|P_Ce\|_{2,G} \lesssim \frac{1}{(\epsilon R)^\alpha}\frac{1}{\epsilon^2R^{2+\delta}}\|e\|_{0,C}.
$$
Consider now $f\in\Cal{E}$. Over $S\minter(A(R,\infty)\times\Bbb{R})$,
$$
P_G f = -\chi_u\Phi B Q_F f + \Psi Q_F f - (WQ_Ff)\chi'_\epsilon.
$$
By \eqnref{NormOfPhiInUpperRegion} and \eqnref{NormOfB},
$$
\left\|\chi_u\Phi B Q_F f\right\|_{2,G}\lesssim \frac{1}{(\epsilon R)^\alpha}\frac{1}{\epsilon^2 R^\delta(\epsilon R)}\|f\|_{0,G}.
$$
By \eqnref{NormOfPsi},
$$
\left\|\Psi Q_F f\right\|_{2,G} \lesssim \frac{1}{\epsilon^2}\|f\|_{0,G}.
$$
By \eqnref{NormOfV},
$$
\|WQ_F f\| \lesssim \frac{R^2}{(\epsilon R)}\|f\|_{0,G},
$$
so that
$$
\|(WQ_F f)\chi'_\epsilon\|_{2,G} \lesssim \frac{1}{\epsilon^2}\|f\|_{0,G}.
$$
Combining these relations yields,
$$
\left\|P_G f\right\|_{2,G} \lesssim \frac{1}{(\epsilon R)^\alpha}\frac{1}{\epsilon^2R^\delta(\epsilon R)}\|f\|_{0,G}.
$$
The result now follows by \eqnref{NormOfMultiplicationByChi}.\qed
\newhead{Existence and Embeddedness}[ExistenceAndEmbeddedness]
\newsubhead{The Schauder Fixed-Point Theorem}[TheSchauderFixedPointTheorem]
It remains only to perturb the approximate MCF solitons constructed in Chapter \headref{SurgeryAndThePerturbationFamily} into actual MCF solitons. This perturbation will be carried out using the Schauder fixed point theorem. It will first be convenient to modify slightly the norms introduced in \eqnref{SeminormsOfJoinedSurface}. We thus define,
$$\eqalign{
\|f\|'_{m,G,H} &:= \|f|_{A(2R,\infty)}\|_{C^{m,\alpha}_{\gamma,\epsilon}(G)},\cr
\|f\|'_{m,G,S} &:= \|f|_{A(2R,\infty)}\|_{H^m_{\gamma,\epsilon}(G)},\ \text{and}\cr
\|f\|'_{m,G} &:= \|f\|'_{m,G,H} + \frac{1}{(\epsilon R)}\|f\|'_{m,G,S}.\cr}\eqnum{\nexteqnno[ModifiedGrimNorms]}
$$
By \eqnref{DefinitionOfRightInverse}, this does not affect \eqnref{NormOfJoinedLAndV}, \eqnref{NormOfW}, \eqnref{NormOfPInLowerRegion} and \eqnref{NormOfPInUpperRegion}. In addition, we will also ignore the factor $\langle\hat{N}_S,N_S\rangle^{-1}$ used in the definitions \eqnref{PerturbationOperators} and \eqnref{GeneralModifiedJacobiOperator} of $(X,Y,Z,\hat{J}_S)$. Indeed, we readily show that the operator of multiplication by this function is uniformly bounded, independent of $\Lambda$, with respect to the norms $\|\cdot\|_{0,C}$ and $\|\cdot\|_{0,G}$, for which reason it also does not affect the above estimates.
\par
For all non-negative, integer $m$, for all $\alpha\in[0,1]$ and for all real $\gamma$, let $E_{m,\alpha,\gamma}$ be the space of $m$-times differentiable functions $f:S\rightarrow\Bbb{R}$ which are invariant under all horizontal symmetries of $C$ and which satisfy
$$
\|f\|_{m,C}, \|f\|'_{m,G} < \infty.
$$
Observe that $E_{m,\alpha,\gamma}$ furnished with these norms is a Frechet space. Now let $M:U\oplus V\oplus W\oplus E_{2,\alpha,\gamma}\rightarrow E_{0,\alpha,\gamma}$ be the MCFS functional about $S$, as defined in Sections \subheadref{TheBasicSurgeryOperation} and \subheadref{ModifiedJacobiOperators}. It only remains to study how $M$ varies up to second order about $S$. As before, throughout this section, we apply \eqnref{ControlOfParameters} without comment.
\noskipproclaim{Lemma \nextprocno}
$$\eqalign{
\|M(0,0,0,0)\|_{0,C} &\lesssim R^{\delta-2},\ \text{and}\cr
\|M(0,0,0,0)\|'_{0,G} &=0.\cr}\eqnum{\nexteqnno[EstimateOfMCFSOfJoinedSurface]}
$$
\endproclaim
\proof Denote $\psi:=M(0,0,0,0)$. Since $C$ is minimal, over $B(R)$,
$$
\psi = \epsilon\mu.
$$
Thus, by \eqnref{ProfileOfCatenoidalEnd} and \eqnref{AppBasicFormulaInCoordinateCharts},
$$
\|\psi|_{B(R)}\|_{C^{0,\alpha}_{2+\delta,\opCyl}(C)} \lesssim \epsilon R^{2+\delta} \lesssim R^{\delta-2}.
$$
By \eqnref{FormulaForJoinedSurface}, over $A(R,2R)$,
$$\eqalign{
H_i &= \frac{cx^i}{r^2} + \opO(R^{-(3+k)}),\ \text{and}\cr
H_{ij} &= \frac{c}{r^2}\bigg(\delta_{ij} - \frac{x^ix^j}{2r^2}\bigg) + \opO(R^{-(4+k)}).\cr}
$$
Thus, by \eqnref{AppBasicFormulaInCoordinateCharts}, over this annulus,
$$\eqalign{
\mu &= 1 + \opO(R^{-{2+k}}),\ \text{and}\cr
g^{ij} &= \delta_{ij} + \opO(R^{-{2+k}}),\cr}
$$
so that, by \eqnref{AppMCFSFunctionalForGraphs},
$$
\psi = \opO(R^{-(4+k)}).
$$
Consequently,
$$
\|\psi|_{A(R,2R)}\|_{C^{0,\alpha}_{2+\delta,\opcyl}(C)} \lesssim R^{\delta-2},
$$
and the first estimate follows upon combining these relations. Finally, by construction, $\psi$ vanishes over $A(2R,\infty)$, so that $\|\psi\|'_{0,G}=0$, and this completes the proof.\qed
\medskip
It is straightforward to show that for $\|u\|$, $\|v\|$, $\|w\|$ and $\|f\|_{2,C}$ sufficiently small, independent of $\Lambda$,
$$
\|M(u,v,w,f) - M(0,0,0,0) - \hat{J}_S f - \hat{X}u - \hat{Y}v\|_{0,C}\lesssim \|f\|_{2,C}^2 + \|u\|^2 + \|v\|^2.\eqnum{\nexteqnno[SecondOrderPerturbationOverCHM]}
$$
The corresponding estimate over rotationally symmetric Grim ends is more subtle.
\proclaim{Lemma \nextprocno}
\noindent There exists $\eta>0$ such that, for sufficiently large $\Lambda$, if $\epsilon(\epsilon R)^{1-2\alpha}\|f\|'_{2,G}<\eta$, then
$$\multiline{
\|M(u,v,w,f) - M(0,0,0,0) - \hat{J}_Sf - \hat{Z}w\|'_{0,G}\cr
\qquad\qquad\qquad\qquad\lesssim\frac{\epsilon^2}{R}\|u\|^2 + \frac{\epsilon^2}{R}\|v\|^2 + \frac{\epsilon^2}{R}\|w\|^2 + \epsilon^3(\epsilon R)^{1-2\alpha}\left(\|f\|'_{2,G}\right)^2.\cr}\eqnum{\nexteqnno[SecondOrderPerturbationOverGrim]}
$$
\endproclaim
\proclabel{QuadraticError}
\remark Before continuing, it is worth reflecting on the terms that will appear in the following proof. First, on the scale of the rotationally symmetric Grim end, the perturbation that we make is of order $\epsilon$ so that, since this perturbation is quadratic, it introduces a factor of $\epsilon^2$. Second, returning to the scale of the joined surface introduces a further factor of $\epsilon$, this explaining the factor of $\epsilon^3$ in the formulae below.
\medskip
\proof Since $M$ is a second-order quasi-linear functional, upon rescaling, we obtain, for all $u$, for all $v$, and for all $g$ with $\|\epsilon g\|'_{1,G,H}$ sufficiently small,
$$\multiline{
\|M(u,v,0,g)-M(u,v,0,0)-\hat{J}_{S,u,v}g\|'_{0,G}\cr
\qquad\qquad\qquad\qquad\lesssim \epsilon^3\|g\|'_{1,G,H}\|g\|'_{2,G}\cr
\qquad\qquad\qquad\qquad\lesssim \frac{\epsilon^2}{R}\left(\|g\|'_{1,G,H}\right)^2 + \epsilon^3(\epsilon R)\left(\|g\|'_{2,G}\right)^2.\cr}
$$
Next, for all sufficiently small $u$ and $v$, and for all $g$,
$$\eqalign{
\left\|\left(\hat{J}_{S,u,v}-\hat{J}_S\right)g\right\|'_{0,G}
&\lesssim \epsilon^3\left(\|u\| + \|v\|\right)\|g\|'_{2,G}\cr
&\lesssim \frac{\epsilon^2}{R}\|u\|^2 + \frac{\epsilon^2}{R}\|v\|^2 + \epsilon^3(\epsilon R)\left(\|g\|'_{2,G}\right)^2.\cr}
$$
Now, bearing in mind the definition of the macroscopic perturbation in the direction of $w$,
$$\multiline{
\left\|M(u,v,w,f) - M(0,0,0,0) - \hat{J}_Sf - \hat{Z}w\right\|'_{0,G}\cr
\qquad\qquad\qquad\lesssim\left\|M(u,v,w,f) - M(u,v,0,0) - \hat{J}_Sf - \hat{Z}w\right\|'_{0,G}\cr
\qquad\qquad\qquad\lesssim\left\|M(u,v,0,f+w(1-\chi_\epsilon')) - M(u,v,0,0) - \hat{J}_S(f+w(1-\chi_{\epsilon'}))\right\|'_{0,G}\cr
\qquad\qquad\qquad\lesssim\left\|M(u,v,0,f+w(1-\chi_\epsilon')) - M(u,v,0,0) - \hat{J}_{S,u,v}(f+w(1-\chi_{\epsilon'}))\right\|'_{0,G}\cr
\qquad\qquad\qquad\qquad\qquad+ \left\|\left(\hat{J}_{S,u,v}-\hat{J}_S\right)(f+w(1-\chi_{\epsilon'}))\right\|'_{0,G}\cr
\qquad\qquad\qquad\lesssim\frac{\epsilon^2}{R}\|u\|^2 + \frac{\epsilon^2}{R}\|v\|^2 + \frac{\epsilon^2}{R}\left(\|f+w(1-\chi_\epsilon')\|'_{1,G,H}\right)^2\cr
\qquad\qquad\qquad\qquad\qquad+\epsilon^3(\epsilon R)\left(\|f+w(1-\chi_\epsilon')\|_{2,G}'\right)^2.\phantom{\left(\hat{J}_S\right)}\cr}
$$
Finally,
$$\eqalign{
\left\|(1-\chi_\epsilon')|_{A(1/2\epsilon,1/\epsilon)}\right\|'_{1,G,H} &\lesssim 1,\ \text{and}\cr
\left\|(1-\chi_\epsilon')|_{A(1/2\epsilon,1/\epsilon)}\right\|'_{2,G} &\lesssim \frac{1}{(\epsilon R)},\cr}
$$
and the result now follows by Lemma \procref{HybridProperty} and the subsequent remark.\qed
\medskip
This concludes our analysis of $M$ up to second order about $S$. We are now ready to prove existence.
\proclaim{Theorem \nextprocno}
\noindent For $\gamma$ sufficiently small, for all $\delta\in]1,2[$, for $\alpha\in]0,1[$ sufficiently small, and for $\Lambda$ sufficiently large, there exist $u$, $v$, $w$ and $f$ such that
$$
M(u,v,w,f)=0.
$$
Furthermore,
$$
\|u\|,\|v\|,\|w\|,\|f\|_{2,C}\lesssim R^{\delta-2},\ \|f\|_{2,G}\lesssim\frac{1}{(\epsilon R)^\alpha\epsilon^2R^4}.\eqnum{\nexteqnno[EstimateOfError]}
$$
\endproclaim
\proclabel{ExistenceTheorem}
\proof Fix $\gamma\ll1$, $\delta\in]1,2[$ and $\alpha\in]0,1[$ small. Set $\psi_0:=M(0,0,0,0)$ and denote
$$
(u_0,v_0,w_0,f_0) := \phi_0 := -(\hat{U}\psi_0,\hat{V}\psi_0,\hat{W}\psi_0,\hat{P}\psi_0).
$$
By \eqnref{NormOfJoinedLAndV}, \eqnref{NormOfW}, \eqnref{NormOfPInLowerRegion}, \eqnref{NormOfPInUpperRegion} and \eqnref{EstimateOfMCFSOfJoinedSurface}, there exists a constant $B>0$, such that, for all large $\Lambda$,
$$
\|u_0\|,\|v_0\|,\|w_0\|,\|f_0\|_{2,C}\leq B R^{\delta-2},\ \|f_0\|'_{2,G}\leq \frac{B}{(\epsilon R)^\alpha\epsilon^2 R^4}.
$$
Define $\Omega\subseteq\Bbb{R}^3\oplus\Bbb{R}^3\oplus\Bbb{R}^3\oplus E_{2,\alpha,\gamma}$ to be the set of all quadruplets $(u,v,w,f)$ such that,
$$
\|u\|,\|v\|,\|w\|,\|f\|_{2,C}\leq 2BR^{\delta-2},\ \|f\|'_{2,G}\leq \frac{2B}{(\epsilon R)^\alpha\epsilon^2R^4}.
$$
Observe that $\Omega$ is convex and, by the Arzela-Ascoli Theorem, for all $\alpha'<\alpha$ and $\gamma'<\gamma$, $\Omega$ is a compact subset of $\Bbb{R}^3\oplus\Bbb{R}^3\oplus\Bbb{R}^3\oplus E_{2,\alpha',\gamma'}$. For $\phi:=(u,v,w,f)$ in $\Omega$, define
$$
\Phi(\phi) := \phi_0 - (\hat{U}\psi,\hat{V}\psi,\hat{W}\psi,\hat{P}\psi),
$$
where
$$
\psi := M(u,v,w,f) - M(0,0,0,0) - \hat{J}_S f - \hat{X}u - \hat{Y}v - \hat{Z}w.
$$
By \eqnref{SecondOrderPerturbationOverCHM}, \eqnref{SecondOrderPerturbationOverGrim} and \eqnref{ControlOfParameters},
$$
\|\psi\|_{0,C} \lesssim R^{2\delta-4},\ \|\psi\|'_{0,G} \lesssim \frac{1}{(\epsilon R)^{6\alpha}R^7},
$$
so that, by \eqnref{NormOfJoinedLAndV}, \eqnref{NormOfW}, \eqnref{NormOfPInLowerRegion} and \eqnref{NormOfPInUpperRegion}, for sufficiently large $\Lambda$, $\Phi$ maps $\Omega$ to itself. Furthermore, for all $\alpha'<\alpha$ and $\gamma'<\gamma$, $\Phi$ is continuous with respect to the topology of $E_{2,\alpha',\gamma'}$. It follows by the Schauder fixed point-theorem (c.f. \cite{GilbTrud}) that there exists a fixed point $\phi$ of $\Phi$ in $\Omega$. We readily verify that $M(\phi)=0$, and this completes the proof.\qed
\proclaim{Theorem \nextprocno}
\noindent Let $(u,v,w,f)$ be as in Theorem \procref{ExistenceTheorem}. For sufficiently large $\Lambda$, the surface $\tilde{E}(u,v,w,f)$ is embedded.
\endproclaim
\proclabel{EmbeddednessTheorem}
\proof We denote the joined surface by $S$, we denote the image of $\tilde{E}(u,v,w,f)$ by $S'$, and we rescale both $S$ and $S'$ by $\epsilon$. Observe that the intersection of $S$ with $A(2\epsilon R,\infty)\times\Bbb{R}$ consists of $3$ distinct rotationally symmetric Grim ends which we denote by $G_+$, $G_0$ and $G_-$ respectively. Let $u_+$, $u_0$ and $u_-$ be their respective profiles of these ends, and let $v_+$, $v_0$ and $v_-$ be the respective derivatives of these functions in the radial direction. Observe that
$$
u_+(\epsilon R)>u_0(\epsilon R)>u_-(\epsilon R),
$$
and
$$
v_+(\epsilon R)>v_0(\epsilon R)>v_-(\epsilon R).
$$
Since $v_+$, $v_0$ and $v_-$ are all solutions of the same first-order ODE, it follows that $v_+(r)>v_0(r)>v_-(r)$ for all $r$. In particular, the ends $G_+$, $G_0$ and $G_+$ are separated vertically by a distance of no less than $\eta$, where $\eta\sim\epsilon\opLog(R)$. Let $\Omega_+$, $\Omega_0$ and $\Omega_-$ denote the open sets of points lying at a vertical distance of no more than $\eta/2$ from $G_+$, $G_0$ and $G_-$ respectively. Observe, in particular, that these $3$ sets are disjoint.
\par
Now let $G_+'$, $G_0'$ and $G_-'$ be the three ends of $S'$. Over the annulus $A(\epsilon R,2\epsilon R)$, by \eqnref{EstimateOfError},
$$
\left\|\epsilon f|_{A(R,2R)}\right\|_{C^0} \lesssim \epsilon R^{-\delta}\left\|f|_{(A(R,2R)}\right\|_{2,C}
\lesssim \epsilon R^{-2},
$$
so that, over this annulus, $G_+'$ lies strictly above $G_0'$, and $G_0'$ lies strictly above $G_-'$.
However, by Lemma \procref{HybridProperty} and the subsequent remark and \eqnref{EstimateOfError} again,
$$
\|\epsilon f\|'_{1,G,H} \lesssim \frac{1}{(\epsilon R)^{3\alpha}R^3}.
$$
Bearing in mind the definition of the norm $\|\cdot\|_{1,G,H}$, it follows that for sufficiently large $\Lambda$, $G_+'$, $G_0'$ and $G_-'$ are all graphs over $A(\epsilon R,\infty)$. Furthermore, for some large $R'$, the intersections of $G_+'$, $G_0'$ and $G_-'$ with $A(R',\infty)\times\Bbb{R}$ are contained in $\Omega_+$, $\Omega_0$ and $\Omega-$ respectively. In particular, outside $B(R')\times\Bbb{R}$, $G_+'$ lies strictly above $G_0'$ and $G_0'$ lies strictly above $G_-'$. Since vertical translates of mean curvature flow solitons are also mean curvature flow solitons, it now follows by the strong maximum principle that, over the whole of $A(\epsilon R,\infty)$, $G_+'$ lies strictly above $G_0'$ and $G_0'$ lies strictly above $G_-$. This completes the proof.\qed
\goodbreak\bigskip
\global\headno=0
\inappendicestrue
\newhead{Terminology, Conventions and Standard Results}[TerminologyConventionsEtc]
\newsubhead{General Definitions}[GeneralDefinitions] Let $\Bbb{R}^2$ and $\Bbb{R}^3$ denote respectively $2$ and $3$ dimensional Euclidean space. We consider $\Bbb{R}^2$ as the $x-y$ plane in $\Bbb{R}^3$. Let $\pi:\Bbb{R}^3\rightarrow\Bbb{R}^2$ denote the canonical projection. Let $r$ denote a smooth positive function over $\Bbb{R}^2$ which is equal to the distance to the origin outside some (suitably large) compact set. We denote the composition of $r$ with $\pi$ also by $r$. Let $e_x$, $e_y$ and $e_z$ denote the vectors of the canonical basis of $\Bbb{R}^3$. Let $e_r$, $e_\theta$ denote respectively the unit radial and unit angular vector fields about the origin over $\Bbb{R}^2$ and about the $z$-axis over $\Bbb{R}^3$. Let $D$ denote the canonical differentiation operator over $\Bbb{R}^2$ and $\Bbb{R}^3$. Let $\Delta$ denote the canonical Laplacian over $\Bbb{R}^2$ (not to be confused with $\Delta^\Sigma$, defined below). Let $C(a)$ denote the circle of radius $a$ about the origin in $\Bbb{R}^2$. Let $B(a)$ denotes the closed disk of radius $a$ about the origin in $\Bbb{R}^2$. Let $A(a,b)$ denote the closed annulus of inner radius $a$ and outer radius $b$ about the origin in $\Bbb{R}^2$. Let $\chi:[0,\infty[\rightarrow\Bbb{R}$ be a non-negative, non-increasing function such that $\chi=1$ over $[0,1]$ and $\chi=0$ over $[2,\infty[$. For all $a$, define $\chi_a:\Bbb{R}^2\rightarrow\Bbb{R}$ by $\chi_a(x)=\chi(\|x\|/a)$. We call $\chi_a$ the {\bf cut-off function} of the {\bf transition region} $A(a,2a)$. Composing with $\pi$, we likewise consider $\chi_a$ as a function over $\Bbb{R}^3$.
\newsubhead{Surface Geometry}[SurfaceGeometry]
Let $\Sigma$ be an embedded surface in $\Bbb{R}^3$. Let $N_\Sigma$ denotes the unit normal vector field over $\Sigma$. Let $\pi^\Sigma$ denote the orthogonal projection onto the tangent space of $\Sigma$. Let $\nabla^\Sigma$ denote the gradient operator as well as the Levi-Civita covariant derivative of $\Sigma$. Let $\opHess^\Sigma$ denote the intrinsic Hessian operator of $\Sigma$. Let $\Delta^\Sigma$ denote the intrinsic Laplacian of $\Sigma$. Let $\opII^\Sigma$ denote the second fundamental form of $\Sigma$. Let $A_\Sigma$ denote the shape operator of $\Sigma$. Let $H_\Sigma$ denote the mean curvature of $\Sigma$ (taken to be the {\sl sum} of the principle curvatures, or the trace of the shape operator). Let  $M_\Sigma$ denote the MCFS operator of $\Sigma$ (with speed $\epsilon$). It is given by
$$
M_\Sigma := H_\Sigma + \epsilon\langle N_\Sigma,e_z\rangle.\eqnum{\nexteqnno[AppMCFSOperator]}
$$
Let $J_\Sigma$ denote the MCFS Jacobi operator (with speed $\epsilon$) of $\Sigma$. That is, $J_\Sigma$ is the linearisation of the MCFS operator of $\Sigma$. It is given by
$$
J_\Sigma f = \Delta^\Sigma f + \opTr(A_\Sigma^2)f + \epsilon\langle\nabla^G f,e_z\rangle.\eqnum{\nexteqnno[AppMCFSJacobiOperator]}
$$
Finally, we recall the following elementary relations. For any function $f$ defined over a neighbourhood of $\Sigma$,
$$\eqalign{
\nabla^\Sigma f &= Df - \langle Df,N_\Sigma\rangle N_\Sigma,\cr
\opHess^\Sigma(f) &= \opHess(f) - \langle Df, N_\Sigma\rangle \opII^\Sigma.\cr}\eqnum{\nexteqnno[AppGradientAndHessian]}
$$
Given any positive function $\phi$ defined over $\Sigma$, if $\hat{J}_\Sigma:=M_\phi^{-1}J_\Sigma M_\phi$ denotes the conjugate of $J_\Sigma$ with the operator of multiplication by $\phi$, then
$$
\hat{J}_\Sigma f = \Delta^\Sigma f + 2\phi^{-}\langle\nabla^\Sigma\phi,\nabla^\Sigma f\rangle + \epsilon \langle \nabla^\Sigma f, e_z\rangle + (\phi^{-1}J_\Sigma\phi)f.\eqnum{\nexteqnno[AppCommutator]}
$$
\newsubhead{Surface Geometry of Graphs}[SurfaceGeometryOfGraphs]
If $\Sigma$ is the graph of a function $u$ over a subset of $\Bbb{R}^2$, then we call $u$ the {\bf profile} of $\Sigma$. In this case, $\pi$ defines a coordinate chart of $\Sigma$ in $\Bbb{R}^2$. It will be more convenient to work, sometimes over $\Sigma$, and sometimes over $\Bbb{R}^2$, and we will move freely between these two perspectives. Let $g_{ij}$ denote the intrinsic metric of $\Sigma$. Its inverse is denoted by $g^{ij}$. Let $\Gamma^k_{ij}$ denote the Christoffel symbols of the Levi-Civita covariant derivative of $g_{ij}$. Denoting
$$
\mu := \langle e_z,N_\Sigma\rangle,\eqnum{\nexteqnno[ThirdComponentOfNormalVector]}
$$
\noindent we readily verify the following relations.
$$\eqalign{
\mu &= \frac{1}{\sqrt{1+\|Du\|^2}},\cr
g_{ij} &= \delta_{ij} + u_iu_j,\cr
g^{ij} &= \delta_{ij} - \mu^2u^iu^j,\cr
\Gamma^k_{ij} &= g^{kp}u_{ij}u_p,\cr
\opHess^\Sigma(f)_{ij} &= f_{ij} - g^{kp}u_{ij}u_pf_k,\cr
\Delta^\Sigma(f) &= g^{ij}f_{ij} - g^{ij}g^{kp}u_{ij}u_pf_k,\cr
\opII^\Sigma_{ij} &= -\mu u_{ij},\cr
(A^\Sigma)^i_j &= -\mu g^{ip}u_{pj},\cr
H^\Sigma &= -\mu g^{ij}u_{ij},\cr
\pi^T(e_z)_i &= \mu^2 u_i.\cr}\eqnum{\nexteqnno[AppBasicFormulaInCoordinateCharts]}
$$
Finally, when $\Sigma$ is a graph, the MCFS functional is given by,
$$
M_\Sigma = \mu g^{ij} u_{ij} + \epsilon \mu.\eqnum{\nexteqnno[AppMCFSFunctionalForGraphs]}
$$
\newsubhead{Function Spaces}[FunctionSpaces]
\noindent Let $X$ be a metric space. For all $\alpha\in[0,1]$, we define the {\bf H\"older seminorm} of order $\alpha$ over $X$ by
$$
[f]_\alpha := \msup_{x\neq y\in X}\frac{\left|f(x) - f(y)\right|}{d(x,y)^\alpha}.\eqnum{\nexteqnno[AppHolderSeminorm]}
$$
Observe that $[f]_0$ measures the {\bf total oscillation} of $f$. In particular,
$$
[f]_0 \leq 2\|f\|_{C^0}.\eqnum{\nexteqnno[AppTotalVariationBoundedByTwiceSupremum]}
$$
For all $\alpha\in[0,1]$,
$$
[f]_\alpha \leq [f]_0^{1-\alpha}[f]_1^\alpha \leq 2^{1-\alpha}\|f\|_{C^0}^{1-\alpha}[f]_1^\alpha.\eqnum{\nexteqnno[AppFirstInterpolationInequality]}
$$
If $X$ is a complete manifold, and if $f$ is differentiable over $X$, then for all $\alpha\in[0,1[$ and for all $\beta\in]0,1]$,
$$
\|Df\|_{C^0} \leq 2[f]_\alpha^{\frac{\beta}{1+(\beta-\alpha)}}[Df]_\beta^{\frac{1-\alpha}{1+(\beta-\alpha)}}.
\eqnum{\nexteqnno[AppSecondInterpolationInequality]}
$$
For all $\alpha$,
$$
[fg]_\alpha \leq \|f\|_{C^0}[g]_\alpha + [f]_\alpha\|g\|_{C^0}.\eqnum{\nexteqnno[AppProductRule]}
$$
Finally, if $X=X_1\munion...\munion X_m$, then, for all $\alpha$,
$$
[f]_\alpha \leq m^{1-\alpha}\msup_{1\leq k\leq m}[f|_{X_i}]_\alpha.\eqnum{\nexteqnno[AppHolderNormOverUnionI]}
$$
If, in particular, $X=[0,m+1]\times S^1$ is a cylinder and $X_i=[i,i+1]\times S^1$ for all $i$, then \eqnref{AppHolderNormOverUnionI} refines to
$$
[f]_\alpha \leq \sum_{i=1}^m [f|_{X_i}]_\alpha.\eqnum{\nexteqnno[AppHolderNormOverUnionII]}
$$
For a continuous function $f$ over $X$, for all $\alpha$, we define
$$
\delta^\alpha f(x) := [f|_{B_1(x)}]_\alpha.\eqnum{\nexteqnno[AppDefinitionOfFractionalDifferentialOperator]}
$$
Now suppose that $X$ is a smooth Riemannian manifold. For all $k,\alpha$, we define the $C^{k,\alpha}$-H\"older norm over $C^\infty(M)$ by
$$
\|f\|_{C^{k,\alpha}} := \sum_{i=0}^k\|D^i f\|_{C^0} + \|\delta^\alpha D^k f\|_{C^0}.\eqnum{\nexteqnno[AppDefinitionOfHolderNorm]}
$$
We define the space $C^{k,\alpha}(X)$ to be the closure of $C^\infty(X)$ with respect to this norm. For all $p$, we define the $L^p$-norm over $C_0^\infty(M)$ by
$$
\|f\|_{L^p}^p := \int_X \left|f\right|^p\opdVol.\eqnum{\nexteqnno[AppDefinitionOfLPNorm]}
$$
We define the space $L^p(X)$ to be the closure of $C_0^\infty(X)$ with respect to this norm. For all $k$, we define the $H^k$-Sobolev norm over $C_0^\infty(M)$ by
$$
\|f\|_{H^k} := \sum_{i=0}^k\|D^if\|_{L^2}.\eqnum{\nexteqnno[AppDefinitionOfSobolevNorm]}
$$
The reader may verify that all surfaces studied in this paper are sufficiently regular at infinity for the Sobolev embedding theorem to hold. That is for all $l$, and for all $k+\alpha<l-1$,
$$
\|f\|_{C^{k,\alpha}} \lesssim \|f\|_{H^l}.\eqnum{\nexteqnno[AppSobolevEmbeddingTheorem]}
$$
The following formulae are readily verified.
$$
\msup_{t\in[1,T]}\opLog(t)t^\alpha \lesssim \left\{\matrix \opLog(T)T^\alpha\hfill&\text{if}\ \alpha>0,\cr
1\hfill&\text{if}\ \alpha<0,\cr\endmatrix\right.\eqnum{\nexteqnno[MaximumOfExponential]}
$$
and
$$
\int_{A(1,T)}\opLog(r)^m r^\alpha \opdVol_{\opCyl} \lesssim \left\{\matrix \opLog(T)^mT^\alpha\hfill&\text{if}\ \alpha>0,\cr
1\hfill&\text{if}\ \alpha<0.\cr\endmatrix\right.\eqnum{\nexteqnno[IntegralOfExponential]}
$$
\newsubhead{Elliptic Estimates}[EllipticEstimates]
\noindent Let $E$ and $F$ be Banach spaces and let $A:E\rightarrow F$ be a bounded linear map. We say that $A$ satisfies an {\bf elliptic estimate} whenever there exists a normed vector space $G$, a compact map $K:E\rightarrow G$, and a constant $C$ such that for all $e$ in $E$,
$$
\|e\| \leq C\left(\|Ke\| + \|Ae\|\right).\eqnum{\nexteqnno[AppEllipticEstimate]}
$$
The following straightforward result plays an important r\^ole in Fredholm theory.
\proclaim{Theorem \nextprocno}
\noindent If $A$ satisfies an elliptic estimate, then the kernel of $A$ is finite-dimensional and its image is a closed subset of $F$.
\endproclaim
\proclabel{ThmEllipticToFredholm}
\goodbreak
\newhead{Bibliography}[Bibliography]
{\leftskip = 5ex \parindent = -5ex
\leavevmode\hbox to 4ex{\hfil \cite{AltschulerWu}}\hskip 1ex{Altschuler S., Wu L. F., Translating surfaces of the non-parametric mean curvature flow with prescribed contact angle, {\sl Calc. Var. Partial Differential Equations}, {\bf 2}, (1994), 101--111}
\medskip
\leavevmode\hbox to 4ex{\hfil \cite{ClutterbuckSchnurerSchulze}}\hskip 1ex{Clutterbuck J., Schn\"urer O., Schulze F., Stability of translating solutions
to mean curvature flow, {\sl Calc. Var. Partial Differential Equations}, {\bf 29}, (2007), 281--293}
\medskip
\leavevmode\hbox to 4ex{\hfil \cite{ChengZhou}}\hskip 1ex{Cheng X., Zhou D., Stability properties and gap theorem for complete $f$-minimal hypersurfaces, arXiv 1307.5099}
\medskip
\leavevmode\hbox to 4ex{\hfil \cite{ChengMejiaZhouI}}\hskip 1ex{Cheng X., Mejia T., Zhou D., Simons' type equation for $f$-minimal hypersurfaces and applications, arXiv:1305.2379}
\medskip
\leavevmode\hbox to 4ex{\hfil \cite{ChengMejiaZhouII}}\hskip 1ex{Cheng X., Mejia T., Zhou D., Eigenvalue estimate and compactness for closed $f$-minimal surfaces, arXiv:1210.8448}
\medskip
\leavevmode\hbox to 4ex{\hfil \cite{ChengMejiaZhouIII}}\hskip 1ex{Cheng X., Mejia T., Zhou D., Stability and compactness for complete $f$-minimal surfaces, to appear in {\sl Transactions of Amer. Math. Soc.}, arXiv:1210.8076}
\medskip
\leavevmode\hbox to 4ex{\hfil \cite{DavilaDelPinoNguyen}}\hskip 1ex{D\'avila J., Nguyen X. H., Del Pino M., Finite topology self-translating surfaces for the mean curvature flow in $\Bbb{R}^3$, {\sl Adv. Math.}, {\bf 320}, 674--729}
\medskip
\leavevmode\hbox to 4ex{\hfil \cite{GilbTrud}}\hskip 1ex{Gilbarg D., Trudinger N. S., {\sl Elliptic partial differential equations of second order}, Classics in Mathematics, Springer-Verlag, Berlin, 2001}
\medskip
\leavevmode\hbox to 4ex{\hfil \cite{Halldorsson}}\hskip 1ex{Halldorsson H. P., Helicoidal surfaces rotating/translating under the mean curvature flow, {\sl Geom. Dedicata}, {\bf 162}, (2013), 45--65}
\medskip
\leavevmode\hbox to 4ex{\hfil \cite{HauswirthPacard}}\hskip 1ex{Hauswirth L., Pacard F., Minimal surfaces of finite genus with two limit ends, {\sl Invent. Math.}, {\bf 169}, no. 3, (2007), 569--620}
\medskip
\leavevmode\hbox to 4ex{\hfil \cite{HoffmanMeeks}}\hskip 1ex{Hoffman D., Meeks W.H., Embedded minimal surfaces of finite topology, {\sl Ann. of Math.}, {\bf 131}, (1990), 1--34}
\medskip
\leavevmode\hbox to 4ex{\hfil \cite{KapouleasI}}\hskip 1ex{Kapouleas N., Complete constant mean curvature surfaces in Euclidean three space, {\sl Ann. of Math.}, {\bf 121}, (1990), 239--330}
\medskip
\leavevmode\hbox to 4ex{\hfil \cite{KapouleasII}}\hskip 1ex{Kapouleas N., Constant mean curvature surfaces constructed by fusing Wente tori, {\sl Invent. Math.}, {\bf 119}, (1995), 443--518}
\medskip
\leavevmode\hbox to 4ex{\hfil \cite{KapouleasIII}}\hskip 1ex{Kapouleas N., Complete embedded minimal surfaces of finite total curvature, {\sl J. Diff Geom.}, {\bf 47}, (1997), 95--169}
\medskip
\leavevmode\hbox to 4ex{\hfil \cite{MartinSavasSmoczyk}}\hskip 1ex{Mart\'\i n F., Savas-Halilaj A., Smoczyk K., On the topology of translating solitons of the mean curvature flow, arXiv:1404.6703}
\medskip
\leavevmode\hbox to 4ex{\hfil \cite{MazzeoPacard}}\hskip 1ex{Mazzeo R., Pacard F., Constant mean curvature surfaces with Delaunay ends, {\sl Comm. Anal. Geom.}, {\bf 9}, no. 1, (2001), 169--237}
\medskip
\leavevmode\hbox to 4ex{\hfil \cite{MazzeoPollack}}\hskip 1ex{Mazzeo R., Pollack D., Gluing and moduli for noncompact geometric problems, in {\sl Geometric theory of singular phenomena in partial differential equations}, Symposia Mathematica Vol. XXXII, Cambridge Univ. Press, (1998), pp. 17--51}
\medskip
\leavevmode\hbox to 4ex{\hfil \cite{Morabito}}\hskip 1ex{Morabito F., Index and nullity of the Gauss map of the Costa-Hoffman-Meeks surfaces, {\sl  Indiana university mathematics journal}, {\bf 58}, no. 2, (2009), 677--708}
\medskip
\leavevmode\hbox to 4ex{\hfil \cite{Nayatani}}\hskip 1ex{Nayatani S., Morse Index and Gauss maps of complete minimal surfaces in Euclidean 3-space, {\sl Comment. Math. Helv.}, {\bf 68}, no. 4, (1993), 511--537}
\medskip
\leavevmode\hbox to 4ex{\hfil \cite{NguyenI}}\hskip 1ex{Nguyen X. H., Complete embedded self-translating surfaces under mean curvature flow, {\sl J. Geom. Anal.}, {\bf 23}, (2013), 1379--1426}
\medskip
\leavevmode\hbox to 4ex{\hfil \cite{NguyenII}}\hskip 1ex{Nguyen X. H., Translating tridents, {\sl Comm. Partial Differential Equations}, {\bf 34}, (2009), 257--280.}
\medskip
\leavevmode\hbox to 4ex{\hfil \cite{Pacard}}\hskip 1ex{Pacard F., {\sl Connected sum constructions in geometry and nonlinear analysis}, lecture notes available online}
\medskip
\leavevmode\hbox to 4ex{\hfil \cite{Schwarz}}\hskip 1ex{Schwarz M., {\sl Morse homology}, Progress in Mathematics, {\bf 111}, Birkh\"auser Verlag, Basel, (1993)}
\medskip
\leavevmode\hbox to 4ex{\hfil \cite{Wang}}\hskip 1ex{Wang X. J., Convex solutions to the mean curvature flow, {\sl Ann. of Math.}, {\bf 173}, (2011), 1185--1239}
\medskip
\leavevmode\hbox to 4ex{\hfil \cite{Weber}}\hskip 1ex{Weber M., Classical Minimal Surfaces in Euclidean Space by Examples: geometric and computational aspects of the weierstrass representation, in {\sl Global Theory of Minimal Surfaces} (David Hoffman Ed.), Clay Mathematics Proceedings, Vol. 2, (2004)}
\par}
%
%
%
%
\enddocument